\documentclass{article}
\usepackage{a4}
\usepackage{latexsym,amsfonts,amssymb,amsbsy}
\usepackage{array}
\usepackage{lscape}
\usepackage{citesort}
\usepackage{longtable}
\newcommand{\swap}[2]{\let\abc#1\let#1#2\let#2\abc}
\newcommand{\tcal}[1]{\ensuremath{\mathcal{#1}}}
\newcommand{\A}{\tcal{A}}
\newcommand{\C}{\tcal{C}}
\newcommand{\D}{\tcal{D}}
\newcommand{\E}{\ensuremath{\mathbf{E}}}
\newcommand{\F}{\tcal{F}}
\renewcommand{\H}{\tcal{H}}
\newcommand{\I}{\tcal{I}}
\newcommand{\J}{\tcal{J}}
\newcommand{\PJ}[1]{\ensuremath{{}_{#1}{\mathbf J}}}
\newcommand{\pJ}{\PJ{p}}
\newcommand{\RingelpJ}{\ensuremath{{}_p\Ringel{\mathbf J}}}
\newcommand{\JD}{\tcal{JD}}
\newcommand{\K}{\tcal{K}}
\renewcommand{\L}{\tcal{L}}
\newcommand{\M}{\tcal{M}}
\renewcommand{\P}{\tcal{P}}
\renewcommand{\S}{\tcal{S}}
\newcommand{\T}{\tcal{T}}
\newcommand{\W}{\tcal{W}}
\newcommand{\AlsoDoctoraldissertation}{Also Doctoral Dissertation}

\renewcommand{\arraystretch}{\myarraystretch}
{\renewcommand{\arraystretch}{#1}}%
{\renewcommand{\arraystretch}{\myarraystretch}}%
\DeclareMathAlphabet{\mathsc}{OT1}{cmr}{m}{sc}
\newcommand{\Ent}[1]{[\![#1]\!]}
\providecommand{\Cset}{\mathbb{C}}
\providecommand{\Fset}{\mathbb{F}}
\providecommand{\Kset}{\mathbb{K}}
\providecommand{\Nset}{\mathbb{N}}
\providecommand{\Oset}{\mathbb{O}}
\providecommand{\Pset}{\mathbb{P}}
\providecommand{\Rset}{\mathbb{R}}
\providecommand{\Sset}{\mathbb{S}}
\providecommand{\Yset}{\mathbb{Y}}
\providecommand{\Zset}{\mathbb{Z}}
\newcommand{\Seq}[1]{\{\!\!\{#1\}\!\!\}}
\newcommand{\SEQ}[1]{\left\{\!\!\!\left\{#1\right\}\!\!\!\right\}}
\newcommand{\Ele}[1]{\langle#1\rangle}
\newcommand{\Set}[1]{\{#1\}}

\newcommand{\Ringel}[1]{\mathop{\kern0pt{#1}}\limits^{\circ}\!}
\newcommand{\abs}[1]{\left\vert #1\right\vert}
\providecommand{\binom}[2]{{#1 \choose #2}}
\renewcommand{\d}{\;d\;}
\newcommand{\di}[2]{{}_{#1}^{({#2})}}
\newcommand{\forwarddiff}{\triangle}

\newtheorem{theorem}{Theorem}
\newtheorem{corollar}{Corollary}
\newtheorem{lemma}{Lemma}
\makeatletter
\newdimen\@bls                              
\@bls=\baselineskip                         
\newcommand{\qed}{\ensuremath{\Box}}
\newenvironment{proof}%
  {\par\addvspace{\@bls \@plus 0.5\@bls \@minus 0.1\@bls}\noindent
   {\bfseries\proofname}\enspace\ignorespaces}%
  {\qed\par\addvspace{\@bls \@plus 0.5\@bls \@minus 0.1\@bls}}
\def\proofname{PROOF.}
\makeatother

\newcommand{\float}[1]{\,\mbox{\tt fl}\left(#1\right)}
\newlength{\mylength}
\newlength{\mylengthA}
\begin{document}
\title{SCALAR LEVIN-TYPE SEQUENCE TRANSFORMATIONS%
\thanks{Invited article for J. Comput. Appl. Math.}}
\author{Herbert H.\ H.\
Homeier\thanks{E-mail: Herbert.Homeier@na-net.ornl.gov}%
\ \thanks{WWW: http://www.chemie.uni-regensburg.de/$\sim$hoh05008}\\
Institut f\"ur Physikalische und Theoretische Chemie\\
Universit\"at Regensburg, D-93040 Regensburg, Germany}

\maketitle

\begin{abstract}
  Sequence transformations are important tools for the convergence
  acceleration of slowly convergent scalar sequences or series and for
  the summation of divergent series. The basic idea is to construct
  from a given sequence $\Seq{s_n}$ a new sequence
  $\Seq{s_n'}=\T(\Seq{s_n})$ where each $s_n'$ depends on a finite
  number of elements $s_{n_1},\dots,s_{n_m}$. Often, the $s_n$ are the
  partial sums of an infinite series. The aim is to find
  transformations such that $\Seq{s_n'}$ converges faster than (or sums)
  $\Seq{s_n}$.

  Transformations $\T(\Seq{s_n},\Seq{\omega_n})$ that depend not only on
  the sequence elements or partial sums $s_n$ but also on an auxiliary
  sequence of so-called remainder estimates $\omega_n$ are of
  Levin-type if they are linear in the $s_n$, and nonlinear in the
  $\omega_n$. Such remainder estimates provide an easy-to-use
  possibility to use asymptotic information on the problem sequence for
  the construction of highly efficient sequence transformations. As
  shown first by Levin, it is possible to obtain such asymptotic
  information easily for large classes of sequences in such a way that
  the $\omega_n$ are simple functions of a few sequence elements $s_n$.
  Then, nonlinear sequence transformations are obtained. Special cases
  of such Levin-type transformations belong to the most powerful
  currently known extrapolation methods for scalar sequences and series.

  Here, we review known Levin-type sequence transformations and
  put them in a common theoretical framework. It is discussed
  how such transformations may be constructed by either a model
  sequence approach or by iteration of simple transformations. As
  illustration, two new sequence transformations are derived. Common
  properties and results on convergence acceleration and stability are
  given. For important special cases, extensions of the general results
  are presented. Also, guidelines for the application of Levin-type
  sequence transformations are discussed, and a few numerical examples are
  given.
\end{abstract}

\newpage
\tableofcontents
\listoftables
\newpage
{\quad}\section{Introduction}

In applied mathematics and the numerate sciences, extrapolation methods
are often used for the convergence acceleration of slowly convergent
sequences or series and for the summation of divergent series. For an
introduction to such methods, and also further information that cannot
be covered here, see the books of Brezinski and Redivo Zaglia
\cite{BrezinskiRedivoZaglia91} and Wimp \cite{Wimp81} and the work of
Weniger \cite{Weniger89,Weniger94} and Homeier
\cite{Homeier96Hab}, but also the books of Baker \cite{Baker75}, Baker
and Graves-Morris \cite{BakerGravesMorris96}, Brezinski
\cite{Brezinski77,Brezinski78,Brezinski80a,Brezinski91b,Brezinski91c},
Graves-Morris \cite{GravesMorris72,GravesMorris73}, Graves-Morris, Saff
and Varga \cite{GravesMorrisSaffVarga84}, Khovanskii
\cite{Khovanskii63}, Lorentzen and Waadeland
\cite{LorentzenWaadeland92}, Nikishin and Sorokin
\cite{NikishinSorokin91}, Petrushev and Popov \cite{PetrushevPopov87},
Ross \cite{Ross87}, Saff and Varga \cite{SaffVarga77}, Wall
\cite{Wall73}, Werner and Buenger \cite{WernerBuenger84} and Wuytack
\cite{Wuytack79a}.

For the discussion of extrapolation methods, one considers a sequence
$\Seq{s_n}=\Seq{s_0,s_1,\dots}$ with  elements $s_n$ or the terms
$a_n=s_{n}-s_{n-1}$ of a series $\sum_{j=0}^{\infty} a_j$ with
partial sums $s_n=\sum_{j=0}^{n} a_j$ for large $n$.
A common approach is to
rewrite $s_n$ as
  \begin{equation} \label{eqbasic}
    s_n = s + R_n
  \end{equation}
where $s$ is the limit
(or antilimit in the case of divergence) and $R_n$ is the remainder or
tail. The aim then is to find a new sequence $\Seq{s_n'}$ such that
  \begin{equation} \label{eqtransformed}
    s_n' = s + R_n'\>, \qquad R_n'/R_n \to 0 \mbox{\ for }
    n\to\infty\>.
  \end{equation}
  Thus, the sequence $\Seq{s_n'}$ converges faster to the limit $s$ (or
  diverges less violently) than $\Seq{s_n}$.

To find the sequence $\Seq{s_n'}$, i.e., to construct a sequence
transformation $\Seq{s_n'}=\T(\Seq{s_n})$, one needs asymptotic
information about the $s_n$ or the terms $a_n$ for large $n$, and
hence about the $R_n$. This information then allows to eliminate the
remainder at least asymptotically, for instance by substracting the
dominant part of the remainder. Either
such information is obtained by a careful mathematical analysis of the
behavior of the $s_n$ and/or $a_n$, or it has to be extracted
numerically from the values of a finite number of the $s_n$ and/or
$a_n$ by some method that ideally can be proven to work for a large
class of problems.

Suppose that one knows quantities $\omega_n$ such that $R_n/\omega_n =
O(1)$ for $n\to \infty$, for instance
\begin{equation}\label{eqomegan}
\lim_{n\to\infty}
R_n/\omega_n=c\ne 0
\end{equation}
where $c$ is a constant. Such quantities are
called remainder estimates. Quite often, such remainder estimates can
be found with relatively low effort but the exact value of $c$ is often
quite hard to calculate. Then, it is rather natural to rewrite the rest
as $R_n=\omega_n \mu_n$ where $\mu_n\to c$. The problem is how to
describe or model the $\mu_n$. Suppose that one has a system of
known functions $\psi_j(n)$ such that $\psi_0(n)=1$ and
$\psi_{j+1}=o(\psi_j(n))$ for  $j\in \Nset_0$. An example of such a
system is $\psi_j(n)=(n+\beta)^{-j}$ for some $\beta\in\Rset_{+}$. Then, one may model
$\mu_n$ as a linear combination of the $\psi_j(n)$ according to
  \begin{equation} \label{eqmun}
    \mu_n \sim \sum_{j=0}^{\infty} c_j \psi_j(n)\>\qquad \mbox{for } n\to\infty\>,
  \end{equation}
whence the problem sequence is modelled according to
    \begin{equation} \label{eqpsijn}
      s_n \sim s + \omega_n \sum_{j=0}^{\infty} c_j \psi_{j}(n)\>.
    \end{equation}
The idea now is to eliminate the leading terms of the remainder with
the unknown constants
$c_j$ up to
$j=k-1$, say. Thus, one uses a model sequence with elements
    \begin{equation} \label{eqpsijnmod}
      \sigma_m = \sigma + \omega_m \sum_{j=0}^{k-1} c_j
      \psi_{j}(m)\>,\qquad m\in \Nset_0
    \end{equation}
and calculates $\sigma$ exactly by solving the system of $k+1$
equations resulting for
$m=n,n+1,\dots,n+k$
for
the unknowns $\sigma$ and $c_j$, $j=0,\dots,k-1$. The solution for
$\sigma$ is a ratio of determinants (see below) and may be denoted symbolically as
\begin{equation}
\sigma=T(\sigma_n,\dots,\sigma_{n+k};\omega_{n},\dots,\omega_{n+k};
\psi_j(n),\dots,\psi_j(n+k))\>.
\end{equation}
The resulting sequence transformation is
\begin{equation}
\T(\Seq{s_n},\Seq{\omega_n})=\Seq{\T\di{n}{k}(\Seq{s_n},\Seq{\omega_n})}
\end{equation}
with
\begin{equation}\label{eqT0}
\T\di{n}{k}(\Seq{s_n},\Seq{\omega_n})=
T(s_n,\dots,s_{n+k};\omega_{n},\dots,\omega_{n+k};
\psi_j(n),\dots,\psi_j(n+k))\>.
\end{equation}
It eliminates the leading terms of the
asymptotic expansion (\ref{eqpsijn}). The model sequences
(\ref{eqpsijnmod}) are in the kernel of the sequence transformation \T,
defined as the set of all sequences such that \T\ reproduces their
(anti)limit exactly.

A somewhat more general approach is based on model sequences of the
form
\begin{equation}\label{eqbekaverf_Emodel}
\sigma_n = \sigma +  \sum_{j=1}^{k} c_j g_j(n), \qquad n\in\Nset_0, \>
k\in\Nset\>.
\end{equation}
Virtually all known sequence transformations can be derived using such
model sequences. This leads to the \E\ algorithm as described below in
Section \ref{secEalgo}. Also, some further important examples of sequence
transformations are described in Section \ref{secBasSeq}.

However, the introduction of remainder estimates proved to be an
important theoretical step since it allows to make use of asymptotic
information of the remainder easily. The most prominent of the
resulting sequence transformations $\T(\Seq{s_n},\Seq{\omega_n})$ is
the Levin transformation \cite{Levin73} that corresponds to the
asymptotic system of functions given by $\psi_{j}(n)=(n+\beta)^{-j}$,
and thus, to Poincare-type expansions of the $\mu_n$. But also other
systems are of importance, like $\psi_{j}(n)=1/(n+\beta)_{j}$ leading
to factorial series, or $\psi_{j}(n)=t_n^j$ corresponding to Taylor
expansions of $t$--dependent functions at the abscissae $t_n$ that tend
to zero for large $n$. The question which asymptotic system is best,
cannot be decided generally. The answer to this question depends on the
extrapolation problem. To obtain efficient extrapolation procedures for
large classes of problems requires to use various asymptotic systems,
and thus, a larger number of different sequence transformations. Also,
different choices of $\omega_n$ lead to different variants of such
transformations. Levin \cite{Levin73} has pioneered this question and
introduced three variants that are both simple and rather successful
for large classes of problems. These variants and some further ones
will be discussed. The question which variant is best, also cannot be
decided generally. There are, however, a number of results that favor
certain variants for certain problems. For example, for Stieltjes
series, the choice $\omega_n=a_{n+1}$ can be theoretically justified
(see Appendix \ref{appstieltjes}.).

Thus, we will focus on sequence transformations that involve an
auxiliary sequence $\Seq{\omega_n}$. To be more specific, we consider
transformations of the form
$\T(\Seq{s_n},\Seq{\omega_n})=\Seq{\T\di{n}{k}}$ with
      \begin{equation} \label{eqT}
          \T\di{n}{k} = \frac{\displaystyle\sum_{j=0}^{k} \lambda\di{n,j}{k}
                       \frac{s_{n+j}}{\omega_{n+j}}}
                      {\displaystyle\sum_{j=0}^{k} \lambda\di{n,j}{k}
                       \frac{1}{\omega_{n+j}}}
      \end{equation}
This will be called a Levin-type transformation. The known sequence
transformations that involve remainder estimates, for instance the \C,
\S, and \M\ transformation of Weniger \cite{Weniger89}, the $W$
algorithm of Sidi \cite{Sidi82}, and the \J\ transformation of Homeier
with its many special cases like the important \pJ\ transformations
\cite{Homeier93,Homeier94ahc,Homeier95,Homeier96aan,Homeier96Hab,Homeier98ots},
are all of this type. Interestingly, also the \H, \I, and \K\
transformations of Homeier
\cite{Homeier92,Homeier93,Homeier94nca,Homeier96Hab,%
TC-NA-97-1,TC-NA-97-3,TC-NA-97-4,Homeier98aah} for the extrapolation of
orthogonal expansions are of this type although the $\omega_n$ in some
sense cease to be remainder estimates as defined in Eq.
(\ref{eqomegan}).

The Levin transformation was also generalized in a different way by Levin and Sidi 
 \cite{LevinSidi81} who introduced the
$d^{(m)}$ transformations. This is an important class of transformations that would deserve
a thorough review itself.  This, however, is outside the scope of the present review.
We collect some important facts regarding this class of transformations in Section \ref{appdm}. 

  Levin-type transformations as defined in Eq.\ (\ref{eqT})  have been used for the solution of a
  large variety of problems. For instance, Levin-type sequence
  transformations have be applied for the convergence
  acceleration of infinite series representations of molecular
  integrals
  \cite{GrotendorstSteinborn86,GrotendorstWenigerSteinborn86,%
  WenigerSteinborn87,%
  WenigerSteinborn88,WenigerSteinborn89a,SteinbornWeniger90,%
  Homeier90,PelzlKing98}, for the calculation of the lineshape of spectral
  holes \cite{TC-PC-95-1}, for the extrapolation of cluster- and
  crystal-orbital calculations of one-dimensional polymer chains
  to infinite chain length
  \cite{CioslowskiWeniger93,WenigerLiegener90,Weniger94}, for
  the calculation of special functions
  \cite{GrotendorstSteinborn86,WenigerSteinborn89a,%
  SteinbornWeniger90,Weniger94,Homeier96Hab,Weniger96cotw,%
  WenigerCizek90}, for the summation of divergent and
  acceleration of convergent quantum mechanical perturbation
  series
  \cite{Grotendorst91,Weniger90,WenigerCizekVinette91,CizekVinetteWeniger91,%
  CizekVinetteWeniger93,WenigerCizekVinette93,Weniger96acr,%
  Weniger96cots,Weniger96nst,Weniger96nste}, for the evaluation
  of semiinfinite integrals with oscillating integrands and
  Sommerfeld integral tails
  \cite{Sidi88c,SteinbornHomeierFernandezRicoEmaLopezRamirez99,Michalski98,Mosig89}, and for the
  convergence acceleration of multipolar and orthogonal
  expansions and Fourier series
  \cite{SmithFord82,Homeier92,Homeier93,Homeier94nca,%
  Homeier96Hab,TC-NA-97-1,TC-NA-97-3,%
  TC-NA-97-4,Homeier98oca,Homeier98aah,Oleksy96,Sidi95}.
  This list is clearly not complete but sufficient to
  demonstrate the possibility of successful application of these
  transformations.

  The outline of this survey is as follows: After listing some
  definitions and notations, we discuss some basic sequence
  transformations in order to provide some background
  information. Then, special definitions relevant for Levin-type
  sequence transformations are given, including variants
  obtained by choosing specific remainder estimates $\omega_n$.
  After this, important examples of Levin-type sequence
  transformations are introduced. In Section \ref{secMetCon}, we
  will discuss approaches for the construction of Levin-type
  sequence transformations, including model sequences, kernels and
  annihilation operators, and also  the concept of hierarchical
  consistency. In Section \ref{secProper}, we derive basic
  properties, those of limiting transformations and discuss the
  application to power series. In Section \ref{secConRes},
  results on convergence acceleration are presented, while in
  Section \ref{secStaRes}, results on the numerical stability of
  the transformations are provided. Finally, we discuss
  guidelines for the application of the transformations and some
  numerical examples in Section \ref{secApplic}.

{\quad}\section{Definitions and Notations}

  \subsection{General Definitions}

    \subsubsection{Sets}
       \paragraph{Natural Numbers:}
       \begin{equation}\label{eq1}
                       \Nset = \Set{1,2,3,...}, \qquad \Nset_0 = \Nset\cup\Set{0}
       \end{equation}

       \paragraph{Integer Numbers:}
       \begin{equation}\label{eq2}
                       \Zset = \Nset \cup \Set{0,-1,-2,-3,...}
       \end{equation}

       \paragraph{Real Numbers and Vectors:}
       \begin{equation}\label{eq3}
         \begin{array}{>{\displaystyle}r@{}>{\displaystyle}l}
           &\Rset = \Set{x:\, x \mbox{ real}},\\
           &\Rset_{+} =\Set{x\in\Rset: x>0},\\
           &\Rset^n = \Set{(x_1,...,x_n) \>\vert\> x_j\in\Rset,j=1,...,n}
         \end{array}
       \end{equation}

       \paragraph{Complex Numbers:}
       \begin{equation}\label{eq4}
         \begin{array}{>{\displaystyle}r@{}>{\displaystyle}l}
           &\Cset = \Set{z=x+\mbox{i}\, y:\, x\in \Rset,\, y\in\Rset,\,
             \mbox{i}^2=-1}\\
           &\Cset^n = \Set{(z_1,...,z_n) \>\vert\> z_j\in\Cset,j=1,...,n}
         \end{array}
       \end{equation}
       For $z=x+\mbox{i}\, y$, real and imaginary parts are denoted as
       $x=\Re(z)$, $y=\Im(z)$.
       We use $\Kset$ to denote $\Rset$ or $\Cset$.

       \paragraph{Vectors with nonvanishing components:}
       \begin{equation}\label{eq5}
           \Fset^n = \Set{(z_1,...,z_n) \>\vert\> z_j\in\Cset, z_j\ne 0,j=1,...,n}
       \end{equation}

        \paragraph{Polynomials:}
        \begin{equation}
            \Pset^{k} = \left.\Set{P: z\mapsto \sum_{j=0}^{k} c_j z^j
            \>\right\vert\>
            z\in\Cset, (c_0,\dots,c_k) \in \Kset^{k+1}}
        \label{eqpoly}\end{equation}

       \paragraph{Sequences:}
       \begin{equation}\label{eqseq}
           \Sset^{\Kset} = \Set{\Seq{s_0,s_1,...,s_n,\dots} \>\vert\>
           s_n\in\Kset,\>
           n\in\Nset_0}
       \end{equation}

       \paragraph{Sequences with nonvanishing terms:}
       \begin{equation}\label{eqnonzeroseq}
           \Oset^{\Kset} = \Set{\Seq{s_0,s_1,...,s_n,\dots} \>\vert\>
           s_n\ne 0,\> s_n\in\Kset,\>
           n\in\Nset_0}
       \end{equation}

    \subsubsection{Special Functions and Symbols}
       \paragraph{Gamma function}
       \cite[p. 1]{MagnusOberhettingerSoni66}:
       \begin{equation}\label{eq6}
         \Gamma(z) = \int_0^{\infty} t^{z-1}\> \exp({-t}) \d t \quad
         (z\in\Rset_{+})
       \end{equation}

       \paragraph{Factorial:}
       \begin{equation}\label{eq7}
         n! = \Gamma(n+1) = \prod_{j=1}^{n} j
       \end{equation}


       \paragraph{Pochhammer Symbol} \cite[p.
       2]{MagnusOberhettingerSoni66}:
       \begin{equation}\label{eqpoch}
         (a)_n = \frac{\Gamma(a+n)}{\Gamma(a)} = \prod_{j=1}^{n} (a+j-1)
       \end{equation}

       \paragraph{Binomial Coefficients} \cite[p. 256, Eq.
       (6.1.21)]{AbramowitzStegun70}:
       \begin{equation}\label{eq10}
         \binom{z}{w} = \frac{\Gamma(z+1)}{\Gamma(w+1)\Gamma(z-w+1)}
       \end{equation}

       \paragraph{Entier Function:}
       \begin{equation}\label{eq11}
       \Ent{x} = \max \Set{j\in\Zset:\, j\le x,\, x\in\Rset}
       \end{equation}

  \subsection{Sequences, Series and Operators}
    \subsubsection{Sequences and Series}\label{secseqser}
      For Stieltjes series see Appendix \ref{appstieltjes}.

      \paragraph{Scalar Sequences} with Elements $s_n$, Tail $R_n$,
      and Limit $s$:
      \begin{equation}\label{eq31}
        \Seq{s_n} =
        \Seq{s_n}_{n=0}^{\infty}=\Seq{s_0,s_1,s_2,\dots}
        \in\Sset^{\Kset}\>, \quad R_n = s_n -s \>,\quad
        \lim_{n\to\infty} s_n =s\>.
      \end{equation}
      If the sequence is not convergent but summable to $s$, $s$ is
      called the Antilimit. The $n$-th element $s_n$ of a sequence
      $\sigma=\Seq{s_n}\in \Sset^{\Kset}$ is also denoted by $\Ele{
      \sigma}_n$. A sequence is called a \emph{constant
      sequence}, if all elements are constant, i.e., if there is a
      $c\in\Kset$ such that $s_n=c$ for all $n\in\Nset_0$, in which
      case it is denoted by $\Seq{c}$. The constant sequence $\Seq{0}$
      is called the \emph{zero sequence}.

      \paragraph{Scalar Series} with Terms $a_j\in \Kset$, Partial Sums
      $s_n$, Tail $R_n$, and Limit/Antilimit $s$:
      \begin{equation}\label{eq32}
         s = \sum_{j=0}^{\infty} a_j\>, \quad s_n = \sum_{j=0}^{n}
         a_j\>, \quad R_n = -\sum_{j=n+1}^{\infty}
         a_j = s_n - s
      \end{equation}
      We say that $\hat a_n$ are
      Kummer-related to the $a_n$ with limit or antilimit $\hat s$ if
      $\hat
      a_n=\forwarddiff \hat s_{n-1}$ satisfy
      $a_n\sim \hat a_n$ for $n\to\infty$ and $\hat s$ is the limit (or
      antilimit) of $\hat s_n=\sum_{j=0}^{n}  \hat a_j$.

      \paragraph{Scalar Power Series} in $z\in\Cset$ with coefficients
      $c_j\in \Kset$, Partial Sums $f_n(z)$, Tail $R_n(z)$, and Limit/Antilimit $f(z)$:
      \begin{equation}\label{eq33}
         f(z) = \sum_{j=0}^{\infty} c_j z^j\>, \quad f_n(z) = \sum_{j=0}^{n}
         c_j z^j\>, \quad R_n(z) = \sum_{j=n+1}^{\infty}
         c_j z^j = f(z)- f_n(z)
      \end{equation}


    \subsubsection{Types of Convergence}\label{secconvtyp}
      Sequences $\Seq{s_n}$ satisfying the equation
      \begin{equation}\label{eqbekaverf_rho}
        \lim_{n\to\infty}
        ({s_{n+1}-s})/({s_n-s}) = \rho
      \end{equation}
      are called \emph{linearly convergent} if $0<\abs{\rho}<1$,
      {\em logarithmically convergent} for $\rho=1$ and
      \emph{hyperlinearly convergent} for $\rho=0$. For $\abs{\rho}>1$,
      the sequence diverges.

      A sequence $\Seq{u_n}$ accelerates a sequence
      $\Seq{v_n}$ to $s$ if
      \begin{equation}\label{eq36}
        \lim_{n\to\infty}
        ({u_{n}-s})/({v_n-s}) = 0\>.
      \end{equation}
      If $\Seq{v_n}$ converges to $s$ then we also say that $\Seq{u_n}$
      converges faster than $\Seq{v_n}$.

      A sequence $\Seq{u_n}$ accelerates a sequence
      $\Seq{v_n}$ to $s$ with order $\alpha>0$ if
      \begin{equation}\label{eq37}
        ({u_{n}-s})/({v_n-s}) = O(n^{-\alpha})\>.
      \end{equation}
      If $\Seq{v_n}$ converges to $s$ then we also say that $\Seq{u_n}$
      converges faster than $\Seq{v_n}$ with order $\alpha$.

    \subsubsection{Operators}

       \paragraph{Annihilation Operator:}
         An operator $\A \>:\> \Sset^{\Kset} \longrightarrow
         {\Kset}$ is called an annihilation operator for a given
         sequence $\Seq{\tau_n}$ if it satisfies
         \begin{equation}
           \begin{array}{>{\displaystyle}l}
                 \A (\Seq{s_n + z t_n}) = \A(\Seq{s_n})+ z\A(\Seq{t_n}) \mbox{
                for all }  \Seq{s_n} \in \Sset^{\Kset}\>, \> \Seq{t_n} \in
                \Sset^{\Kset}\>, \> z\in\Kset \\
                 \A(\Seq{\tau_n}) = 0
           \end{array}
           \label{eqannihil}
         \end{equation}

       \paragraph{Forward Difference Operator:}
       \begin{equation}\label{eq28}
         \begin{array}{>{\displaystyle}r@{}>{\displaystyle}l}
          & \forwarddiff_m g(m) = g(m+1) -g(m)\>,
            \forwarddiff_m g_m = g_{m+1} -g_{m}\\
          & \forwarddiff_m^k = \forwarddiff_m \forwarddiff_m^{k-1}\\
          & \forwarddiff = \forwarddiff_n \\
          & \forwarddiff^k g_n = \sum_{j=0}^{k} (-1)^{k-j} \binom{k}{j}
               g_{n+j}
         \end{array}
       \end{equation}

       \paragraph{Generalized Difference Operator} $\nabla\di{n}{k}$ for
       given quantities $\delta\di{n}{k}\ne 0$:
       \begin{equation}\label{eq29}
          \nabla\di{n}{k} = (\delta\di{n}{k})^{-1} \forwarddiff
       \end{equation}

       \paragraph{Generalized Difference Operator} $\widetilde\nabla\di{n}{k}$ for
       given quantities $\zeta\di{n}{k}\ne 0$:
       \begin{equation}\label{eq29a}
          \widetilde\nabla\di{n}{k} = (\zeta\di{n}{k})^{-1} \forwarddiff^2
       \end{equation}

       \paragraph{Generalized Difference Operator}
       $\bigtriangledown\di{n}{k}[\alpha]$ for
       given quantities $\Delta\di{n}{k}\ne 0$:
       \begin{equation}\label{eq29b}
          \bigtriangledown\di{n}{k}[\alpha] f_n =
          (\Delta\di{n}{k})^{-1} (f_{n+2} - 2\cos\alpha f_{n+1} + f_n)
       \end{equation}

       \paragraph{Generalized Difference Operator}
       $\partial\di{n}{k}[\zeta]$ for
       given quantities $\widetilde\Delta\di{n}{k}\ne 0$:
       \begin{equation}\label{eq29c}
          \partial\di{n}{k}[\zeta] f_n
           =
           (\widetilde\Delta\di{n}{k})^{-1} (\zeta\di{n+k}{2} f_{n+2}
           + \zeta\di{n+k}{1}  f_{n+1} + \zeta\di{n+k}{0}f_n)
       \end{equation}

       \paragraph{Weighted Difference Operators}  for
       given  $P\di{}{k-1}\in
       \Pset^{k-1}$:
       \begin{equation}\label{eqweighteddiff}
          \W\di{n}{k} =\W\di{n}{k}[P\di{}{k-1}] = \forwarddiff^k
          P\di{}{k-1}(n)
       \end{equation}

       \paragraph{Polynomial Operators} \P for given $P\di{}{k}\in
       \Pset^{k}$: Let $P\di{}{k}(x)=\sum_{j=0}^{k} p\di{j}{k} x^j$.
       Then put
       \begin{equation} \label{eqPolyOp}
         \P[P\di{}{k}] g_n = \sum_{j=0}^{k} p\di{j}{k} g_{n+j}\>.
       \end{equation}

       \paragraph{Divided Difference Operator:}
       For given $\Seq{x_n}$ and $k,n\in \Nset_0$, put
       \begin{equation}\label{eq30}
         \begin{array}{>{\displaystyle}r@{}>{\displaystyle}l}
           & \square\di{n}{k}[\Seq{x_n}](f(x)) =\square\di{n}{k}(f(x)) = f[x_n,\dots,x_{n+k}]
                              = \sum_{j=0}^{k} f(x_{n+j})
                                \prod_{i=0 \atop i\ne j }^{k}
                                   \frac{1}{x_{n+j}-x_{n+i}}\>\\
           & \square\di{n}{k}[\Seq{x_n}] g_n = \square\di{n}{k} g_n =
                              \sum_{j=0}^{k} g_{n+j}
                                \prod_{i=0 \atop i\ne j }^{k}
                                   \frac{1}{x_{n+j}-x_{n+i}}
          \end{array}
       \end{equation}

  \section{Some Basic Sequence Transformations}\label{secBasSeq}
      \subsection{\E\ Algorithm}\label{secEalgo}
Putting for  sequences $\Seq{y_n}$ and $\Seq{g_j(n)}$, $j=1,\dots,k$
\begin{equation}
E\di{n}{k}[\Seq{y_n};\Seq{g_j(n)}]=\left\vert
\begin{array}{>{\displaystyle}cc>{\displaystyle}c}
 y_n & \cdots & y_{n+k} \\
g_1(n)   & \cdots & g_1(n+k)     \\
\vdots   & \ddots & \vdots       \\
g_k(n)   & \cdots & g_k(n+k)     \\
\end{array}
\right\vert
\end{equation}
one may define the sequence transformation
\begin{equation}
\E\di{n}{k}(\Seq{s_n})  =  \frac{E\di{n}{k}[\Seq{s_n};\Seq{g_j(n)}]}
                                {E\di{n}{k}[\Seq{1};\Seq{g_j(n)}]}
\end{equation}
As is plain using Cramer's rule, we have $\E\di{n}{k}(\Seq{\sigma_n})=\sigma$
if the $\sigma_n$ satisfy Eq.\ (\ref{eqbekaverf_Emodel}).
Thus, the sequence transformation yields the limit $\sigma$ exactly for
model sequences (\ref{eqbekaverf_Emodel}).

The sequence transformation \E\ is known as the \E\ algorithm or also as
{Brezinski--H{\aa}vie--Protocol}  \cite[Sec.\ 10]{Wimp81} after two of
its main investigators, H{\aa}vie \cite{Haavie79} and Brezinski
\cite{Brezinski80b}. A good introduction to this transformation is also
given in the book of
Brezinski and Redivo Zaglia \cite[Sec.\
2.1]{BrezinskiRedivoZaglia91}. Compare also Ref.\
\cite{BrezinskiRedivoZaglia93a}.

Numerically, the computation of the
$\E\di{n}{k}(\Seq{s_n})$ can be performed recursively using either the
algorithm of Brezinski  \cite[p.
58f]{BrezinskiRedivoZaglia91}
\begin{eqnarray}
&{}&\E\di{n}{0}(\Seq{s_n})=s_n,\qquad
g\di{0,i}{n}=g_i(n), \qquad n\in\Nset_0, \> i\in\Nset
\nonumber\\
&{}&\E\di{n}{k}(\Seq{s_n}) = \E\di{n}{k-1}(\Seq{s_n}) -
\frac{\E\di{n+1}{k-1}(\Seq{s_n})-\E\di{n}{k-1}(\Seq{s_n})}
{g\di{k-1,k}{n+1}-g\di{k-1,k}{n}} \cdot g\di{k-1,k}{n}
\nonumber\\
&{}&g\di{k,i}{n} = g\di{k-1,i}{n} -
\frac{g\di{k-1,i}{n+1}-g\di{k-1,i}{n}}
{g\di{k-1,k}{n+1}-g\di{k-1,k}{n}} \cdot g\di{k-1,k}{n}, \qquad
i=k+1,k+2,\dots
\nonumber\\
\label{eqbekaverf_Ealgo}
\end{eqnarray}
or the algorithm of Ford and Sidi \cite{FordSidi87} that requires
additionally the quantities
$g_{k+1}(n+j)$, $j=0,\dots,k$ for the computation of
$\E\di{n}{k}(\Seq{s_n})$. The algorithm of Ford and Sidi involves the
quantities
\begin{equation}\label{eqbekaverf_Psi}
\Psi_{k,n}(u)=
\frac{E\di{n}{k}[\Seq{u_n};\Seq{g_j(n)}]}
     {E\di{n}{k}[\Seq{g_{k+1}(n)};\Seq{g_j(n)}]}
\end{equation}
for any sequence $\Seq{u_0,u_1,\ldots}$ where the $g_i(n)$ are not
changed even if they depend on the
$u_n$ and the $u_n$ are changed. Then we have
\begin{equation}\label{eqbekaverf_EFSalgo}
\E\di{n}{k}(\Seq{s_n}) =
\frac{\Psi\di{k}{n}(s)}{\Psi\di{k}{n}(1)}\>,
\end{equation}
and the  $\Psi$ are calculated  recursively via
\begin{equation}\label{eqbekaverf_FSalgo}
\Psi_{k,n}(u) =
\frac{\Psi_{k-1,n+1}(u)-\Psi_{k-1,n}(u)}
{\Psi_{k-1,n+1}(g_{k+1})-\Psi_{k-1,n}(g_{k+1})}\>.
\end{equation}

Of course, for $g_j(n)=\omega_n \psi_{j-1}(n)$, i.e., in the context of
sequences modelled via the expansion (\ref{eqpsijn}),  the \E\ algorithm may
be used to obtain an explicit representation for any Levin-type
sequence transformation of the form (compare Eq.\ (\ref{eqT0}))
\begin{equation}
\T\di{n}{k}=
T(s_n,\dots,s_{n+k};\omega_{n},\dots,\omega_{n+k};
\psi_j(n),\dots,\psi_j(n+k))
\end{equation}
as ratio of two determinants
\begin{equation}
\T\di{n}{k}(\Seq{s_n},\Seq{\omega_n})  =
  \frac{E\di{n}{k}[\Seq{s_n/\omega_n};\Seq{\psi_{j-1}(n)}]}
       {E\di{n}{k}[\Seq{1/\omega_n};\Seq{\psi_{j-1}(n)}]}\>.
\end{equation}
This follows from the identity \cite{BrezinskiRedivoZaglia91}
\begin{equation}
  \frac{E\di{n}{k}[\Seq{s_n};\Seq{\omega_n\psi_{j-1}(n)}]}
       {E\di{n}{k}[\Seq{1};\Seq{\omega_n\psi_{j-1}(n)}]}
=  \frac{E\di{n}{k}[\Seq{s_n/\omega_n};\Seq{\psi_{j-1}(n)}]}
       {E\di{n}{k}[\Seq{1/\omega_n};\Seq{\psi_{j-1}(n)}]}
\end{equation}
that is an easy consequence of usual algebraic manipulations of determinants.

\subsection{The $d\di{}{m}$ transformations}\label{appdm}
       As noted in the introduction, the  $d^{(m)}$ transformations were introduced by Levin and Sidi
       \cite{LevinSidi81} as a generalization of the $u$ variant of the
       Levin transformation \cite{Levin73}.  
       We describe a slightly modified variant of the
       $d^{(m)}$ transformations \cite{Sidi95}:
       
	Let $s_r$, $r=0,1,\dots$ be a real or complex sequence 
        with limit or antilimit $s$ and terms $a_0=s_0$ and
	$a_r=s_r-s_{r-1}$, $r=1,2,\dots$ such that $s_r=\sum_{r=0}^r a_j$, $r=0,1,\dots$.
      For given $m\in \Nset$ and $\xi_l\in\Nset_0$ with $l\in\Nset_0$
       and $0\le\xi_0<\xi_1<\xi_2<\cdots$ and $\nu=(n_1,\dots,n_m)$ with $n_j\in
       \Nset_0$ the $d^{(m)}$ transformation yields a table of
       approximations $s\di{\nu}{m,j}$ for the (anti-)limit $s$ 
       as
       solution of the linear system of equations
       \begin{equation}\label{eqbekaverf_d}
          s_{\xi_l} = s\di{\nu}{m,j} +
           \sum_{k=1}^{m}
              (\xi_l+\alpha)^{k}
              [\Delta^{k-1} a_{\xi_l}]
              \sum_{i=0}^{n_k}
                \frac{\bar\beta_{ki}}
                     {(\xi_l+\alpha)^{i}}\,, \quad
                                           j\le l\le j+N\,,
       \end{equation}
       with $\alpha>0$, $N=\sum_{k=1}^{m} n_k$ and the $N+1$ unknowns
       $s\di{\nu}{m,j}$ and $\bar\beta_{k\,i}$. The $[\Delta^k a_j]$
       are defined via
       $[\Delta^0 a_j ]= a_j$ and $[\Delta^k a_j]= [\Delta^{k-1} a_{j+1}] -
       [\Delta^{k-1} a_j]$, $k=1,2,\dots$. In most cases, all $n_k$ are chosen equal and one
       puts  $\nu=(n,n,\dots,n)$.  Apart
from the value of $\alpha$, only the input of $m$ and of $\xi_\ell$ is required from the
user.
       As transformed sequence, often one chooses
       the elements $s_{(n,\dots,n)}^{(m,0)}$ for $n=0,1,\dots$.
       The $u$ variant of the Levin transformation is obtained for
       $m=1$, $\alpha=\beta$ and $\xi_l=l$.
        The  definition above differs slightly from the original one
       \cite{LevinSidi81} and was given in Ref.\ \cite{FordSidi87} with $\alpha=1$.

       Ford and Sidi have shown, how these
       transformations can be calculated recursively with the
       ${\mathbf W}^{(m)}$ al\-go\-rithms \cite{FordSidi87}.
       The $d^{(m)}$ transformations are the best known special
       cases of the \emph{generalised Richardson Extrapolation process}
       (GREP) as defined by Sidi \cite{Sidi79b,Sidi82,Sidi95caf}. 

       The $d^{(m)}$ transformations are derived by asymptotic analysis of the remainders
       $s_r-s$ for $r\to\infty$ for the family $\tilde B^{(m)}$ of sequences
       $\Seq{a_r}$ as defined in Ref.\ \cite{LevinSidi81}. For such sequences, the $a_r$
satisfy a difference equation of order $m$ of the form
       \begin{equation}
         a_r = \sum_{k=1}^{m} p_k(r) \Delta^k\,a_r\>.
       \end{equation}
       The $p_k(r)$ satisfy the asymptotic relation
       \begin{equation}
         p_k(r) \sim r^{i_k} \sum_{\ell=0}^{\infty}
\frac{p_{k\ell}}{r^{\ell}}\>\qquad\mbox{for } r\to\infty\>.
       \end{equation}
       The $i_k$ are integers satisfying $i_k\le k$ for $k=1,\dots,m$.
       This family of sequences is very large. But still, Levin and Sidi could prove
\cite[Theorem 2]{LevinSidi81} that
under mild additional assumptions, the remainders for such sequences satisfy
       \begin{equation}
         s_r-s \sim \sum_{k=1}^{m} r^{j_k}(\Delta^{k-1}\,a_r) \sum_{\ell=0}^{\infty}
\frac{\beta_{k\ell}}{r^{\ell}}\>\qquad\mbox{for } r\to\infty\>.
       \end{equation}
       The $j_k$ are integers satisfying $j_k\le k$ for $k=1,\dots,m$. A corresponding
result for $m=1$ was proven by Sidi \cite[Theorem 6.1]{Sidi79}.
   
       The system (\ref{eqbekaverf_d}) now is obtained by truncation of the expansions at
$\ell=n_n$, evaluation at $r=\xi_l$, and some further obvious substitutions.

       The introduction of suitable $\xi_l$ was shown to improve the accuracy and stability
in difficult situations considerably \cite{Sidi95}.

      \subsection{Shanks Transformation and Epsilon Algorithm}
        An important special case of the \E\ algorithm is the choice
        $g_j(n) = \forwarddiff s_{n+j-1}$ leading to the
        Shanks transformation \cite{Shanks55}
       \begin{equation} \label{eqShanks}
       e_k(s_n) = \frac{E\di{n}{k}[\Seq{s_n};\Seq{\forwarddiff s_{n+j-1}}]}
                       {E\di{n}{k}[\Seq{1};\Seq{\forwarddiff s_{n+j-1}}]}
       \end{equation}
       Instead of using one of the recursive schemes for the \E\
       algorithms, the Shanks transformation may be implemented using
        the epsilon algorithm \cite{Wynn56a} that is defined by the recursive scheme
        \begin{equation}\label{eqbekaverf_epsilon}
          \begin{array}{>{\displaystyle}r@{}>{\displaystyle}l}
            \epsilon_{-1}^{(n)}  &{} = 0 \> , \qquad\epsilon_0^{(n)} =  s_n\>, \\
            \epsilon_{k+1}^{(n)} &{} = \epsilon_{k-1}^{(n+1)} \, + \,
            1 / [\epsilon_{k}^{(n+1)} - \epsilon_{k}^{(n)} ]\>. \\
          \end{array}
        \end{equation}
        The relations
        \begin{equation}
          \epsilon_{2 k}^{(n)} \; = \; e_k (s_n)\>,\quad
          \epsilon_{2k+1}^{(n)} \; = \; 1/e_k (\forwarddiff s_n)\>
        \end{equation}
        hold and show that the elements $\epsilon_{2k+1}^{(n)}$ are
        only auxiliary quantities.

        The kernel of the Shanks transformation $e_k$ is given by
        sequences of the form
        \begin{equation}
            s_n \; = \; s + \sum_{j=0}^{k-1} \> c_j \, \forwarddiff s_{n+j}\>.
        \end{equation}
        See also \cite[Theorem 2.18]{BrezinskiRedivoZaglia91}.

        Additionally, one can use the Shanks transformation --- and hence
        the epsilon al\-go\-rithm --- to compute the upper half of the
        Pad\'e table according to  \cite{Shanks55,Wynn56a}
\begin{equation}
e_k ( f_n (z) ) \; = \; [ \, n \, + \, k \, / \, k \, ]_f \, (z) \, ,
\qquad (k\ge 0\;,n\ge 0)\>.
\end{equation}
where
\begin{equation}
f_n (z) \; = \; \sum_{j=0}^n \> c_{j} \, z^{j}
\end{equation}
are the partial sums of a power series of a function $f(z)$.
{Pad\'e approximants} of $f (z)$ are rational functions in $z$ given as
ratio of two polynomials $p_{\ell}\in \Pset\di{}{\ell}$ and $q_m\in
\Pset\di{}{m}$
according to
\begin{equation}\label{eqbekaverf_Pade}
[ \>\ell \> / \> m \>]_f \>(z) \; = \;
p_{\ell} (z) \> / \> q_m (z),
\end{equation}
where the Taylor series of $f$ and $[\ell/m]_f$ are identical to the
highest possible power of $z$, i.e.
\begin{equation}
f (z) \, - \, p_{\ell} (z) \, /  \,  q_m (z)  \; = \; O(z^{\ell + m
+1}).
\end{equation}
Methods for the extrapolation of power series will be treated later.

        \subsection{Aitken Process}

        The special case $\epsilon\di{2}{n}=e_1(s_n)$ is identical to
        the famous $\Delta^2$ method of Aitken \cite{Aitken26}
\begin{equation}\label{eqbekaverf_Aitken}
s\di{n}{1}= s_n - \frac{(s_{n+1}-s_n)^2}{s_{n+2}-2s_{n+1}+s_n}\>
\end{equation}
with kernel
\begin{equation}\label{eqbekaverf_Aitkenkernel}
s_n = s + c\,(s_{n+1}-s_{n}), \qquad n\in\Nset_0\>.
\end{equation}
Iteration of the $\Delta^2$ method yields the { iterated
Aitken process}
\cite{Wimp81,Weniger89,BrezinskiRedivoZaglia91}
\begin{equation}\label{eqbekaverf_itAitken}
\begin{array}{>{\displaystyle}r@{}>{\displaystyle}l}
 &\mathbf{A}\di{n}{0}=s_n\>, \\
&\mathbf{A}\di{n}{k+1}=
\mathbf{A}\di{n}{k} -
\frac{(\mathbf{A}\di{n+1}{k}-\mathbf{A}\di{n}{k})^2}
{\mathbf{A}\di{n+2}{k}-2\mathbf{A}\di{n+1}{k}
+\mathbf{A}\di{n}{k}}\>.\\
\end{array}
\end{equation}
The iterated Aitken process and the epsilon algorithm
accelerate linear convergence and can sometimes be applied successfully
for the summation of alternating divergent series.

      \subsection{Overholt Process}
        The Overholt process is defined by the recursive scheme
        \cite{Overholt65}
        \begin{equation} \label{eqOverholtP}
          \begin{array}{>{\displaystyle}r@{}>{\displaystyle}l}
             V\di{n}{0}(\Seq{s_n}) = {}
             &   s_n\;, \\
             V\di{n}{k}(\Seq{s_n}) {}=
             & \frac{(\forwarddiff s_{n+k-1})^k
                     V\di{n+1}{k-1}(\Seq{s_n})
                     -(\forwarddiff s_{n+k})^kV
                     \di{n}{k-1}(\Seq{s_n})}
                    {(\forwarddiff s_{n+k-1})^k
                     -(\forwarddiff s_{n+k})^k}
           \end{array}
        \end{equation}
        for $k\in\Nset$ and $n\in\Nset_0$. It is important for the
        convergence acceleration of fixed point iterations.

{\quad}\section{Levin-Type Sequence
Transformations}\label{secLevSeq}

   \subsection{Definitions for Levin-Type Transformations}
      A set $\Lambda\di{}{k}=\Set{\lambda\di{n,j}{k}\in \Kset\>\vert\>
      n\in\Nset_0\>,\> 0\le j\le k}$ is called a
      \emph{coefficient set of order k} with $k\in\Nset$  if
      $\lambda\di{n,k}{k}\ne 0$
      for all
      $n\in\Nset_0$.
      Also,
      $\Lambda=\Set{\Lambda\di{}{k}\>\vert\> k\in \Nset}$ is called
      \emph{coefficient set}. Two coefficient sets
      $\Lambda=\Set{\Set{\lambda\di{n,j}{k}}}$ and
      $\widehat\Lambda=\Set{\Set{\widehat\lambda\di{n,j}{k}}}$ are
      called \emph{equivalent}, if for all $n$ and $k$, there is a
      constant $c\di{n}{k}\ne 0$ such that
      $\widehat\lambda\di{n,j}{k}=c\di{n}{k} \lambda\di{n,j}{k}$ for all
      $j$ with $0\le j\le k$.

      For each
      coefficient set
      $\Lambda\di{}{k}=\Set{\lambda\di{n,j}{k}\>\vert\> n\in\Nset_0, 0\le j\le
      k}$ of order $k$, one may define a \emph{Levin-Type sequence
      transformation of order $k$}  by
      \begin{equation} \label{eqLevinType}
         \begin{array}{>{\displaystyle}r>{\displaystyle}l>{\displaystyle}l>{\displaystyle}l}
           \T[\Lambda\di{}{k}] \>:\> & \Sset^{\Kset} \times
           \Yset\di{}{k} & \longrightarrow & \Sset^{\Kset}\\
                       \>:\> & (\Seq{s_n},\Seq{\omega_n}) & \longmapsto
                       &
                           \Seq{s_n'} =\T[\Lambda\di{}{k}]
                                         (\Seq{s_n},\Seq{\omega_n})

          \end{array}
      \end{equation}
      with
      \begin{equation} \label{eqsnd}
          s_n' =\T\di{n}{k}(\Seq{s_n},\Seq{\omega_n})= \frac{\displaystyle\sum_{j=0}^{k} \lambda\di{n,j}{k}
                       \frac{s_{n+j}}{\omega_{n+j}}}
                      {\displaystyle\sum_{j=0}^{k} \lambda\di{n,j}{k}
                       \frac{1}{\omega_{n+j}}}
      \end{equation}
      and
      \begin{equation} \label{eqYnk}
        \Yset\di{}{k} = \Set{ \Seq{\omega_n} \in \Oset^{\Kset} \>:\>
        \sum_{j=0}^{k}\lambda\di{n,j}{k} /\omega_{n+j} \ne 0 \mbox{ for
        all } n\in \Nset_0 }
      \end{equation}
      We call $\T[\Lambda]=\Set{\T[\Lambda\di{}{k}]\>\vert\>
      k\in\Nset}$ the Levin-type sequence transformation
      corresponding to the coefficient set
      $\Lambda=\Set{\Lambda\di{}{k}\>\vert\> k\in \Nset}$. We
      write $\T\di{}{k}$ and \T\ instead of $\T[\Lambda\di{}{k}]$ and
      $\T[\Lambda]$, respectively,
      whenever the
      coefficients $\lambda\di{n,j}{k}$ are
      clear from the context.  Also, if two coefficient sets $\Lambda$ and
      $\widehat\Lambda$ are equivalent, they give rise to the same
      sequence transformation, i.e.,  $\T[\Lambda]=\T[\widehat
      \Lambda]$, since
      \begin{equation} \label{eqsndscal}
           \frac{\displaystyle\sum_{j=0}^{k} \widehat\lambda\di{n,j}{k}
                       \frac{s_{n+j}}{\omega_{n+j}}}
                {\displaystyle\sum_{j=0}^{k} \widehat\lambda\di{n,j}{k}
                       \frac{1}{\omega_{n+j}}}
           = \frac{\displaystyle\sum_{j=0}^{k} \lambda\di{n,j}{k}
                       \frac{s_{n+j}}{\omega_{n+j}}}
                      {\displaystyle\sum_{j=0}^{k} \lambda\di{n,j}{k}
                       \frac{1}{\omega_{n+j}}}
           \quad \mbox{for }\widehat\lambda\di{n,j}{k}=c\di{n}{k} \lambda\di{n}{k}
      \end{equation}
      with arbitrary $c\di{n}{k}\ne 0$.

      The number $\T\di{n}{k}$ are often arranged in a two-dimensional
      table
      \begin{equation}
        \begin{array}{llll}
          \T\di{0}{0} & \T\di{0}{1} & \T\di{0}{2} &\cdots \\
          \T\di{1}{0} & \T\di{1}{1} & \T\di{1}{2} &\cdots \\
          \T\di{2}{0} & \T\di{2}{1} & \T\di{2}{2} &\cdots \\
          \vdots      &  \vdots     & \vdots      & \ddots
        \end{array}
      \end{equation}
      that is called the \T\ table. The transformations $\T\di{}{k}$
      thus correspond to columns, i.e., to following vertical
      paths in the table. The numerators and denominators such
      that $\T\di{n}{k}=N\di{n}{k}/D\di{n}{k}$ also
      are often arranged in analogous $N$ and $D$ tables.

      Note that for fixed $N$, one may also define a transformation
      \begin{equation}
        \T_N : \Seq{s_{n+N}} \longmapsto
        \Seq{\T\di{N}{k}}_{k=0}^{\infty}\>.
      \end{equation}
      This corresponds to horizontal paths in the \T\ table. These are
      sometimes called diagonals, because rearranging
      the table in such a way that elements with constant values of
      $n+k$ are members of the same row, $\T\di{N}{k}$ for fixed $N$
      correspond to diagonals of the rearranged table.

      For a given coefficient set $\Lambda$
      define the \emph{moduli} by
      \begin{equation} \label{eqmunk}
        \mu\di{n}{k} = \max_{0\le j\le k}
        \Set{\abs{\lambda\di{n,j}{k}}}
      \end{equation}
      and the \emph{characteristic polynomials} by
      \begin{equation} \label{eqpink}
        \Pi\di{n}{k}\in \Pset^{k}\>:\> \Pi\di{n}{k}(z) = \sum_{j=0}^{k}
        \lambda\di{n,j}{k} z^j\>
      \end{equation}
      for $n\in\Nset_0$ and $k\in\Nset$.

      Then, $\T[\Lambda]$ is said to be \emph{in normalized form} if
      $\mu\di{n}{k}=1$ for all $k\in\Nset$ and $n\in\Nset_0$. Is is said to
      be \emph{in subnormalized form} if for all $k \in \Nset$ there
      is a constant $\widetilde\mu\di{}{k}$ such that  $\mu\di{n}{k}\le
      \widetilde\mu\di{}{k}$ for all $n\in\Nset_0$.

      Any Levin-type sequence transformation $\T[\Lambda]$
      can rewritten in normalized form.  To see this, use
      \begin{equation} \label{eqsndnormc}
        c\di{n}{k}=1/\mu\di{n}{k}
      \end{equation}
      in Eq.\ (\ref{eqsndscal}).
      Similarly, each Levin-type sequence transformation can be
      rewritten in (many different) subnormalized forms.

      A Levin-type sequence
      transformation of order $k$ is said to be
      \emph{convex} if $\Pi\di{n}{k}(1)=0$ for all $n$ in
      $\Nset_0$. Equivalently, it is convex if $\Seq{1}\not\in
      \Yset\di{}{k}$, i.e., if the transformation vanishes for
      $\Seq{s_n}=\Seq{c\omega_n}$, $c\in \Kset$. Also, $\T[\Lambda]$ is
      called convex, if $\T[\Lambda\di{}{k}]$ is convex for all
      $k\in\Nset$. We will see that this property is important for
      ensuring convergence acceleration for linearly convergent
      sequences.

      A given Levin-type transformation \T\ can also be rewritten as
\begin{equation}\label{eqTgamma}
\T\di{n}{k}(\Seq{s_n},\Seq{\omega_n}) =
\sum_{j=0}^{k} \gamma\di{n,j}{k}(\vec \omega_n)
\, s_{n+j}\>,\qquad \vec \omega_n=(\omega_n,\dots,\omega_{n+k})\>,
\end{equation}
with
\begin{equation}\label{eqgammanjk}
\gamma\di{n,j}{k}(\vec \omega_n)
= \frac{\lambda\di{n,j}{k}}{\omega_{n+j}}
\left[\sum_{j'=0}^{k}\frac{\lambda\di{n,j'}{k}}{\omega_{n+j'}}\right]^{-1}
\>,\qquad \sum_{j=0}^{k} \gamma\di{n,j}{k}(\vec \omega_n)=1\>.
\end{equation}
      Then, one may define \emph{stability
      indices} by
      \begin{equation} \label{eqGammank}
            {\boldsymbol \Gamma}\di{n}{k}(\T) =
            \sum_{j=0}^{k} \vert\gamma\di{n,j}{k}(\vec
            \omega_n)\vert \ge 1\>.
      \end{equation}
Note that any sequence transformation  $\mathcal{Q}$
\begin{equation}
\mathcal{Q}\di{n}{k} = \sum_{j=0}^{k} q\di{n,j}{k} s_{n+j}\>
\end{equation}
with
\begin{equation}\label{eqsumq}
\sum_{j=0}^{k} q\di{n,j}{k} =1
\end{equation}
can formally be rewritten as a Levin-type sequence transformation
according to
$\mathcal{Q}\di{n}{k}=\T\di{n}{k}(\Seq{s_n},\Seq{\omega_n})$ with
coefficients
$\lambda\di{n,j}{k}=\omega_{n+j} q\di{n,j}{k} \rho\di{n}{k}$
where the validity of Eq.\ (\ref{eqsumq}) requires to set
\begin{equation}
\rho\di{n}{k} = \sum_{j=0}^{k} \lambda\di{n,j}{k}/ \omega_{n+j}\>.
\end{equation}

      If for given $k\in\Nset$ and for a transformation
      $\T[\Lambda\di{}{k}]$  the following limits exist and
      have the values
      \begin{equation} \label{eqlimLambda}
      \lim_{n\to\infty} \lambda\di{n,j}{k} = \Ringel \lambda\di{j}{k}
      \end{equation}
      for all $0\le
      j\le k$, and if $\Ringel\Lambda\di{}{k}$ is a coefficient set of
      order $k$ which means that at least the limit
      $\Ringel\lambda\di{k}{k}$ does not vanish,
      then a limiting transformation $\Ringel \T[\Ringel \Lambda\di{}{k}]$
      exists where
      $\Ringel\Lambda\di{}{k}=\Set{\Ringel\lambda\di{j}{k}}$. More
      explicitly, we have
      \begin{equation} \label{eqRingelLevinType}
         \begin{array}{>{\displaystyle}r>{\displaystyle}l>{\displaystyle}l>{\displaystyle}l}
           \Ringel\T[\Lambda\di{}{k}] \>:\> & \Sset^{\Kset} \times
           \Ringel\Yset\di{}{k} & \longrightarrow & \Sset^{\Kset}\\
                       \>:\> & (\Seq{s_n},\Seq{\omega_n}) & \longmapsto
                       &
                           \Seq{s_n'}
          \end{array}
      \end{equation}
      with
      \begin{equation} \label{eqRingelsnd}
          s_n' =\Ringel\T\di{}{k}(\Seq{s_n},\Seq{\omega_n})=
          \frac{\displaystyle\sum_{j=0}^{k} \Ringel\lambda\di{j}{k}
                       \frac{s_{n+j}}{\omega_{n+j}}}
                      {\displaystyle\sum_{j=0}^{k} \Ringel\lambda\di{j}{k}
                       \frac{1}{\omega_{n+j}}}
      \end{equation}
      and
      \begin{equation} \label{eqRingelYk}
        \Ringel\Yset\di{}{k} = \Set{ \Seq{\omega_n} \in \Oset^{\Kset} \>:\>
        \sum_{j=0}^{k}\Ringel\lambda\di{j}{k} /\omega_{n+j} \ne 0 \mbox{ for
        all } n\in \Nset_0 }\>.
      \end{equation}
      Obviously, this limiting transformation  itself is a Levin-type
      sequence transformation and automatically is given in
      subnormalized form.

    \subsubsection{Variants of Levin-Type Transformations}
      For the following, assume that $\beta>0$ is an arbitrary
      constant, $a_n=\forwarddiff s_{n-1}$, and $\hat a_n$ are
      Kummer-related to the $a_n$ with limit or antilimit $\hat
      s$ (cp.\ Section \ref{secseqser}).

      A \emph{variant} of a Levin-Type sequence transformation \T\ is
      obtained by a particular choice $\omega_n$. For
      $\omega_n=f_n(\Seq{s_n})$, the transformation \T\ is nonlinear in the
      $s_n$.  In
      particular, we have \cite{Levin73,SmithFord79,HomeierWeniger95}:

      \paragraph{t Variant}
         \begin{equation} \label{eqtvar}
            {}^t\omega_n = \forwarddiff s_{n-1}=a_{n}\>:\>
            {}^t\T\di{n}{k}(\Seq{s_n}) = \T\di{n}{k}(\Seq{s_n},\Seq{{}^t\omega_n})
         \end{equation}

      \paragraph{u Variant}
         \begin{equation} \label{equvar}
            {}^u\omega_n = (n+\beta) \forwarddiff
            s_{n-1}=(n+\beta)a_n\>:\>
            {}^u\T\di{n}{k}(\beta,\Seq{s_n}) =
            \T\di{n}{k}(\Seq{s_n},\Seq{{}^u\omega_n})
         \end{equation}

      \paragraph{v Variant}
         \begin{equation} \label{eqvvar}
            {}^v\omega_n =
              -\frac{\forwarddiff s_{n-1}\forwarddiff s_{n}}
                   {\forwarddiff^2 s_{n-1}}
              =
              \frac{a_n a_{n+1}}{a_n - a_{n+1}}
            \>:\>
            {}^v\T\di{n}{k}(\Seq{s_n}) = \T\di{n}{k}(\Seq{s_n},\Seq{{}^v\omega_n})
         \end{equation}

      \paragraph{$\widetilde t$ Variant}
         \begin{equation} \label{eqttildvar}
            {}^{\widetilde t}\omega_n =  \forwarddiff s_{n} = a_{n+1}
            \>:\>
            {}^{\widetilde t}\T\di{n}{k}(\Seq{s_n}) =
            \T\di{n}{k}(\Seq{s_n},\Seq{{}^{\widetilde t}\omega_n})
         \end{equation}

      \paragraph{lt Variant}
         \begin{equation} \label{eqltvar}
            {}^{lt}\omega_n = \hat a_{n}\>:\>
            {}^{lt}\T\di{n}{k}(\Seq{s_n}) = \T\di{n}{k}(\Seq{s_n},\Seq{{}^{lt}\omega_n})
         \end{equation}

      \paragraph{lu Variant}
         \begin{equation} \label{eqluvar}
            {}^{lu}\omega_n =
            (n+\beta) \hat a_n\>:\>
            {}^{lu}\T\di{n}{k}(\beta,\Seq{s_n}) =
            \T\di{n}{k}(\Seq{s_n},\Seq{{}^{lu}\omega_n})
         \end{equation}

      \paragraph{lv Variant}
         \begin{equation} \label{eqlvvar}
            {}^{lv}\omega_n =
              \frac{\hat a_n \hat a_{n+1}}{\hat a_n - \hat a_{n+1}}
            \>:\>
            {}^{lv}\T\di{n}{k}(\Seq{s_n}) = \T\di{n}{k}(\Seq{s_n},\Seq{{}^{lv}\omega_n})
         \end{equation}

      \paragraph{l$\widetilde t$ Variant}
         \begin{equation} \label{eqlttildvar}
            {}^{l\widetilde t}\omega_n =  \hat a_{n+1}\>:\>
            {}^{l\widetilde t}\T\di{n}{k}(\Seq{s_n}) =
            \T\di{n}{k}(\Seq{s_n},\Seq{{}^{l\widetilde t}\omega_n})
         \end{equation}

      \paragraph{K Variant}
         \begin{equation} \label{eqKvar}
            {}^K\omega_n = \hat s_{n}-\hat s
            \>:\>
            {}^K\T\di{n}{k}(\Seq{s_n}) = \T\di{n}{k}(\Seq{s_n},\Seq{{}^K\omega_n})
         \end{equation}
         The K variant of a Levin-type transformation \T\ is linear in the
         $s_n$. This holds also for the lt, lu, lv and l$\widetilde t$
         variants.

   \subsection{Important Examples of Levin-Type Sequence Transformations}
     In this section, we present important Levin-type sequence
     transformations. For each transformation, we give the definition,
     recursive algorithms and some background information.

    \subsubsection{\J\ Transformation}
       The \J\ transformation was derived and studied by
       Homeier\cite{Homeier93,Homeier94ahc,Homeier95,Homeier96aan,Homeier96Hab,Homeier98ots}.
       Although the \J\ transformation was derived by hierarchically
       consistent iteration of the simple transformation
       \begin{equation} \label{eqJtranssimple}
         s_n' = s_{n+1} - \omega_{n+1}
         \frac{\forwarddiff s_n}
              {\forwarddiff \omega_n}
       \end{equation}
       it was possible to derive an explicit formula for its kernel
       as is discussed lated. It
       may be defined via the recursive scheme
       \begin{equation} \label{eqJtransRec}
          \begin{array}{>{\displaystyle}r@{}>{\displaystyle}l}
            & N\di{n}{0} = s_n/\omega_n\>,\quad D\di{n}{0} = 1/\omega_n\>,\\
            & N\di{n}{k} = \nabla\di{n}{k-1} N\di{n}{k-1}\>,\quad
              D\di{n}{k} = \nabla\di{n}{k-1} D\di{n}{k-1}\>,\\
            & \J\di{n}{k}(\Seq{s_n},\Seq{\omega_n},\Set{\delta\di{n}{k}}) = N\di{n}{k} / D\di{n}{k}\>,
           \end{array}
       \end{equation}
       where the generalized difference operator defined in Eq.\
       (\ref{eq29}) involves quantities $\delta\di{n}{k}\ne 0$ for
       $k\in\Nset_0$. Special
       cases of the \J\ transformation result from corresponding
       choices of the $\delta\di{n}{k}$. These are summarized in
       Table~\ref{tabJtransSpecial}.

      Using generalized difference operators $\nabla\di{n}{k}$,
      we also have the
      representation \cite[Eq.\ (38)]{Homeier94ahc}
      \begin{equation}\label{eqjtrans_31}
        \J\di{n}{k}(\Seq{s_n},\Seq{\omega_n},\Set{\delta\di{n}{k}}) =
        \frac
          {\nabla\di{n}{k-1}\nabla\di{n}{k-2}\dots \nabla\di{n}{0}
            [s_n/\omega_n]}
          {\nabla\di{n}{k-1}\nabla\di{n}{k-2}\dots \nabla\di{n}{0}
            [1/\omega_n]}
        \>.
      \end{equation}

\setlongtables
\settowidth{\mylength}{%
${}= {\J\di{n}{2k}(\Seq{s_n},\Seq{e^{-i\alpha
n}\omega_n},\Set{\delta_n^{(k)}})\> } $ %
}
\settowidth{\mylengthA}{%
$\displaystyle
\frac{x_{n+k+1}-x_n}{x_n+k-1}
\prod_{j=0}^{n-1}
\frac{(x_j+k)(x_{j+k+1}+k-1)}{(x_j+k-1)(x_{j+k+2}+k)}$%
}
\setlength{\LTleft}{0pt}
\setlength{\LTright}{0pt}
{\small
\begin{longtable}[l]{@{}%
b{\mylength}@{}|@{\extracolsep{0.1cm plus 1 fill}}%
>{$\displaystyle}p{5em}<{$}@{}|@{\extracolsep{0.1cm plus 1 fill}}%
b{\mylengthA}@{}}
\caption{Special cases of the \J\
transformation ${}^{\dagger}$\label{tabJtransSpecial}}\\
\hline\hline
{Case}       & \psi_{j}(n) \>^{a)}\phantom{\frac{1}{1}} & $\delta\di{n}{k}
\>^{b)}$
\\
\hline
\endfirsthead
\multicolumn{3}{c}{Table \ref{tabJtransSpecial} (continued)}\\
\hline
{Case}    & \psi_{j,n} \>^{a)}\phantom{\frac{1}{1}} & $\delta\di{n}{k}
\>^{b)}$
\\
\hline
\endhead
\hline\hline
\multicolumn{3}{@{}p{\textwidth}@{}}{${}^{\dagger}${\ Refs.\
\cite{Homeier94ahc,Homeier95,Homeier96Hab}. }
${}^{a)}${\ For the definition of the $\psi_{j,n}$
see Eq.\ (\ref{eqpsijn}). }
${}^{b)}${\ Factors independent of $n$ are
irrelevant.}}
\endlastfoot
{Drummond transformation} \newline $\D\di{n}{k}(\Seq{s_n},\Seq{\omega_n})$ &
n^j
& 1\\ [0.4cm] {}
{Homeier \I\ transformation} \newline
{$\I\di{n}{k}(\alpha,\Seq{s_n},\Seq{\omega_n},\Set{\Delta_n^{(k)}}) $
\newline
$=
{\displaystyle\J\di{n}{2k}(\Seq{s_n},\Seq{e^{-i\alpha
n}\omega_n},\Set{\delta_n^{(k)}})\> }
$}
&
\mbox{Eq.\ (\ref{eqJtranskernel})}
&
{\quad}\newline $\displaystyle
\displaystyle\delta\di{n}{2\ell} = \exp(2i\alpha n)\>, \atop
{\displaystyle\delta\di{n}{2\ell+1}
= \exp(-2i\alpha n)\, \Delta\di{n}{\ell}}$
\\ [0.4cm] {}
{Homeier \F\ transformation} \newline
{$\F\di{n}{k}(\Seq{s_n},\Seq{\omega_n},\Seq{x_n})$}
&
1/(x_n)_j
&
{\quad}\newline $\displaystyle
\displaystyle
\frac{x_{n+k+1}-x_n}{x_n+k-1}
\prod_{j=0}^{n-1}
\frac{(x_j+k)(x_{j+k+1}+k-1)}{(x_j+k-1)(x_{j+k+2}+k)}
$
\\ [0.4cm] {}
{Homeier \pJ\ transformation} \newline
$\pJ\di{n}{k}(\beta,\Seq{s_n},\Seq{\omega_n})$ &
\mbox{Eq.\ (\ref{eqJtranskernel})}
& {\quad}\newline $\displaystyle
\frac{1}
       {(n+\beta+(p-1)k)_2}
$
\\ [0.4cm] {}
{Levin transformation} \newline $\L\di{n}{k}(\beta,\Seq{s_n},\Seq{\omega_n})$ &
(n+\beta)^{-j}
&{\quad}\newline $\displaystyle  \frac{1}
       {(n+\beta)(n+\beta+k+1)}
       $
\\ [0.4cm] {}
{generalized \L\ transformation} \newline
$\L\di{n}{k}(\alpha,\beta,\Seq{s_n},\Seq{\omega_n})$ &
(n+\beta)^{-j\alpha}
& {\quad}\newline $\displaystyle \frac{(n+\beta+k+1)^{\alpha}-(n+\beta)^{\alpha}}
       {(n+\beta)^{\alpha}(n+\beta+k+1)^{\alpha}}
$ \\ [0.4cm] {}
{Levin-Sidi $d\di{}{1}$ transformation} \newline
\cite{LevinSidi81,FordSidi87,Sidi95} \newline
$(d\di{}{1})\di{n}{k}(\alpha,\Seq{s_n})$
             &
(R_{n}+\alpha)^{-j}
& {\quad}\newline $\displaystyle
\frac{1}
       {R_{n+k+1}+\alpha}
-
\frac{1}{R_{n}+\alpha} $
\\ [0.4cm] {}
{Mosig-Michalski algorithm} \newline
\cite{Mosig89,Michalski98} \newline
$M\di{n}{k}(\Seq{s_n},\Seq{\omega_n},\Seq{x_n})$ &
\mbox{Eq.\ (\ref{eqJtranskernel})}
&
{\quad}\newline $\displaystyle
\frac{1}{x_n^2}\left(1-\frac{\omega_n
x_{n+1}^{2k}}{\omega_{n+1}x_{n}^{2k}}\right)
$
\\ [0.4cm] {}
{Sidi W algorithm \newline (GREP$^{(1)}$)}
\cite{Sidi82,Sidi95,Sidi95caf} \newline
$W\di{n}{k}(\Seq{s_n},\Seq{\omega_n},\Seq{t_n})$ &
t_n^j
&
$t_{n+k+1}-t_n$
\\ [0.4cm] {}
{Weniger \C\ transformation } \newline \cite{Weniger92} \newline
$\C\di{n}{k}(\gamma,\beta/\gamma,\Seq{s_n},\Seq{\omega_n})$ &
\frac{1}{(\gamma n+\beta)_{j}}
& {\quad}\newline $\displaystyle
\frac
  {(n+1+(\beta + k-1)/\gamma)_k}
  {(n+(\beta + k)/\gamma)_{k+2}}
$
\\ [0.4cm] {}
{Weniger \M\ transformation} \newline
$\M\di{n}{k}(\xi,\Seq{s_n},\Seq{\omega_n})$ &
\frac{1}{(-n-\xi)_{j}}
&
{\quad}\newline $\displaystyle
\frac
  {(n+1+\xi - (k-1))_k}
  {(n+\xi - k)_{k+2}}
$
\\ [0.4cm] {}
{Weniger \S\ transformation} \newline
$\S\di{n}{k}(\beta,\Seq{s_n},\Seq{\omega_n})$ &
1/(n+\beta)_{j}
&
{\quad}\newline $\displaystyle \frac{1}
       {(n+\beta+2k)_2}
$
\\ [0.4cm] {}
{Iterated Aitken process} \newline \cite{Aitken26,Weniger89}
\newline
${\displaystyle  \mathbf{A}\di{n}{k}(\Seq{s_n})}$ \newline
${\displaystyle  = \J\di{n}{k}(\Seq{s_n},\Seq{\forwarddiff
s_n},\Set{\delta\di{n}{k}})}
$
&
\mbox{Eq.\ (\ref{eqJtranskernel})}
&
{\quad}\newline $\displaystyle
\frac{(\forwarddiff \mathbf{A}\di{n}{k+1}(\Seq{s_n}))\,
                    (\forwarddiff^2 \mathbf{A}\di{n}{k})(\Seq{s_n})}
                   {(\forwarddiff \mathbf{A}\di{n}{k}(\Seq{s_n}))\,
                    (\forwarddiff \mathbf{A}\di{n+1}{k}(\Seq{s_n}))}
$
\\ [0.4cm] {}
{Overholt process} \newline  \cite{Overholt65} \newline
$  {\displaystyle V\di{n}{k}(\Seq{s_n})}$ \newline
$  {\displaystyle {}=\J\di{n}{k}(\Seq{s_n},\Seq{\forwarddiff
  s_n},\Set{\delta\di{n}{k}})}
$
&
\mbox{Eq.\ (\ref{eqJtranskernel})}
&
{\quad}\newline $\displaystyle
\frac{(\forwarddiff s_{n+k+1}) \forwarddiff\left[(\forwarddiff
  s_{n+k})^{k+1}\right]}
         {(\forwarddiff s_{n+k})^{k+1}}
$
\\ [0.4cm]
\end{longtable}
}


      The \J\ transformation may also be computed using the
      alternative recursive schemes \cite{Homeier94ahc,Homeier98ots}
      \begin{equation} \label{eqJtransRecA}
         \begin{array}{>{\displaystyle}r@{}>{\displaystyle}l}
           &\widehat{D}\di{n}{0} = 1/\omega_n\>,
           \qquad \widehat{N}\di{n}{0} = s_n/\omega_n\>,
           \\
           &\widehat{D}\di{n}{k} =
           \Phi\di{n}{k-1}   \widehat{D}\di{n+1}{k-1}
           -\widehat{D}\di{n}{k-1}\>,\qquad k\in\Nset\>,
           \\
           &\widehat{N}\di{n}{k} =
           \Phi\di{n}{k-1}   \widehat{N}\di{n+1}{k-1}
           -\widehat{N}\di{n}{k-1}\>,\qquad k\in\Nset\>,
           \\
           &\J\di{n}{k}(\Seq{s_n},\Seq{\omega_n},\Set{\delta\di{n}{k}})
              = \frac{\widehat{N}\di{n}{k}}{\widehat{D}\di{n}{k}}
          \end{array}
      \end{equation}
      with
      \begin{equation}\label{eqPhink}
         \Phi\di{n}{0}=1\>,\qquad \Phi\di{n}{k} =  \frac{\delta\di{n}{0} \delta\di{n}{1} \cdots\delta\di{n}{k-1}}
        {\delta\di{n+1}{0} \delta\di{n+1}{1}
        \cdots\delta\di{n+1}{k-1}}\>,\qquad k\in\Nset\>.
      \end{equation}
      and
      \begin{equation} \label{eqJtransRecB}
         \begin{array}{l}
           \widetilde{D}\di{n}{0} \displaystyle{}= 1/\omega_n\>, \qquad \widetilde{N}\di{n}{0} = s_n/\omega_n\>,
           \\
           \widetilde{D}\di{n}{k}  \displaystyle{}=
            \widetilde{D}\di{n+1}{k-1}
               -\Psi\di{n}{k-1}  \widetilde{D}\di{n}{k-1}\>,\qquad
               k\in\Nset\>,
           \\
           \widetilde{N}\di{n}{k}  \displaystyle{}=
           \widetilde{N}\di{n+1}{k-1}
              -\Psi\di{n}{k-1}   \widetilde{N}\di{n}{k-1}\>,\qquad
              k\in\Nset\>,
           \\
           \displaystyle           \J\di{n}{k}(\Seq{s_n},\Seq{\omega_n},\Set{\delta\di{n}{k}})
           = \frac{\widetilde{N}\di{n}{k}}{\widetilde{D}\di{n}{k}}
         \end{array}
      \end{equation}
      with
      \begin{equation}
         \Psi\di{n}{0}=1\>,\qquad \Psi\di{n}{k} =  \frac{\delta\di{n+k}{0}
           \delta\di{n+k-1}{1} \cdots\delta\di{n+1}{k-1}}
           {\delta\di{n+k-1}{0} \delta\di{n+k-2}{1}
           \cdots\delta\di{n}{k-1}}\>,\qquad k\in\Nset\>.
      \end{equation}
      The  quantities $\Psi\di{n}{k}$ should not be mixed up
      with the $\Psi_{k,n}(u)$ as defined in Eq.\
      (\ref{eqbekaverf_Psi}).

      As shown in \cite{Homeier98ots}, the coefficients for the
      algorithm (\ref{eqJtransRecA}) that are defined via $\widehat D\di{n}{k}=\sum_{j=0}^{k} \lambda\di{n,j}{k}
      /\omega_{n+j}$, satisfy the recursion
      \begin{equation}\label{eqreclamda}
        \lambda\di{n,j}{k+1}=
        \Phi\di{n}{k}\lambda\di{n+1,j-1}{k}-\lambda\di{n,j}{k}\>
      \end{equation}
      with starting values $\lambda\di{n,j}{0}=1$. This holds for all
      $j$ if we define $\lambda\di{n,j}{k}=0$ for $j<0$ or $j>k$.
      Because $\Phi\di{n}{k}\ne 0$, we have $\lambda\di{n,k}{k}\ne 0$
      such that $\Set{\lambda\di{n,j}{k}}$ is a coefficient set for all
      $k\in\Nset_0$.

      Similarly, the coefficients for the
      algorithm (\ref{eqJtransRecB}) that are defined via $\widetilde D\di{n}{k}=\sum_{j=0}^{k} \widetilde\lambda\di{n,j}{k}
      /\omega_{n+j}$, satisfy the recursion
      \begin{equation}\label{eqreclamdaB}
        \widetilde\lambda\di{n,j}{k+1}=
        \widetilde\lambda\di{n+1,j-1}{k}-
        \Psi\di{n}{k}\widetilde\lambda\di{n,j}{k}\>
      \end{equation}
      with starting values $\widetilde\lambda\di{n,j}{0}=1$. This holds
      for all $j$ if we define $\widetilde\lambda\di{n,j}{k}=0$ for
      $j<0$ or $j>k$. In this case, we have
      $\widetilde\lambda\di{n,k}{k}=1$   such that
      $\Set{\widetilde\lambda\di{n,j}{k}}$ is a coefficient set for all
      $k\in\Nset_0$.

      Since the \J\ transformation vanishes for
      $\Seq{s_n}=\Seq{c\omega_n}$, $c\in \Kset$ according to Eq.\
      (\ref{eqjtrans_31}) for all $k\in\Nset$, it is convex.
      This may also be shown by using induction in $k$ using
      $\lambda\di{n,1}{1}=-\lambda\di{n,0}{1}=1$ and the
      equation
       \begin{equation} \label{eqjtranssumlambda}
         \sum_{j=0}^{k+1} \lambda\di{n,j}{k+1} =
         \Phi\di{n}{k}\sum_{j=0}^{k} \lambda\di{n+1,j}{k} -
         \sum_{j=0}^{k} \lambda\di{n,j}{k}
       \end{equation}
      that follows from Eq.\ (\ref{eqreclamda}).

      Assuming that the limits $\Phi_k=\lim_{n\to\infty}\Phi\di{n}{k}$
      exist for all $k\in\Nset$ and noting that for $k=0$ always
      $\Phi_0=1$ holds, it follows that there exists a limiting
      transformation $\Ringel\J[\Ringel\Lambda]$ that can be considered
      as special variant of the $\mathcal{J}$ transformation and with
      coefficients given explicitly as \cite[Eq.\ (16)]{Homeier98ots}
      \begin{equation}\label{explicit}
        \Ringel\lambda\di{j}{k} = (-1)^{k-j} \sum_{j_0+j_1+\dots+j_{k-1}=j,
        \atop j_0\in\{0,1\},\dots,j_{k-1}\in\{0,1\}} \prod_{m=0}^{k-1} (\Phi_m)^{j_m}\>.
      \end{equation}
      As characteristic polynomial we obtain
       \begin{equation} \label{eqJtransPi}
         \Ringel\Pi\di{}{k}(z) = \sum_{j=0}^k \Ringel\lambda\di{j}{k} z^j =
         \prod_{j=0}^{k-1} (\Phi_j \,z-1)\>.
       \end{equation}
       Hence, the $\Ringel\J$ transformation is convex since
       $\Ringel\Pi\di{}{k}(1)=0$ due to $\Phi_0=1$.

       \paragraph{The \pJ\ Transformation.} This is the special case
         of the \J transformation corresponding to
         \begin{equation}
           \delta\di{n}{k} = \frac{1}{\bigl(n+\beta+(p-1)k\bigr)_2}\>,
         \end{equation}
         or to \cite[Eq.\ (18)]{Homeier98ots}\footnote{The equation in
         \cite{Homeier98ots} contains an error.}
         \begin{equation}
           \Phi\di{n}{k} = \left\{
           \begin{array}{ll}
           \displaystyle
           \left.
            \left(
              \frac{n+\beta+2}{p-1}
               \right)_k
               \right/
               \left(
                 \frac{n+\beta}{p-1}
                 \right)_k &
                 \mbox{for }p\ne1\\
                  & \\
           \displaystyle
           \left(\frac{n+\beta+2}{n+\beta}\right)^k &
           \mbox{for }p=1
           \end{array}
           \right.
         \end{equation}
         or to
         \begin{equation}
           \Psi\di{n}{k} = \left\{
           \begin{array}{ll}
           \displaystyle
           \left.
            \left(
              \frac{n+\beta+k-1}{p-2}
               \right)_k
               \right/
               \left(
                 \frac{n+\beta+k+1}{p-2}
                 \right)_k &
                 \mbox{for }p\ne2\\
                  & \\
           \displaystyle
           \left(\frac{n+\beta+k-1}{n+\beta+k+1}\right)^k &
           \mbox{for }p=2
           \end{array}
           \right.
           \end{equation}
           that is,
           \begin{equation}\label{eqjtrans_20}
             \pJ\di{n}{k} (\beta,\Seq{s_n},\Seq{\omega_n}) =
             \J\di{n}{k} (\Seq{s_n},\Seq{\omega_n},\Set{1/(n+\beta+(p-1) k)_2})\>.
           \end{equation}

           The limiting transformation \RingelpJ\ of the \pJ\ transformation
           exists for all $p$ and corresponds to the $\Ringel\J$ transformation
           with $\Phi_k=1$ for all $k$ in $\Nset_0$. This is exactly the
           Drummond transformation discussed in Section
           \ref{secDrummond}, i.e., we have
           \begin{equation} \label{eqRingelpJ}
             \RingelpJ\di{n}{k} (\beta,\Seq{s_n},\Seq{\omega_n}) =
             \D\di{n}{k} (\Seq{s_n},\Seq{\omega_n})\>.
           \end{equation}

    \subsubsection{Drummond Transformation}\label{secDrummond}
       This transformation was given by Drummond \cite{Drummond72}. It
       was also discussed by Weniger \cite{Weniger89}. It
       may be defined as
       \begin{equation} \label{eqDrummond}
         \D\di{n}{k}(\Seq{s_n},\Seq{\omega_n}) =
                       \frac{\forwarddiff^k [s_n/\omega_n]}
                            {\forwarddiff^k [1/\omega_n]}
       \end{equation}
       Using the definition (\ref{eq28}) of the forward difference
       operator, the coefficients may be taken as
       \begin{equation} \label{eqDrummondLambda}
          \lambda\di{n,j}{k} = (-1)^{j} \binom{k}{j}\>,
       \end{equation}
       i.e., independent of $n$. As moduli, one has
       $\mu\di{n}{k}=\binom{k}{\Ent{k/2}}=\widetilde\mu\di{}{k}$. Consequently, the
       Drummond
       transformation is given in subnormalized form.
       As characteristic polynomial we obtain
       \begin{equation} \label{eqDrummondPi}
         \Pi\di{n}{k}(z) = \sum_{j=0}^k (-1)^{j} \binom{k}{j} z^j =
         (1-z)^k\>.
       \end{equation}
       Hence, the Drummond transformation is convex since
       $\Pi\di{n}{k}(1)=0$.
       Interestingly, the Drummond transformation is identical to its
       limiting transformation:
       \begin{equation} \label{eqDrummondRingel}
         \Ringel\D\di{}{k}(\Seq{s_n},\Seq{\omega_n}) =
         \D\di{n}{k}(\Seq{s_n},\Seq{\omega_n})\>.
       \end{equation}
       The Drummond transformation may be computed using the recursive
       scheme
       \begin{equation} \label{eqDrummondRec}
          \begin{array}{>{\displaystyle}r@{}>{\displaystyle}l}
            & N\di{n}{0} = s_n/\omega_n\>,\quad D\di{n}{0} = 1/\omega_n\>,\\
            & N\di{n}{k} = \forwarddiff N\di{n}{k-1}\>,\quad
              D\di{n}{k} = \forwarddiff D\di{n}{k-1}\>,\\
            & \D\di{n}{k} = N\di{n}{k} / D\di{n}{k}\>.
           \end{array}
       \end{equation}

    \subsubsection{Levin Transformation}
       This transformation was given by Levin \cite{Levin73}. It
       was also discussed by Weniger \cite{Weniger89}. It
       may be defined as\footnote{Note that the order of indices is
       different from that in the literature.}
       \begin{equation} \label{eqLevin}
         \L\di{n}{k}(\beta,\Seq{s_n},\Seq{\omega_n}) =
                       \frac{(n+\beta+k)^{1-k}\forwarddiff^k [(n+\beta)^{k-1} s_n/\omega_n]}
                            {(n+\beta+k)^{1-k}\forwarddiff^k [(n+\beta)^{k-1}/\omega_n]}
       \end{equation}
       Using the definition (\ref{eq28}) of the forward difference
       operator, the coefficients may be taken as
       \begin{equation} \label{eqLevinLambda}
          \lambda\di{n,j}{k} = (-1)^{j} \binom{k}{j}
          (n+\beta+j)^{k-1}/(n+\beta+k)^{k-1}\>,
       \end{equation}
       The moduli satisfy $\mu\di{n}{k} \le
       \binom{k}{\Ent{k/2}}=\widetilde\mu\di{}{k}$ for given $k$. Consequently, the Levin
       transformation is given in subnormalized form.
       As characteristic polynomial we obtain
       \begin{equation} \label{eqLevinPi}
         \Pi\di{n}{k}(z) = \sum_{j=0}^k (-1)^{j} \binom{k}{j} z^j (n+\beta+j)^{k-1}/(n+\beta+k)^{k-1}
       \end{equation}
       Since $\Pi\di{n}{k}(1)=0$ because $\forwarddiff^{k-1}$
       annihilates any polynomial in $n$ with degree less than $k$, the
       Levin transformation is convex.
       The limiting transformation is identical to the Drummond
       transformation
       \begin{equation} \label{eqLevinRingel}
         \Ringel\L\di{}{k}(\Seq{s_n},\Seq{\omega_n}) =
         \D\di{n}{k}(\Seq{s_n},\Seq{\omega_n})\>.
       \end{equation}

       The Levin transformation may be computed using the recursive
       scheme
       \cite{Longman81,FesslerFordSmith83a,Weniger89}, \cite[Sec.\ 2.7]{BrezinskiRedivoZaglia91}
       \begin{equation} \label{eqLevinRec}
          \begin{array}{>{\displaystyle}r@{}>{\displaystyle}l}
            & N\di{n}{0} = s_n/\omega_n\>,\quad D\di{n}{0} = 1/\omega_n\>,\\
            & N\di{n}{k} =  N\di{n+1}{k-1} -
            \frac{(\beta+n)(\beta+n+k-1)^{k-2}}{(\beta+n+k)^{k-1}}
            N\di{n}{k-1}\>,\\
            & D\di{n}{k} =  D\di{n+1}{k-1} -
            \frac{(\beta+n)(\beta+n+k-1)^{k-2}}{(\beta+n+k)^{k-1}}
            D\di{n}{k-1}\>,\\
            & \L\di{n}{k}(\beta,\Seq{s_n},\Seq{\omega_n}) = N\di{n}{k} / D\di{n}{k}\>.
           \end{array}
       \end{equation}
       This is essentially the same as the recursive scheme (\ref{eqJtransRecB})
       for the \J\ transformation
       with
       \begin{equation}
          \Psi\di{n}{k}=\frac{(\beta+n)(\beta+n+k)^{k-1}}{(\beta+n+k+1)^{k}}\>.
       \end{equation}
       since the Levin transformation is a special case of the \J\
       transformation (see Table \ref{tabJtransSpecial}). Thus, the
       Levin transformation can also be computed recursively using the
       scheme (\ref{eqJtransRec})
       \begin{equation}
         \delta\di{n}{k} = \frac{1}{(n+\beta)(n+\beta+k+1)}\>,
       \end{equation}
       or the scheme (\ref{eqJtransRecA}) with \cite{Homeier98ots}
       \begin{equation}
         \Phi\di{n}{k}=(n+\beta+k+1)\,\frac{(n+\beta+1)^{k-1}}{(n+\beta)^k}\>.
       \end{equation}

    \subsubsection{Weniger Transformations}
       Weniger \cite{Weniger89,Weniger92,Weniger94} derived sequence
       transformations related to factorial series. These may be
       regarded as special cases of the transformation
       \begin{equation} \label{eqWenigerC}
         \C\di{n}{k}(\alpha,\zeta,\Seq{s_n},\Seq{\omega_n}) =
             \frac{((\alpha[n+\zeta+k])_{k-1})^{-1}
                    \forwarddiff^k [(\alpha[n+\zeta])_{k-1} s_n/\omega_n]}
                  {((\alpha[n+\zeta+k])_{k-1})^{-1}
                    \forwarddiff^k
                    [(\alpha[n+\zeta])_{k-1}/\omega_n]}\>.
       \end{equation}
       In particular, the Weniger \S\ transformation may be defined as
       \begin{equation} \label{eqWenigerS}
         \S\di{n}{k}(\beta,\Seq{s_n},\Seq{\omega_n}) =
         \C\di{n}{k}(1,\beta,\Seq{s_n},\Seq{\omega_n})
       \end{equation}
       and the Weniger \M\ transformation as
       \begin{equation} \label{eqWenigerM}
         \M\di{n}{k}(\xi,\Seq{s_n},\Seq{\omega_n}) =
         \C\di{n}{k}(-1,\xi,\Seq{s_n},\Seq{\omega_n})\>.
       \end{equation}
       The parameters $\beta$, $\xi$, and $\zeta$ are taken to be positive real
       numbers. Weniger
       considered the \C\ transformation only for $\alpha>0$
       \cite{Weniger92,Weniger94} and thus, he was not considering the
       \M\ transformation as a special case of the \C\ transformation.
       He also found that one should choose $\xi \ge k-1$. In the $u$
       variant of the \M\ transformation he proposed to choose
       $\omega_n=(-n-\xi) \forwarddiff s_{n-1}$. This variant is
       denoted as ${}^{u}\M$ transformation in the present work.

       Using the definition (\ref{eq28}) of the forward difference
       operator, the coefficients may be taken as
       \begin{equation} \label{eqWenigerCLambda}
          \lambda\di{n,j}{k} = (-1)^{j} \binom{k}{j}
          (\alpha[n+\zeta+j])_{k-1}/(\alpha[n+\zeta+k])_{k-1}\>
       \end{equation}
       in the case of the \C\ transformation,        as
       \begin{equation} \label{eqWenigerSLambda}
          \lambda\di{n,j}{k} = (-1)^{j} \binom{k}{j}
          (n+\beta+j)_{k-1}/(n+\beta+k)_{k-1}\>
       \end{equation}
       in the case of the \S\ transformation,     and as
       \begin{equation} \label{eqWenigerMLambda}
          \lambda\di{n,j}{k} = (-1)^{j} \binom{k}{j}
          (-n-\xi-j)_{k-1}/(-n-\xi-k)_{k-1}\>
       \end{equation}
       in the case of the \M\ transformation.

       The \S\ transformation in (\ref{eqWenigerS}) may be computed
       using the recursive scheme (\ref{eqJtransRecB}) with \cite[Sec.
       8.3]{Weniger89} \begin{equation} \label{eqPsiWenigerS}
       \Psi\di{n}{k} =
       \frac{(\beta+n+k)(\beta+n+k-1)}{(\beta+n+2k)(\beta+n+2k-1)}
       \end{equation} The \M\ transformation in (\ref{eqWenigerM}) may
       be computed using the recursive scheme (\ref{eqJtransRecB}) with
       \cite[Sec. 9.3]{Weniger89} \begin{equation}
       \label{eqPsiWenigerM} \Psi\di{n}{k} =
       \frac{\xi+n-k+1}{\xi+n+k+1} \end{equation} The \C\
       transformation in (\ref{eqWenigerC}) may be computed using the
       recursive scheme (\ref{eqJtransRecB}) with \cite[Eq.\
       (3.3)]{Weniger92}
       \begin{equation} \label{eqPsiWenigerC}
         \Psi\di{n}{k} = (\alpha[\zeta+n]+k-2)
         \frac{(\alpha[n+\zeta+k-1])_{k-2}}{(\alpha[n+\zeta+k])_{k-1}}
       \end{equation}

       Since the operator
       $\forwarddiff^k$ for $k\in\Nset$ annihilates all polynomials in
       $n$ of degree smaller than $k$, the transformations \S, \M, and
       \C\ are convex. The moduli satisfy $\mu\di{n}{k} \le
       \binom{k}{\Ent{k/2}}=\widetilde\mu\di{}{k}$ for given $k$.
       Consequently, the three Weniger transformations are given in
       subnormalized form.

       For $\alpha\to\infty$, the Levin transformation is obtained from
       the \C\ transformation \cite{Weniger92}. The \S\ transformation
       is identical to the $\PJ{3}$ transformation. It is also
       the special case $x_n=n+\beta$ of the \F\ transformation.
       Analogously, the \C\ transformation is obtained for
       $x_n=\alpha[\zeta+n]$. All
       these Weniger transformations are special cases of the \J\
       transformation (compare Table \ref{tabJtransSpecial}).

       The limiting transformation of all these Weniger transformations
       is the Drummond transformation.

    \subsubsection{Levin-Sidi Transformations and W Algorithms}
       As noted above in Section \ref{appdm}, 
       the { $d^{(m)}$ transformations} were introduced by Levin and Sidi
       \cite{LevinSidi81} as a generalization of the $u$ variant of the
       Levin transformation, and these transformations may be 
       implemented recursively using the
       ${\mathbf W}^{(m)}$ algorithms.  

       The case $m=1$ corresponding to the $d\di{}{1}$ transformation and the
       ${\mathbf W}^{(1)}=W$ algorithm is relevant for the present survey of Levin-type
       transformations. In the following, the $k$-th order
       transformation $\T\di{}{k}$ of Levin-type transformation
       \T as given by the $W$ algorithm is denoted by
       $W\di{}{k}$ which should not be confused with the
       ${\mathbf W}^{(m)}$ algorithms of Ford and Sidi \cite{FordSidi87}.

       The $W$ algorithm \cite{Sidi82} was also studied by other
       authors \cite[Sec.\ 7.4]{Weniger89}, \cite[p. 71f,
       116f]{BrezinskiRedivoZaglia91} and may be regarded as a special case of
       the \J\ transformation \cite{Homeier94ahc}. It may be defined as
       (compare \cite[Theorems 1.1 and 1.2]{Sidi95caf})
       \begin{equation} \label{eqWRec}
          \begin{array}{>{\displaystyle}r@{}>{\displaystyle}l}
            & N\di{n}{0} = \frac{s_n}{\omega_n}\>,\quad
              D\di{n}{0} = \frac{1}{\omega_n}\>,\\
            & N\di{n}{k} = \frac{N\di{n+1}{k-1} - N\di{n}{k-1}}
                                  {t_{n+k}-t_{n}}\>,\\
            & D\di{n}{k} = \frac{D\di{n+1}{k-1}-D\di{n}{k-1}}
                                  {t_{n+k}-t_{n}}\>,\\
            & W\di{n}{k}(\Seq{s_{n}},
             \Seq{\omega_n},
             \Seq{t_n}) = N\di{n}{k} / D\di{n}{k}\>,
           \end{array}
       \end{equation}
       and computes
       \begin{equation} \label{eqW}
          W\di{n}{k}(\Seq{s_n},\Seq{\omega_n},\Seq{t_n}) =
                       \frac{\square\di{n}{k} ( s_n/\omega_n)}
                            {\square\di{n}{k} (1/\omega_n)}
       \end{equation}
       where the divided difference operators
       $\square\di{n}{k}=\square\di{n}{k}[\Seq{t_n}]$ are used. The $W$
       algorithm may be used to calculate the Levin transformation on
       putting $t_n=1/(n+\beta)$. Some authors call a linear variant of
       the $W$ algorithm with $\omega_n=(-1)^{n+1} e^{-n\zeta
       q}t_n^\alpha$ the $W$ transformation, while the $\widetilde t$
       variant of the $W$ algorithm \cite{Sidi88a,Sidi88c} is sometimes
       called $mW$ transformation
       \cite{LucasStone95,HasegawaSidi96,Michalski98}.

       If $t_{n+1}/t_{n}\to \tau$ for large $n$, one obtains as
       limiting transformation the $\Ringel\J$ transformation
       with $\Phi_j=\tau^{-j}$ and characteristic polynomial
       \begin{equation} \label{eqRingelPid1}
         \Ringel\Pi\di{}{k}(z) =
         \prod_{j=0}^{k-1} (z/\tau^{j}-1)\>.
       \end{equation}

       For the $d\di{}{1}$ transformation, we write
       \begin{equation} \label{eqd1}
       (d\di{}{1})\di{n}{k}(\alpha,\Seq{s_n},\Seq{\xi_n})
       =W\di{n}{k}(\Seq{s_{\xi_{n}}},
             \Seq{(\xi_{n}+\alpha)(s_{\xi_{n}}-s_{\xi_{n}-1})},
             \Seq{1/(\xi_{n}+\alpha)})
       \end{equation}
       Thus, it corresponds to the variant of the $W$ algorithm with
       remainder estimates chosen as $(\xi_{n}+\alpha)(s_{\xi_{n}}-s_{\xi_{n}-1})$ operating on the
       subsequence $\Seq{s_{\xi_n}}$ of $\Seq{s_n}$ with
       $t_n=1/(\xi_n+\alpha)$. It should be noted that this is not(!)
       identical to the $u$ variant
       \begin{equation} \label{eqWu}
         {}^uW\di{n}{k}(\Seq{s_{\xi_n}},\Seq{1/(\xi_{n}+\alpha)})
           =W\di{n}{k}(\Seq{s_{\xi_{n}}},\Seq{{}^u\omega_n},
             \Seq{1/(\xi_{n}+\alpha)})\>,
       \end{equation}
       neither for
       ${}^u\omega_n={(n+\alpha)(s_{\xi_{n}}-s_{\xi_{n-1}})}$ nor
       for ${}^u\omega_n={(\xi_{n}+\alpha)(s_{\xi_{n}}-s_{\xi_{n-1}})}$,
       since the remainder estimates are chosen differently in Eq.\
       (\ref{eqd1}).

       The $d\di{}{1}$ transformation was
       thoroughly analyzed by Sidi (See \cite{Sidi95,Sidi95caf} and
       references therein).

    \subsubsection{Mosig-Michalski Transformation}
       The Mosig-Michalski transformation --- also known as
       ``weighted--averages algorithm'' --- was introduced by Mosig
       \cite{Mosig89} and  modified later by Michalski who gave the
       $\widetilde t$ variant of the transformation the name $\mathcal{K}$
       transformation (that is used for a different transformation in
       the present article(!)), and applied it to the computation of
       Sommerfeld integrals \cite{Michalski98}.

       The Mosig-Michalski transformation $M$ may be defined via the
       recursive scheme
       \begin{equation} \label{eqMMrec}
          \begin{array}{@{}>{\displaystyle}l}
            s\di{n}{0} = s_n\>, \\
            s\di{n}{k+1} =
              \frac{s\di{n}{k}+\eta\di{n}{k}s\di{n+1}{k}}
                   {1+\eta\di{n}{k}}\>, \\
            M\di{n}{k}(\Seq{s_n},\Seq{\omega_n},\Seq{x_n}) = s\di{n}{k}
          \end{array}
       \end{equation}
       for $n\in\Nset_0$ and $k\in\Nset_0$ where $\Seq{x_n}$ is an
       auxiliary sequence with $\lim_{n\to\infty} 1/x_{n}=0$ such that
       $x_{n+\ell}>x_{n}$ for $\ell \in \Nset$ and $x_0>1$, i.e., a
       diverging sequence  of monotonously increasing positive numbers,
       and
       \begin{equation} \label{eqetank}
                \eta\di{n}{k}=-\frac{\omega_n}{\omega_{n+1}}
         \left(\frac{x_{n+1}}{x_n}\right)^{2k}\>.
       \end{equation}
       Putting $\omega\di{n}{k}=\omega_n/x_n^{2k}$,
       $N\di{n}{k}=s\di{n}{k}/\omega\di{n}{k}$, and
       $D\di{n}{k}=1/\omega\di{n}{k}$, it is easily seen that the
       recursive scheme (\ref{eqMMrec}) is equivalent to the scheme
       (\ref{eqJtransRec}) with
       \begin{equation} \label{eqMMdelta}
         \delta\di{n}{k} = \frac{1}{x_n^2}
                           \left(
                              1-\frac{\omega_n x_{n+1}^{2k}}
                                     {\omega_{n+1}x_{n}^{2k}}
                           \right)\>.
       \end{equation}
       Thus, the Mosig-Michalski transformation is a special case of
       the \J\ transformation. Its character as a Levin-type
       transformation is somewhat formal since the $\delta\di{n}{k}$
       and, hence, the coefficients $\lambda\di{n,j}{k}$ depend on the
       $\omega_n$.

       If $x_{n+1}/x_n\sim
       \xi>1$ for large $n$, then a limiting transformation
       exists, namely \\
       $M(\Seq{s_n},\Seq{\omega_n},\Seq{\xi^{n+1}})$.
       It corresponds to the $\Ringel\J$ transformation with
       $\Phi_k=\xi^{2k}$. This may be seen by putting $\widehat
       D\di{n}{k}=1/\omega_n$, $\widehat
       N\di{n}{k}=s\di{n}{k}D\di{n}{k}$ and
       $\Phi\di{n}{k}=\xi^{2k}$ in Eq.\ (\ref{eqJtransRecA}).

     \subsubsection{\F\ Transformation}
       This transformation is seemingly new. It will be derived in a
       later section.
       It
       may be defined as
       \begin{equation} \label{eqFtrans}
         \F\di{n}{k}(\Seq{s_n},\Seq{\omega_n},\Seq{x_n}) =
                       \frac{\square\di{n}{k} ( (x_n)_{k-1} s_n/\omega_n)}
                            {\square\di{n}{k} ((x_n)_{k-1}/\omega_n)}
              =         \frac{x_n^k/(x_n)_{k-1}\square\di{n}{k} ( (x_n)_{k-1} s_n/\omega_n)}
                             {x_n^k/(x_n)_{k-1}\square\di{n}{k} ((x_n)_{k-1}/\omega_n)}
       \end{equation}
       where $\Seq{x_n}$ is an auxiliary sequence with $\lim_{n\to\infty} 1/x_{n}=0$ such that
       $x_{n+\ell}>x_{n}$ for $\ell \in \Nset$ and $x_0>1$, i.e., a
       diverging sequence  of monotonously
       increasing positive numbers.
       Using the definition (\ref{eq30}) of the divided difference
       operator $\square\di{n}{k}=\square\di{n}{k}[\Seq{x_n}]$, the
       coefficients may be taken as
       \begin{equation} \label{eqFtransLambda}
          \lambda\di{n,j}{k} = \frac{(x_{n+j})_{k-1}}
                                    {(x_n)_{k-1}}
                              \prod_{i=0 \atop i\ne j }^{k}
                                   \frac{x_n}
                                        {x_{n+j}-x_{n+i}}
                             = \prod_{m=0}^{k-2}
                               \frac{x_{n+j}+m}
                                    {x_n+m}
                                   \left(\frac{x_n}
                                        {x_{n+j}}\right)^k
                              \prod_{i=0 \atop i\ne j }^{k}
                                   \frac{1}
                                        {1- x_{n+i}/x_{n+j}}\>.
       \end{equation}
       Assuming that the following limit exists such that
       \begin{equation} \label{eqlimxn}
         \lim_{n\to\infty} \frac{x_{n+1}}{x_n} = \xi >1
       \end{equation}
       holds, we see that one can define a limiting transformation
       $\Ringel\F\di{}{k}$ with coefficients
       \begin{equation} \label{eqRingelFtransLambda}
         \Ringel\lambda\di{j}{k}=\lim_{n\to\infty} \lambda\di{n,j}{k} =
         \frac{1}{\xi^{j}} \prod_{\ell=0 \atop \ell\ne j }^{k}
         \frac{1}{1-\xi^{\ell-j}} = (-1)^k\xi^{-k(k+1)/2} \prod_{\ell=0 \atop \ell\ne j }^{k}
         \frac{1}{\xi^{-j}-\xi^{-\ell}}
       \end{equation}
       since
       \begin{equation} \label{eqaux1}
         \prod_{m=0}^{k-2}
             \frac{x_{n+j}+m}
                  {x_n+m}
         \left(
           \frac{x_n}
                {x_{n+j}}
         \right)^k
         \prod_{\ell=0 \atop \ell\ne j }^{k}
             \frac{1}
                  {1- x_{n+\ell}/x_{n+j}}
         \to
         \xi^{(k-1)j}  \xi^{k(-j)}
                              \prod_{\ell=0 \atop \ell\ne j }^{k}
                                   \frac{\xi^{-l}}
                                        {\xi^{-\ell}- \xi^{-j}}
       \end{equation}
       for $n\to \infty$.
       Thus, the limiting transformation is given by
       \begin{equation} \label{eqFtransRingel}
         \Ringel\F\di{}{k}(\Seq{s_n},\Seq{\omega_n},\xi) =
         \frac{\displaystyle
               \sum_{j=0}^k
                 \frac{s_{n+j}}
                      {\omega_{n+j}}
                 \prod_{\ell=0 \atop \ell\ne j }^{k}
                    \frac{1}
                         {\xi^{-j}-\xi^{-\ell}}
               }
              {\displaystyle
               \sum_{j=0}^k
                 \frac{1}
                      {\omega_{n+j}}
                 \prod_{\ell=0 \atop \ell\ne j }^{k}
                    \frac{1}
                         {\xi^{-j}-\xi^{-\ell}}
               }
         \>.
       \end{equation}
       Comparison with the definition (\ref{eq30}) of the divided difference
       operators reveals that the limiting transformation can be
       rewritten as
       \begin{equation} \label{eqRingelFtrans}
         \Ringel\F\di{}{k}(\Seq{s_n},\Seq{\omega_n},\xi) =
          \frac{\square\di{n}{k}[\Seq{\xi^{-n}}] ( s_n/\omega_n)}
               {\square\di{n}{k}[\Seq{\xi^{-n}}]
               (1/\omega_n)}\>.
       \end{equation}
       Comparison to Eq.\ (\ref{eqW}) shows that  the limiting
       transformation is nothing but the W algorithm for
       $t_n=\xi^{-n}$.
       As characteristic polynomial we obtain
       \begin{equation} \label{eqFtransPi}
         \Ringel\Pi\di{}{k}(z) =
            \sum_{j=0}^k
              {z}^{j}
              \prod_{\ell=0 \atop \ell\ne j }^{k}
                 \frac{1}
                      {\xi^{-j}-\xi^{-\ell}} = \xi^{k(k+1)/2}
                      \prod_{j=0}^{k-1}\frac{1-z\xi^j}{\xi^{j+1}-1}
         \>.
       \end{equation}
       The last equality is easily proved by induction.
       Hence, the $\Ringel\F$ transformation is convex since
       $\Ringel\Pi\di{}{k}(1)=0$.

       As shown in Appendix \ref{appB}, the \F\ transformation may be computed using the recursive
       scheme
       \begin{equation} \label{eqFtransRec}
          \begin{array}{>{\displaystyle}r@{}>{\displaystyle}l}
            & N\di{n}{0} = \frac{1}{x_n-1}\frac{s_n}{\omega_n}\>,\quad
              D\di{n}{0} = \frac{1}{x_n-1}\frac{1}{\omega_n}\>,\\
            & N\di{n}{k} = \frac{(x_{n+k}+k-2)N\di{n+1}{k-1}
                                  -(x_{n}+k-2)N\di{n}{k-1}}
                                  {x_{n+k}-x_{n}}\>,\\
            & D\di{n}{k} = \frac{(x_{n+k}+k-2)D\di{n+1}{k-1}
                                  -(x_{n}+k-2)D\di{n}{k-1}}
                                  {x_{n+k}-x_{n}}\>,\\
            & \F\di{n}{k} = N\di{n}{k} / D\di{n}{k}\>.
           \end{array}
       \end{equation}

       It follows directly from Eq.\ (\ref{eqRingelFtrans}) and the recursion
       relation for divided differences that the limiting
       transformation can be computed via the recursive scheme
       \begin{equation} \label{eqRingelFtransRec}
          \begin{array}{>{\displaystyle}r@{}>{\displaystyle}l}
            & \Ringel N\di{n}{0} = \frac{s_n}{\omega_n}\>,\quad
              \Ringel D\di{n}{0} = \frac{1}{\omega_n}\>,\\
            & \Ringel N\di{n}{k} = \frac{\Ringel N\di{n+1}{k-1}
                                  -\Ringel N\di{n}{k-1}}
                                  {\xi^{-(n+k)}-\xi^{-n}}\>,\\
            & \Ringel D\di{n}{k} = \frac{\Ringel D\di{n+1}{k-1}
                                  -\Ringel D\di{n}{k-1}}
                                  {\xi^{-(n+k)}-\xi^{-n}}\>,\\
            & \Ringel\F\di{n}{k} = \Ringel N\di{n}{k} /\Ringel D\di{n}{k}\>.
           \end{array}
       \end{equation}

    \subsubsection{\JD\ Transformation}
      This transformation is newly introduced in this article. In
      Section \ref{secJD}, it is
      derived via (asymptotically) hierarchically consistent iteration
      of the $\D\di{}{2}$ transformation, i.e., of
       \begin{equation} \label{eqDrummond2}
         s_n'  = \frac{\forwarddiff^2 (s_n/\omega_n)}
                      {\forwarddiff^2 (1/\omega_n)}
                      \>.
       \end{equation}
      The \JD\ transformation may be defined via the recursive scheme
       \begin{equation} \label{eqJDRec}
          \begin{array}{>{\displaystyle}r@{}>{\displaystyle}l}
            & N\di{n}{0} = s_n/\omega_n\>,\quad D\di{n}{0} = 1/\omega_n\>,\\
            & N\di{n}{k} = \widetilde\nabla\di{n}{k-1} N\di{n}{k-1}\>,\quad
              D\di{n}{k} = \widetilde\nabla\di{n}{k-1} D\di{n}{k-1}\>,\\
            & \JD\di{n}{k}(\Seq{s_n},\Seq{\omega_n},\Set{\zeta\di{n}{k}}) = N\di{n}{k} / D\di{n}{k}\>,
           \end{array}
       \end{equation}
      where the generalized difference operator defined in Eq.\
      (\ref{eq29a}) involves quantities $\zeta\di{n}{k}\ne 0$ for
      $k\in\Nset_0$. Special cases of the \JD\ transformation result
      from corresponding choices of the $\zeta\di{n}{k}$. From Eq.\
      (\ref{eqJDRec}) one easily obtains the alternative representation
      \begin{equation}\label{eqJD_31}
        \JD\di{n}{k}(\Seq{s_n},\Seq{\omega_n},\Set{\zeta\di{n}{k}}) =
        \frac
          {\widetilde\nabla\di{n}{k-1}\widetilde\nabla\di{n}{k-2}\dots \widetilde\nabla\di{n}{0}
            [s_n/\omega_n]}
          {\widetilde\nabla\di{n}{k-1}\widetilde\nabla\di{n}{k-2}\dots \widetilde\nabla\di{n}{0}
            [1/\omega_n]}
        \>.
      \end{equation}
      Thus, the $\JD\di{}{k}$ is a Levin-type sequence transformation
      of order $2k$.

    \subsubsection{\H\ Transformation and Generalized \H\
    Transformation}
      The \H\ transformation was introduced by Homeier \cite{Homeier92}
      and used or studied in a series of articles
      \cite{Homeier93,Oleksy96,TC-NA-97-1,TC-NA-97-3,TC-NA-97-4,Homeier98aah}.
      Target of the \H\ transformation are Fourier series
      \begin{equation}\label{eqFourierreell}
         s = A_0/2+ \sum_{j=1}^{\infty} \bigl(A_j \cos(j\alpha) + B_j
         \sin(j\alpha)\bigr)
      \end{equation}
      with partial sums
      $ 
      s_n = A_0/2+ \sum_{j=1}^{n} \bigl(A_j \cos(j\alpha) + B_j
      \sin(j\alpha)\bigr)
      $ 
      where the Fourier coefficients $A_n$ and $B_n$ have asymptotic
      expansions of the form
      \begin{equation} \label{eqaccept}
        C_n \sim \rho^n n^{\epsilon} \sum_{j=0}^{\infty} c_j n^{-j}
      \end{equation}
      for $n\to\infty$ with $\rho\in\Kset$, $\epsilon\in\Kset$ and
      $c_0\ne 0$.

      The \H\ transformation was critized by Sidi \cite{Sidi95} as very
      unstable and useless near singularities of the Fourier series.
      However, Sidi failed to notice
      that --- as in the case of the $d\di{}{1}$ transformation with
      $\xi_n=\tau n$ --- one can apply  also the \H\ transformation
      (and also most other Levin-type sequence transformations) to the
      subsequence $\Seq{s_{\xi_n}}$ of $\Seq{s_n}$. The new sequence
      elements
      $s_{\xi_n}=s_{\tau n}$ can be regarded as the partial sums of a
      Fourier series with $\tau$--fold frequency. Using this
      $\tau$-fold frequency approach, one can obtain stable and
      accurate convergence acceleration even in the vicinity of
      singularities
      \cite{TC-NA-97-1,TC-NA-97-3,TC-NA-97-4,Homeier98aah}.

      The \H\ transformation may be defined as
      \begin{equation}\label{eqFourierH}
        \renewcommand{\arraystretch}{1.1}
        \begin{array}{l@{}l}
        \displaystyle{} N\di{n}{0} &\displaystyle{}{}= (n+\beta)^{-1}\,s_n/\omega_{n}\>, \qquad
         D\di{n}{0} {}= (n+\beta)^{-1}/\omega_{n}\>,\\
         \displaystyle{} N\di{{n}}{k} &\displaystyle{}{}=
         (n+\beta)N\di{{n}}{k-1}
         + (n+2k+\beta) N\di{{n+2}}{k-1}
                    \\ &{}
                    - 2\cos(\alpha)(n+k+\beta)
                      N\di{{n+1}}{k-1}\>,\\
         \displaystyle{} D\di{{n}}{k} &\displaystyle{}{}=
         (n+\beta)D\di{{n}}{k-1}
                    + (n+2k+\beta) D\di{{n+2}}{k-1}
                    \\ &{}
                    - 2\cos(\alpha)(n+k+\beta)
                      D\di{{n+1}}{k-1}\>,\\
         \multicolumn{2}{l}{\displaystyle \H\di{n}{k}(\alpha,\beta,\Seq{s_n},\Seq{\omega_n})=
                     {N\di{n}{k}}/{D\di{n}{k}}\>,}\\
        \end{array}
      \end{equation}
      where $\cos\alpha \ne \pm 1$ and $\beta\in\Rset_{+}$.

      It can also be represented in the explicit form \cite{Homeier92}
      \begin{equation} \label{eqHexpl}
 \H\di{n}{k}(\alpha,\beta,\Seq{s_n},\Seq{\omega_n})
= \frac{\P[P\di{}{2k}(\alpha)] [(n+\beta)^{k-1} s_n/\omega_n]}
                {\P[P\di{}{2k}(\alpha)] [(n+\beta)^{k-1} /\omega_n]}
      \end{equation}
      where the $p\di{m}{2k}(\alpha)$ and the polynomial
      $P\di{}{2k}(\alpha)\in\Pset\di{}{2k}$ are defined via
      \begin{equation} \label{eqHpmk}
        P\di{}{2k}(\alpha)(x)=(x^2-2x \cos\alpha +1)^{k} =
        \sum_{m=0}^{2k} p\di{m}{2k}(\alpha) x^m\>
      \end{equation}
      and $\P$ is the polynomial operator defined in Eq.\
      (\ref{eqPolyOp}).
      This shows that the $\H\di{}{k}$ transformation is a Levin-type
      transformation of order $2k$. It is not convex.

      A subnormalized form is
      \begin{equation} \label{eqHexplA}
        \H\di{n}{k}(\alpha,\beta,\Seq{s_n},\Seq{\omega_n})=
        \frac{\displaystyle\sum_{m=0}^{2k} p\di{m}{2k}(\alpha)
        \frac{(n+\beta+m)^{k-1}}{(n+\beta+2k)^{k-1}}
        \frac{s_{n+m}}{\omega_{n+m}}}{\displaystyle\sum_{m=0}^{2k} p\di{m}{2k}(\alpha)
        \frac{(n+\beta+m)^{k-1}}{(n+\beta+2k)^{k-1}}
        \frac{1}{\omega_{n+m}}
        }\>.
      \end{equation}
      This relation shows that the limiting transformation
       \begin{equation} \label{eqRingelH}
       \Ringel\H\di{}{k}=
           \frac{\P[P\di{}{2k}(\alpha)] [ s_n/\omega_n]}
                {\P[P\di{}{2k}(\alpha)] [ 1/\omega_n]}
       \end{equation}
       exists,
       and has characteristic polynomial
      $P\di{}{2k}(\alpha)$.

      A generalized \H\ transformation was defined by Homeier
      \cite{Homeier96Hab,TC-NA-97-1}. It is given in terms of the
      polynomial $P\di{}{k,M}(\mathbf{e})\in\Pset\di{}{kM}$ with
           \begin{equation}\label{eqpolFourierplus}
              P\di{}{k,M}(\mathbf{e}) (x) = \prod_{m=1}^M (x-e_m)^k
               = \sum_{\ell=0}^{k M} p\di{\ell}{k,M}(\mathbf{e}) x^\ell
           \end{equation}
      where $\mathbf{e}=(e_1,\dots,e_M)\in\Kset^{M}$ is a vector of
      constant parameters. Then, the generalized \H\ transformation is
      defined as
      \begin{equation}\label{eqFourierplus}
         \H\di{n}{k,M}(\beta,\Seq{s_n},\Seq{\omega_n},\mathbf{e})
         = \frac{\P[P\di{}{k,M}(\mathbf{e})] [(n+\beta)^{k-1} s_n/\omega_n]}
                {\P[P\di{}{k,M}(\mathbf{e})] [(n+\beta)^{k-1} /\omega_n]}
      \end{equation}
      This shows that the generalized $\H\di{}{k,M}$ is a Levin-type
      sequence transformation of order $kM$.
      The generalized \H\ transformation can be computed recursively
      using the scheme \cite{Homeier96Hab,TC-NA-97-1}
      \begin{equation}\label{eqgenHalgo}
        \begin{array}{r@{}l}
         N\di{n}{0} &\displaystyle{}= (n+\beta)^{-1}\,s_n/\omega_{n}\>, \qquad
         D\di{n}{0} {}= (n+\beta)^{-1}/\omega_{n}\>,\\
         N\di{{n}}{k} &\displaystyle{}=\sum_{j=0}^{M} q_j \,(n+\beta+j\, k)
         N\di{n+j}{k-1}
                      \>,\\
         D\di{{n}}{k} &\displaystyle{}=\sum_{j=0}^{M} q_j \,(n+\beta+j\, k)
         D\di{n+j}{k-1}
                      \>,\\
         \multicolumn{2}{l}{\displaystyle
           \H\di{n}{k,M}(\beta,\Seq{s_n},\Seq{\omega_n},\mathbf{e})
           =\frac{N\di{n}{k}}
                 {D\di{n}{k}}
           \>.}
         \end{array}
      \end{equation}
      Here, the $q_j$ are defined by
      \begin{equation}
        \prod_{m=1}^M (x-e_m) = \sum_{j=0}^{M} q_j x^j\>.
      \end{equation}
      The algorithm (\ref{eqFourierH}) is a special case of the
      algorithm (\ref{eqgenHalgo}). To see this, one observes that
      $M=2$, $e_1=\exp(\mbox{i}\alpha)$ und $e_2=\exp(-\mbox{i}\alpha)$
      imply $q_0=q_2=1$ and $q_1=-2\cos(\alpha)$.

      For $M=1$ and $e_1=1$, the Levin transformation is recovered.

    \subsubsection{\I\ Transformation}
      The \I\ transformation was in a slightly different form
      introduced by Homeier \cite{Homeier93}. It was
      derived via (asymptotically) hierarchically consistent iteration
      of the $\H\di{}{1}$ transformation, i.e., of
       \begin{equation} \label{eqH1}
         s_n'  =
          \frac
             {\displaystyle{s_{n+2}}/{\omega_{n+2}} - 2 \cos(\alpha)\,
                {s_{n+1}}/{\omega_{n+1}} + {s_n}/{\omega_n}}
             {\displaystyle{1}/{\omega_{n+2}} - 2
                \cos(\alpha)/{\omega_{n+1}}       + {1}/{\omega_n}}
       \end{equation}
      For the derivation and an analysis of the properties of the \I\
      transformation see \cite{Homeier96Hab,Homeier98aah}.
      The \I\ transformation may be defined via the recursive scheme
       \begin{equation} \label{eqIRec}
          \begin{array}{l@{}l}
          \displaystyle{} N\di{n}{0} &\displaystyle{}= s_n /
          \omega_n\>, \qquad D\di{n}{0} = 1/ \omega_n\>, \\
          \displaystyle{} N\di{n}{k+1} &
          \displaystyle{}
          =\bigtriangledown\di{n}{k}[\alpha] N\di{n}{k}  \>, \\
          \displaystyle{} D\di{n}{k+1} &\displaystyle{}
          = \bigtriangledown\di{n}{k}[\alpha] D\di{n}{k}\>, \\
          \multicolumn{2}{l}{\displaystyle \I\di{n}{k}(\alpha,\Seq{s_n},\Seq{\omega_n},
                       \Set{\Delta\di{n}{k}})
          =\frac{{N}\di{n}{k} }
          {{D}\di{n}{k}}\>.}
          \end{array}
       \end{equation}
      where the generalized difference operator $\bigtriangledown\di{n}{k}[\alpha]$
      defined in Eq.\
      (\ref{eq29b}) involves quantities $\Delta\di{n}{k}\ne 0$ for
      $k\in\Nset_0$. Special cases of the \I\ transformation result
      from corresponding choices of the $\Delta\di{n}{k}$. From Eq.\
      (\ref{eqIRec}) one easily obtains the alternative representation
      \begin{equation}\label{eqI_31}
        \I\di{n}{k}(\Seq{s_n},\Seq{\omega_n},\Set{\Delta\di{n}{k}}) =
        \frac
          {\bigtriangledown\di{n}{k-1}[\alpha]
           \bigtriangledown\di{n}{k-2}[\alpha]\dots
            \bigtriangledown\di{n}{0}[\alpha]
            [s_n/\omega_n]}
          {\bigtriangledown\di{n}{k-1}[\alpha]
           \bigtriangledown\di{n}{k-2}[\alpha]\dots
            \bigtriangledown\di{n}{0}[\alpha]
            [1/\omega_n]}
        \>.
      \end{equation}
      Thus, $\I\di{}{k}$ is a Levin-type sequence transformation
      of order $2k$.  It is not convex.

      Put $\Theta_n^{(0)}=1$ and for $k>0$ define
\begin{equation}
\Theta_n^{(k)}= \frac{\Delta_n^{(0)}\dots \Delta_n^{(k-1)}}
                   {\Delta_{n+1}^{(0)}\dots \Delta_{n+1}^{(k-1)}}\>.
\end{equation}
       If  for all $k\in \Nset$ the limits
\begin{equation}\label{eqFourierPhi}
\lim_{n\to\infty} \Theta_{n}^{(k)} = \Theta_k
\end{equation}
exist (we have always $\Theta_0=1$), then one can define a limiting
trans\-for\-ma\-ti\-on $\Ringel\I$ for large $n$. It is a special case
of the \I\ transformation according to \cite{Homeier98aah}
\begin{equation}
\Ringel{\I}\di{n}{k}(\alpha,\Seq{s_n},\Seq{\omega_n},\Seq{\Theta_k})
=
\I\di{n}{k}(\alpha,\Seq{s_n},\Seq{\omega_n},
                    \Seq{(\Theta_{k}/\Theta_{k+1})^n })\>.
\end{equation}
This is a transformation of  order $2k$.
    The characteristic polynomials of $\Ringel\I$ are known
    \cite{Homeier98aah} to be
    \begin{equation}
     Q\di{}{2k}(\alpha)\in \Pset^{2k}\>:\> Q\di{}{2k}(\alpha)(z) =
    \prod_{j=0}^{k-1} \left[(1-z\Theta_j\exp(i\alpha))(1-z\Theta_j\exp(-i\alpha)
)\right]\>.
      \end{equation}

    \subsubsection{\K\ Transformation}
      The \K\ transformation was introduced by Homeier
      \cite{Homeier94nca} in a slightly different form. It was obtained
      via iteration of the simple transformation
       \begin{equation} \label{eqK1}
         s_n'  =
        \frac{\displaystyle
        \zeta_{n}^{(0)} \frac{s_{n}}{\omega_{n}} +
        \zeta_{n}^{(1)} \frac{s_{n+1}}{\omega_{n+1}} +
        \zeta_{n}^{(2)} \frac{s_{n+2}}{\omega_{n+2}}
        }
       {\displaystyle
        \zeta_{n}^{(0)} \frac{1}{\omega_{n}} +
        \zeta_{n}^{(1)} \frac{1}{\omega_{n+1}} +
        \zeta_{n}^{(2)} \frac{1}{\omega_{n+2}}
        }
       \end{equation}
       that is exact for sequences of the form
       \begin{equation}\label{eqOrthomodseq}
         s_n = s + \omega_n\, (c\, P_n + d\, Q_n )\>,
      \end{equation}
      where $c$ and $d$ are arbitrary constants, while $P_n$ and $Q_n$ are two linearly independent
      solutions of
      the three-term recurrence
      \begin{equation}\label{eqOrthorec}
        \zeta_{n}^{(0)} v_{n} + \zeta_{n}^{(1)} v_{n+1} +
        \zeta_{n}^{(2)} v_{n+2} = 0\>.
      \end{equation}
      The \K\ transformation may be defined via the recursive scheme
       \begin{equation} \label{eqKRec}
          \begin{array}{l@{}l}
          \displaystyle{} N\di{n}{0} &\displaystyle{}= s_n /
          \omega_n\>, \qquad D\di{n}{0} = 1/ \omega_n\>, \\
          \displaystyle{} N\di{n}{k+1} &\displaystyle{}=
          \partial\di{n}{k}[\zeta] N\di{n}{k}  \>, \\
          \displaystyle{} D\di{n}{k+1} &\displaystyle{}=
          \partial\di{n}{k}[\zeta] D\di{n}{k}\>, \\
          \multicolumn{2}{l}{\displaystyle \K\di{n}{k}(\Seq{s_n},\Seq{\omega_n},
          \Set{\widetilde\Delta\di{n}{k}},\Set{\zeta\di{n}{j}})
          =\frac{{N}\di{n}{k} }
          {{D}\di{n}{k}}\>.}
          \end{array}
       \end{equation}
      where the generalized difference operator
      $\partial\di{n}{k}[\zeta]$ defined in Eq.\
      (\ref{eq29c}) involves recursion coefficients $\zeta\di{n+k}{j}$ with $j=0,1,2$
      and  quantities $\widetilde\Delta\di{n}{k}\ne 0$
      for
      $k\in\Nset_0$. Special cases of the \K\ transformation for given
      recursion, i.e., for given $\zeta\di{n}{j}$, result
      from corresponding choices of the $\widetilde\Delta\di{n}{k}$. From Eq.\
      (\ref{eqKRec}) one easily obtains the alternative representation
      \begin{equation}\label{eqK_31}
        \K\di{n}{k}(\Seq{s_n},\Seq{\omega_n},\Set{\widetilde\Delta\di{n}{k}},\Set{\zeta\di{n}{j}}) =
        \frac
          {\partial\di{n}{k-1}[\zeta]
           \partial\di{n}{k-2}[\zeta]\dots
            \partial\di{n}{0}[\zeta]
            [s_n/\omega_n]}
          {\partial\di{n}{k-1}[\zeta]
           \partial\di{n}{k-2}[\zeta]\dots
            \partial\di{n}{0}[\zeta]
            [1/\omega_n]}
        \>.
      \end{equation}
      Thus, $\K\di{}{k}$ is a Levin-type sequence transformation
      of order $2k$.  It is not convex.

      For applications of the \K\
      transformation see
      \cite{Homeier94nca,Homeier96Hab,Homeier98oca,TC-NA-97-4}.

{\quad}\section{Methods for the Construction of Levin-Type
Transformations}\label{secMetCon}

  In this section, we discuss approaches for the construction of
  Levin-type sequence transformations and point out the relation to
  their kernel.

  \subsection{Model Sequences and Annihilation Operators}\label{secmodann}

    As discussed in the introduction, the derivation of sequence
    transformations may be based on model sequences. These may be of
    the form (\ref{eqbekaverf_Emodel})  or of the form
    (\ref{eqpsijnmod}). Here, we consider model sequences of the latter
    type that involves remainder estimates $\omega_n$. As described in
    Section \ref{secEalgo}, determinantal representations for the
    corresponding sequence transformations can be derived using
    Cramer's rule, and one of the recursive schemes of the \E\
    algorithm may be used for the computation. However, for important
    special choices of the functions $\psi_j(n)$, simpler recursive
    schemes and more explicit representations in the form (\ref{eqT})
    can be obtained using the annihilation operator approach of Weniger
    \cite{Weniger89}. This approach was also studied by Brezinski and
    Matos \cite{BrezinskiMatos96} who showed that it leads to a unified
    derivation of many extrapolation algorithms and related devices and
    general results about their kernels. Further, we mention the work
    of Matos \cite{Matos97} who analysed the approach further and
    derived a number of convergence acceleration results for Levin-type
    sequence transformations.

    In this approach, an annihilation operator $\A=\A\di{n}{k}$ as defined in Eq.\
    (\ref{eqannihil}) is needed that annihilates the sequences
    $\Seq{\psi_j(n)}$, i.e., such that
       \begin{equation} \label{eqannipsi}
         \A\di{n}{k}(\Seq{\psi_j(n)}) = 0 \quad \mbox{for }
         j=0,\dots,k-1\>.
       \end{equation}
    Rewriting Eq.\ (\ref{eqpsijnmod}) in the form
    \begin{equation} \label{eqpsijnmodmod}
      \frac{\sigma_n - \sigma}{ \omega_n} = \sum_{j=0}^{k-1} c_j
      \psi_{j}(n)\>,
    \end{equation}
    and applying \A\ to both sides of this equation, one sees that
    \begin{equation} \label{eqpsijnmod1}
      \A\di{n}{k}\SEQ{\frac{\sigma_n - \sigma}{ \omega_n}} = 0
    \end{equation}
    This equation may be solved for $\sigma$ due to the linearity of
    \A. The result is
        \begin{equation} \label{eqpsijnmod2}
          \sigma=\frac{\A\di{n}{k}(\Seq{\sigma_n/ \omega_n})}
                      {\A\di{n}{k}(\Seq{1/\omega_n})}
        \end{equation}
     leading to a sequence transformation
        \begin{equation} \label{eqTanni}
          \T\di{n}{k}(\Seq{s_n},\Seq{\omega})
                =\frac{\A\di{n}{k}\left(\Seq{s_n/ \omega_n}\right)}
                      {\A\di{n}{k}\left(\Seq{1/\omega_n}\right)}
        \end{equation}
      Since \A\ is linear, this transformation can be rewritten in the
      form (\ref{eqT}), i.e., a Levin-type transformation has been
      obtained.

      We note that this process can be reversed, that is, for each
      Levin-type sequence transformation $\T[\Lambda\di{}{k}]$
      of order $k$ there is an annihilation operator, namely the
      polynomial operator $\P[\Pi\di{n}{k}]$ as defined in Eq.\
      (\ref{eqPolyOp}) where $\Pi\di{n}{k}$ are the characteristic
      polynomials as defined in Eq.\ (\ref{eqpink}). Using this
      operator, the defining equation (\ref{eqsnd}) can be rewritten as
      \begin{equation} \label{eqsndmod}
          \T\di{n}{k}(\Seq{s_n},\Seq{\omega_n})=
                 \frac{\P[\Pi\di{n}{k}](s_{n}/\omega_{n})}
                      {\P[\Pi\di{n}{k}](1/\omega_{n})}
      \end{equation}
      Let $\phi\di{n,m}{k}$ for $m=0,\dots,k-1$ be $k$ linearly
      independent solutions of the linear $(k+1)$--term recurrence
       \begin{equation} \label{eqPrec}
         \sum_{j=0}^{k} \lambda\di{n,j}{k} v_{n+j} = 0\>.
       \end{equation}
      Then $\P[\Pi\di{n}{k}] \phi\di{n,m}{k}=0$ for $m=0,\dots,k-1$,
      i.e., $\P[\Pi\di{n}{k}]$ is an annihilation operator for all
      solutions of Eq.\ (\ref{eqPrec}). Thus, all sequences that are
      annihilated by this operator are linear combinations of the $k$
      sequences $\Seq{\phi\di{n,m}{k}}$.

      If $\Seq{\sigma_n}$ is a sequence in the kernel of $\T\di{}{k}$ with
      (anti)limit $\sigma$, we must have
      \begin{equation} \label{eqsigmod}
          \sigma=
                 \frac{\P[\Pi\di{n}{k}](\sigma_{n}/\omega_{n})}
                      {\P[\Pi\di{n}{k}](1/\omega_{n})}
      \end{equation}
      or after some rearrangement using the linearity of \P
      \begin{equation} \label{eqsigmod1}
             \P[\Pi\di{n}{k}] \left(
                 \frac{\sigma_{n}-\sigma}
                      {\omega_{n}}
                             \right ) =0 \>.
      \end{equation}
      Hence, we must have
      \begin{equation} \label{eqsigmod2}
                 \frac{\sigma_{n}-\sigma}
                      {\omega_{n}} = \sum_{m=0}^{k-1} c_m \phi\di{n,m}{k} \>,
      \end{equation}
      or, equivalently
      \begin{equation} \label{eqsigmod3}
                 \sigma_{n}=\sigma
                      +\omega_{n} \sum_{m=0}^{k-1} c_m \phi\di{n,m}{k} \>,
      \end{equation}
      for some constants $c_m$. Thus, we have determined the kernel of
      $\T\di{}{k}$ that can also be considered as the set of model
      sequences for this transformation. Thus, we have proved the
      following theorem:

      \begin{theorem}\label{thTkernel}
        Let $\phi\di{n,m}{k}$ for $m=0,\dots,k-1$ be the $k$ linearly
        independent solutions of the linear $(k+1)$--term recurrence
        (\ref{eqPrec}). The kernel of
        $\T[\Lambda\di{}{k}](\Seq{s_n},\Seq{\omega_n})$ is given by all
        sequences $\Seq{\sigma_n}$ with (anti)limit $\sigma$ and
        elements $\sigma_n$ of the form (\ref{eqsigmod3}) for arbitrary
        constants $c_m$.
      \end{theorem}

      We note that the $\psi_j(n)$ for $j=0,\dots,k-1$ can essentially
      be identified with the $\phi\di{n,j}{k}$. Thus, we have
      determinantal representations for known $\psi_j(n)$ as noted
      above in the context of the E algorithm. See also
      \cite{Homeier95} for determinantal representations of the \J\
      transformations and the relation to its kernel.

      Examples of annihilation operators and the functions $\psi_j(n)$
      that are annihilated are given in Table \ref{tabanni}. Examples
      for the Levin-type sequence transformations that have been
      derived using the approach of model sequences are discussed in
      Section \ref{secEXanni}.

      Note that the annihilation operators used by Weniger
      \cite{Weniger89,Weniger92,Weniger94} were weighted difference
      operators $\W\di{n}{k}$ as defined in Eq.\
      (\ref{eqweighteddiff}). Homeier
      \cite{Homeier94ahc,Homeier95,Homeier96aan} discussed operator
      representations for the \J\ transformation that are equivalent to
      many of the annihilation operators and related sequence
      transformations as given by Brezinski and Matos
      \cite{BrezinskiMatos96}. The latter have been further discussed
      by Matos \cite{Matos97} who considered among others Levin-type sequence
      transformations with constant coefficients,
      $\lambda\di{n,j}{k}=const.$, and with polynomial coefficients
      $\lambda\di{n,j}{k}=\lambda_j(n+1)$, with $\lambda_j\in \Pset$, and
      $n\in\Nset_0$, in
      particular annihilation operators of the form
       \begin{equation} \label{eqMatos1}
         L(u_n) = (\Omega^l + \lambda_1 \Omega^{l-1} + \dots +
         \lambda_l) (u_n)
       \end{equation}
       with the special cases
       \begin{equation} \label{eqMatos11}
         L_1(u_n) = (\Omega -\alpha_1) (\Omega-\alpha_2) \cdots
         (\Omega -\alpha_l)
         (u_n)\>,\qquad (\alpha_i\ne\alpha_j\>\mbox{ for all } i\ne
         j\>,
       \end{equation}
       and
       \begin{equation} \label{eqMatos12}
         L_2(u_n) = (\Omega -\alpha)^{l}
         (u_n)
       \end{equation}
       where
       \begin{equation} \label{eqMatos2}
         \Omega^r(u_n) = (n+1)_r u_{n+r} \>,\qquad n\in\Nset_0
       \end{equation}
      and
       \begin{equation} \label{eqMatos3}
         \widetilde L(u_n) = (\pi -\alpha_1) (\pi-\alpha_2) \cdots
         (\pi -\alpha_l)
         (u_n)
       \end{equation}
       where
       \begin{equation} \label{eqMatos4}
         \pi(u_n) = (n+1)\forwarddiff u_{n}\>,\qquad \pi^{r}(u_n) =
         \pi(\pi^{r-1}(u_n))\>,\qquad n\in\Nset_0
       \end{equation}
       and the $\lambda$'s and $\alpha$'s are constants. Note that $n$
       is shifted in comparison to \cite{Matos97} where the convention
       $n\in\Nset$ was used. See also Table
       \ref{tabanni} for the corresponding annihilated functions
       $\psi_j(n)$.

       Matos \cite{Matos97} also considered difference operators of the
       form
       \begin{equation} \label{eqMatos5}
         L(u_n) = \forwarddiff^{k} + p_{k-1}(n)\forwarddiff^{k-1} + \cdots
          + p_1(n)\forwarddiff + p_0(n)
       \end{equation}
       where the functions $f_j$ given by $f_j(t)=p_j(1/t) t^{-k+j}$
       for $j=0,\dots,k-1$ are analytic in the neighborhood of 0. For
       such operators, there is no explicit formula for the functions
       that are annihilated. However, the asymptotic behavior of such
       functions is known \cite{BenderOrszag87,Matos97}. We will later
       return to such annihilation operators and state some convergence
       results.

\begin{table}[htp]
\caption{Examples of Annihilation Operators \ensuremath{{}^{\dagger}}}
\label{tabanni}
$$
          \begin{array}{l>{\displaystyle}l|>{\displaystyle}l}
            \hline\hline
            \mbox{Type} & \mbox{Operator} & \psi_j(n)\>,j=0,\dots,k-1 \\
            \hline
            \mbox{Differences}
                        & \forwarddiff^k     & (n+\beta)^j \\
                        &                   & (n+\beta)_j \\
                        &                   & (\alpha[n+\zeta])^j \\
                        &                   & (\alpha[n+\zeta])_j \\
                        &                   & p_j(n)\>,p_j\in\Pset\di{}{j} \\
            \mbox{Weighted Differences}
                        & \forwarddiff^k (n+\beta)^{k-1}    & 1/(n+\beta)^j \\
                        &\forwarddiff^k (n+\beta)_{k-1}    & 1/(n+\beta)_j \\
                        &\forwarddiff^k (\alpha[n+\zeta])^{k-1} & 1/(\alpha[n+\zeta])^j \\
                        &\forwarddiff^k (\alpha[n+\zeta])_{k-1} & 1/(\alpha[n+\zeta])_j \\
            \mbox{Divided Differences}
                        &\square\di{n}{k}[\Seq{t_n}] & t_n^{j} \\
                        &\square\di{n}{k}[\Seq{t_n}] & p_j(t_n)\>,p_j\in\Pset\di{}{j} \\
                        &\square\di{n}{k}[\Seq{x_n}] (x_n)_{k-1} & 1/(x_n)_j \\
            \mbox{Polynomial}
                        & \P[P\di{}{2k}(\alpha)]   & \exp(+i\alpha n)
                                                     p_j(n)\>,p_j\in\Pset\di{}{j} \\
                        &                          & \exp(-i\alpha n)
                                                     p_j(n)\>,p_j\in\Pset\di{}{j} \\
                        & \P[P\di{}{2k}(\alpha)] (n+\beta)^{k-1}  &
                        \exp(+ i\alpha n) / (n+\beta)^j \\
                        &                                         &
                        \exp(- i\alpha n) / (n+\beta)^j \\
                        & \P[P\di{}{2k}(\alpha)] (n+\beta)_{k-1}  &
                          \exp(+ i\alpha n) / (n+\beta)_j \\
                        &                                         &
                          \exp(- i\alpha n) / (n+\beta)_j \\
                        & \P[P\di{}{k}] & \psi_j(n) \mbox{\ is solution of }
                             \\
                        &               &
                        \sum_{m=0}^{k} p\di{n}{k} v_{n+j} = 0    \\
                        & \P[P\di{}{k}] (n+\beta)^{m} & (n+\beta)^m \psi_j(n) \mbox{\ is solution of }
                         \\
                        &                             &
                        \sum_{m=0}^{k} p\di{n}{k} v_{n+j} =0\\
                        & L_1   \mbox{\ (see (\ref{eqMatos11}))}                      &
                        \frac{\alpha_j^{n+1}}{n!} \\
                        & L_2   \mbox{\ (see (\ref{eqMatos12}))}                      &
                        \frac{n^j\alpha^{n+1}}{n!} \\
                        & \widetilde L   \mbox{\ (see (\ref{eqMatos3}))}                      &
                        \frac{\Gamma(n+\alpha_j+1)}{n!} \\
            \hline\hline
            \multicolumn{3}{l}{\mbox{\ensuremath{{}^{\dagger}} See also Section \ref{secEXanni}.}}
           \end{array}
$$
\end{table}

    \subsubsection{Derivation of the \F\ Transformation}
      As an example for the application of the annihilation operator
      approach, we derive the \F\ transformation. Consider the model
      sequence
       \begin{equation} \label{eqFmod}
         \sigma_n = \sigma + \omega_n \sum_{j=0}^{k-1} c_j
         \frac{1}{(x_n)_j}
       \end{equation}
       that may be rewritten as
       \begin{equation} \label{eqFmodA}
         \frac{\sigma_n - \sigma}{\omega_n} = \sum_{j=0}^{k-1} c_j
         \frac{1}{(x_n)_j}\>.
       \end{equation}
       We note that Eq. (\ref{eqFmod}) corresponds to modeling $\mu_n=R_n/\omega_n$
       as a truncated factorial series in $x_n$ (instead as a truncated
       power series as in the case of the W algorithm). The $x_n$ are elements
       of  $\Seq{x_n}$ an auxiliary sequence  $\Seq{x_n}$ such that
       $\lim_{n\to\infty} 1/x_{n}=0$ and also $x_{n+\ell}>x_{n}$ for
       $\ell \in \Nset$ and $x_0>1$, i.e., a diverging sequence  of
       monotonously increasing positive numbers. To find an
       annihilation operator
       for the $\psi_j(n)=1/(x_n)_j$, we make use of the fact that the divided
       difference operator $\square\di{n}{k}=\square\di{n}{k}[\Seq{x_n}]$ annihilates
       polynomials in $x_n$ of degree less than $k$. Also, we observe
       that the definition of the Pochhammer symbols entails that
       \begin{equation} \label{eqxnpoch}
         (x_n)_{k-1}/ (x_n)_j = (x_n+j)_{k-1-j}\>,\qquad
       \end{equation}
       is a polynomial of degree less than $k$ in $x_n$ for $0\le j\le
       k-1$. Thus, the sought annihilation operator is $\A=\square\di{n}{k}
       (x_n)_{k-1}$ because
       \begin{equation} \label{eqxnpochann}
         \square\di{n}{k}(x_n)_{k-1}  \frac{1}{(x_n)_j} = 0\>,\quad
         0\le j\le k-1\>.
       \end{equation}
       Hence, for the model sequence (\ref{eqFmod}), one can calculate
       $\sigma$ via
        \begin{equation} \label{eqFmod2}
          \sigma=\frac{\square\di{n}{k}({(x_n)_{k-1}\sigma_n/ \omega_n})}
                      {\square\di{n}{k}({(x_n)_{k-1}/\omega_n})}\>
        \end{equation}
        and the \F\ transformation (\ref{eqFtrans}) results by
        replacing $\sigma_n$ by $s_n$ in
        the right hand side of Eq.\ (\ref{eqFmod2}).

    \subsubsection{Important Special Cases}\label{secEXanni}
      Here, we collect model sequences and annihilation operators for
      some important Levin-type sequence transformations that were
      derived using the model sequence approach. For further examples
      see also \cite{BrezinskiMatos96}. The model sequences
      are the kernels by construction. In Section \ref{secEXiter},
      kernels and annihilation operators are stated for important
      Levin-type transformation that were derived using iterative
      methods.

      \paragraph{Levin transformation}
        The model sequence for $\L\di{}{k}$ is
       \begin{equation} \label{eqmodLevin}
         \sigma_n = \sigma + \omega_n\sum_{j=0}^{k-1} c_j /
         (n+\beta)^{j}\>.
       \end{equation}
       The annihilation operator is
       \begin{equation} \label{eqannLevin}
         \A\di{n}{k} =\forwarddiff^k (n+\beta)^{k-1}\>.
       \end{equation}

      \paragraph{Weniger transformations}
        The model sequence for $\S\di{}{k}$ is
       \begin{equation} \label{eqmodS}
         \sigma_n = \sigma + \omega_n\sum_{j=0}^{k-1} c_j /
         (n+\beta)_{j}\>.
       \end{equation}
       The annihilation operator is
       \begin{equation} \label{eqannS}
         \A\di{n}{k} =\forwarddiff^k (n+\beta)_{k-1}\>.
       \end{equation}
        The model sequence for $\M\di{}{k}$ is
       \begin{equation} \label{eqmodM}
         \sigma_n = \sigma + \omega_n\sum_{j=0}^{k-1} c_j /
         (-n-\xi)_{j}\>.
       \end{equation}
       The annihilation operator is
       \begin{equation} \label{eqannM}
         \A\di{n}{k} =\forwarddiff^k (-n-\xi)_{k-1}\>.
       \end{equation}
        The model sequence for $\C\di{}{k}$ is
       \begin{equation} \label{eqmodC}
         \sigma_n = \sigma + \omega_n\sum_{j=0}^{k-1} c_j /
         (\alpha[n+\zeta])_{j}\>.
       \end{equation}
       The annihilation operator is
       \begin{equation} \label{eqannC}
         \A\di{n}{k} =\forwarddiff^k (\alpha[n+\zeta])_{k-1}\>.
       \end{equation}

      \paragraph{W algorithm}
        The model sequence for $W\di{}{k}$ is
       \begin{equation} \label{eqmodW}
         \sigma_n = \sigma + \omega_n\sum_{j=0}^{k-1} c_j t_n^j\>.
       \end{equation}
       The annihilation operator is
       \begin{equation} \label{eqannW}
         \A\di{n}{k} =\square\di{n}{k}[\Seq{t_n}]\>.
       \end{equation}

      \paragraph{\H\ Transformation}
        The model sequence for $\H\di{}{k}$ is
       \begin{equation} \label{eqmodH}
         \sigma_n = \sigma + \omega_n\left(
         \exp(i\alpha n)
            \sum_{j=0}^{k-1} c^{+}_j /(n+\beta)^{j}
         + \exp(-i\alpha n)\sum_{j=0}^{k-1} c^{-}_j
                  /(n+\beta)^{j}
                  \right)\>.
       \end{equation}
       The annihilation operator is
       \begin{equation} \label{eqannH}
         \A\di{n}{k} =\P[P\di{}{2k}(\alpha)](n+\beta)^{k-1}\>.
       \end{equation}

       \paragraph{Generalized \H\ Transformation}
          The model sequence for $\H\di{}{k,M}$ is
           \begin{equation}\label{eqmodFourierplus}
             \sigma_n=\sigma+\omega_n \sum_{m=1}^{M} e_m^n \sum_{j=0}^{k-1} c_{m,j}
               (n+\beta)^{-j}\>.
           \end{equation}
          The annihilation operator is
           \begin{equation}\label{eqannFourierplus}
             \A\di{n}{k}=\P[P\di{}{k,M}(\mathbf{e})] (n+\beta)^{k-1}\>.
           \end{equation}

  \subsection{Hierarchically Consistent Iteration}\label{seciter}
    As alternative to the derivation of sequence transformations using
    model sequences and possibly annihilation operators, one may take
    some simple sequence transformation $T$ and iterate it $k$ times to obtain a
    transformation $T\di{}{k}=T\circ \dots \circ T$. For the iterated
    transformation, by
    construction one has a simple algorithm by
    construction, but the theoretical
    analysis is complicated since usually no kernel is known. See for
    instance the iterated Aitken process where the $\Delta^2$ method
    plays the r{\^o}le of the simple transformation. However, as is
    discussed at length in Refs.\ \cite{Weniger91,Homeier94ahc}, there
    are usually several possibilities for the iteration. Both problems
    -- unknown kernel and arbitrariness of iteration -- are overcome
    using the concept of hierarchical consistency
    \cite{Homeier94ahc,Homeier96Hab,Homeier98aah} that was shown  to
    give rise to powerful
    algorithms like the \J\ and the \I\ transformations
    \cite{Homeier96aan,Homeier96Hab,Homeier98aah}. The basic idea of
    the concept is to provide a hierarchy of model sequences such that
    the simple transformation provides a mapping between neighboring
    levels of the hierarchy. To ensure the latter, normally one has to
    fix some parameters in the simple transformation to make the
    iteration consistent with the hierarchy.

    A formal description of the concept is given in the following taken
    mainly from the literature.
    \cite{Homeier98aah}. As an example, the concept is later used to derive the
    \JD\ transformation in Section \ref{secJD}.

Let  $\Seq{\sigma_n(\vec c,\vec p)}_{n=0}^{\infty}$ be a simple
``basic'' model sequence that depends on a vector $\vec c
\in \Kset^{a}$ of constants, and further parameters
$\vec p$. Assume that its (anti)limit $\sigma(\vec p)$
exists and is independent of $\vec c$. Assume that the
basic transformation $T=T(\vec p)$ allows to compute the (anti)limit exactly
according to
\begin{equation}\label{eqhier_1}
T(\vec p\,): \Seq{\sigma_n(\vec c,\vec p\,)}
\longrightarrow
\Seq{\sigma(\vec p\,)}\>.
\end{equation}
Let
the hierarchy of model sequences be given by
\begin{equation}
\Set{\Seq{\sigma_n^{(\ell)}(\vec c\,{}^{(\ell)},\vec
p^{(\ell)})\vert \vec c^{(\ell)}\in \Kset^{a^{(\ell)}}
}}_{\ell=0}^{L}
\end{equation}
with $a^{(\ell)}>a^{(\ell')}$ for $\ell>\ell'$. Here, $\ell$ numbers
the levels of the hierarchy. Each of the model sequences
$\Seq{\sigma_n^{(\ell)}(\vec c\,{}^{(\ell)},\vec p^{(\ell)})}$ depends on an
$a^{(\ell)}$--dimensional complex vector $\vec c^{(\ell)}$
and further parameters $\vec p^{(\ell)}$.
Assume that the model sequences of lower levels are also contained in
those of higher levels: For all $\ell<L$ and
all $\ell'>\ell$ and $\ell'\le L$, every sequence $\Seq{\sigma_n^{(\ell)}(\vec
c^{(\ell)},\vec p\,{}^{(\ell)})}$ is assumed to be representable as a
model sequence
$\Seq{\sigma_n^{(\ell')}(\vec c\,{}^{(\ell')},\vec p\,{}^{(\ell')})}$
where
$\vec c\,{}^{(\ell')}$ is obtained from $\vec c\,{}^{(\ell)}$ by the
natural injection
$\Kset^{a^{(\ell)}}\to\Kset^{a^{(\ell')}}$.
Assume that for all $\ell$ with $0<\ell\le L$
\begin{equation}\label{eqhier_2}
T(\vec p\,{}^{(\ell)}) : \Seq{
   \sigma_n^{(\ell)}
       (\vec c\,{}^{(\ell)},\vec p\,{}^{(\ell)})}
\longrightarrow
   \Seq{
     \sigma_n^{(\ell-1)}
       (\vec c\,{}^{(\ell-1)},\vec p\,{}^{(\ell-1)})}
\end{equation}
is a mapping between
neighboring levels of the hierarchy.
Composition yields an iterative transformation
\begin{equation}
T^{(L)} = T(\vec p\,{}^{(0)}) \circ T(\vec p\,{}^{(1)}) \circ \dots
\circ T(\vec p\,{}^{(L)})
\end{equation}
This transformation  is called ``hierarchically consistent'' or
``consistent with the hierarchy''. It maps model sequences
$\sigma_n^{(\ell)}(\vec c\,{}^{(\ell)},\vec p^{(\ell)})$ to constant
sequences if Eq.\ (\ref{eqhier_1}) holds with
\begin{equation}
\Seq{\sigma_n^{(0)}(\vec
c\,{}^{(0)},\vec p\,{}^{(0)})}= \Seq{\sigma_n(\vec c,\vec p\,)}\>.
\end{equation}
If instead of Eq.\ (\ref{eqhier_2}) we have
\begin{equation}\label{eqhier_4}
T(\vec p\,{}^{(\ell)})\left(\Seq{\sigma_n^{(\ell)}(\vec c\,{}^{(\ell)},\vec
p^{(\ell)})}\right) \sim
\Seq{\sigma_n^{(\ell-1)}(\vec c\,{}^{(\ell-1)},\vec p\,{}^{(\ell-1)})}\>
\end{equation}
for $n\to\infty$ for all $\ell>0$ then the iterative transformation
$T^{(L)}$ is called ``asymptotically consistent with the hierarchy'' or
``asymptotically hierarchy-consistent''.

    \subsubsection{Derivation of the \JD\ Transformation}\label{secJD}
      The simple transformation is the $\D\di{}{2}$ transformation
       \begin{equation} \label{eqJDsimp}
         s_n'=T(\Seq{\omega_n})(\Seq{s_n})  = \frac{\forwarddiff^2 (s_n/\omega_n)}
                      {\forwarddiff^2 (1/\omega_n)}
       \end{equation}
       depending on the ``parameters'' $\Seq{\omega_n}$,
       with basic model sequences
       \begin{equation} \label{eqJDsimpmod}
         \frac{\sigma_n}{\omega_n} = \sigma\frac{1}{ \omega_n} +(a n +b)\>.
       \end{equation}

       The more complicated model sequences of the next level are taken to be
       \begin{equation} \label{eqJDs1mod}
         \frac{\sigma_n}{\omega_n} = \sigma\frac{1}{\omega_n} + (a n +b + (a_1 n+b_1)r_n)\>.
       \end{equation}
       Application of $\forwarddiff^2$ eliminates the terms involving $a$
       and $b$. The result is
       \begin{equation} \label{eqJDs1modA}
         \frac{\displaystyle \forwarddiff^2\frac{\sigma_n}{\omega_n}}
              {\forwarddiff^2 r_n}
                =
         \sigma \frac{\displaystyle\forwarddiff^2\frac{1}{\omega_n}}
                     {\forwarddiff^2 r_n}
                     +
               + \left( a_1 n + b_1 + 2 a_1
                    \frac{{\forwarddiff r_n}}{{\forwarddiff^2 r_n}}
               \right)\>
       \end{equation}
       for $\forwarddiff^2 r_n\ne 0$. Assuming that for large $n$
       \begin{equation} \label{eqrnass}
             \frac{\forwarddiff r_n}
                  {\forwarddiff^2 r_n}
             = A n + B + o(1)
       \end{equation}
       holds, the result is asymptotically of the same \emph{form} as the model sequence
       in Eq.\ (\ref{eqJDsimpmod}), namely
       \begin{equation} \label{eqJDsimpmodd}
         \frac{\sigma_n'}{\omega_n'} = \sigma\frac{1}{ \omega_n'} +(a'
         n +b' + o(1))\>.
       \end{equation}
       with renormalized ``parameters''
       \begin{equation} \label{eqomd}
         1/\omega_n' = \frac
                           {\forwarddiff^2(1/\omega_n)}
                           {\forwarddiff^2 r_n}\>
       \end{equation}
       and obvious identifications for $a'$ and $b'$.

       We now assume that this mapping between two neighboring levels
       of the hierarchy can be extended to any two neighboring levels,
       provided that one introduces $\ell$-dependent quantities,
       especially $r_n\to r\di{n}{\ell}$ with
       $\zeta\di{n}{\ell}=\forwarddiff^2 r\di{n}{\ell}\ne 0$,
       $s_n/\omega_n \to N\di{n}{\ell}$, $1/\omega_n\to D\di{n}{\ell}$
       and
       $s_n'/\omega_n' \to N\di{n}{\ell+1}$, $1/\omega_n'\to
       D\di{n}{\ell+1}$.

       Iterating in this way leads to the algorithm (\ref{eqJDRec}).

       The condition (\ref{eqrnass}) or more generally
       \begin{equation} \label{eqrnlass}
             \frac{\forwarddiff r\di{n}{\ell}}
                  {\forwarddiff^2 r\di{n}{\ell}}
             = A_\ell n + B_\ell + o(1)
       \end{equation}
       for given $\ell$ and for large $n$ is satisfied in many cases.
       For instance, it is satisfied if there are constants
       $\beta_\ell\ne 0$, $\gamma_\ell$ and $\delta_\ell\ne 0$ such that
       \begin{equation} \label{eqrnlg}
             \forwarddiff r\di{n}{\ell} \sim
             \beta_\ell \left\{
                        \begin{array}{ll}
                          \displaystyle
                          \left(\frac{\delta_\ell+1}{\delta_\ell}\right)^n
                          & \mbox{for } \gamma_\ell= 0 \\
                          \displaystyle
                          \left.
                          \left(\frac{\delta_\ell+1}{\gamma_\ell}\right)_n
                          \right/
                          \left(\frac{\delta_\ell}{\gamma_\ell}\right)_n
                          & \mbox{otherwise}
                        \end{array}
                      \right.
                      \>.
       \end{equation}
       This is for instance the case for $r\di{n}{\ell}=n^{\zeta_\ell}$
       with $\zeta_\ell(\zeta_\ell-1)\ne 0$.

       The kernel of $\JD\di{}{k}$ may be found inductively in the
       following way:
       \begin{equation} \label{eqJDindu}
          \begin{array}{l>{\displaystyle}r@{}>{\displaystyle}l}
           &   N\di{n}{k} -\sigma D\di{n}{k} &{}=0  \\
         \Longrightarrow  &  \forwarddiff^2 (N\di{n}{k-1} -\sigma D\di{n}{k-1}) &{}=0  \\
         \Longrightarrow  &   N\di{n}{k-1} -\sigma D\di{n}{k-1} &{}=a_{k-1} n + b_{k-1} \\
         \Longrightarrow  &   \forwarddiff^2 (N\di{n}{k-2} -\sigma D\di{n}{k-2})
            &{}=(a_{k-1} n + b_{k-1})\zeta\di{n}{k-2} \\
         \Longrightarrow  &   N\di{n}{k-2} -\sigma D\di{n}{k-2})
            &{}=a_{k-2} n + b_{k-2} + \sum_{j=0}^{n-2}
            \sum_{n'=0}^{j}(a_{k-1} n' + b_{k-1})\zeta\di{n'}{k-2}
           \end{array}
       \end{equation}
       yielding the result
       \begin{equation} \label{eqJDkern}
          \begin{array}{>{\displaystyle}r@{}>{\displaystyle}l}
            N\di{n}{0}-\sigma D\di{n}{0} ={}& \left. a_0n+b_0 + \sum_{j=0}^{n-2}
            \sum_{n_1=0}^{j}\zeta\di{n_1}{0}\right(a_{1} n_1 + b_{1} + \dots \\
            & \left.\dots (a_{k-2} n +b_{k-2} + \sum_{j_{k-2}=0}^{n_{k-2}-2}
            \sum_{n_{k-1}=0}^{j_{k-2}}\zeta\di{n_{k-1}}{k-2}(a_{k-1}
            n_{k-1} + b_{k-1}) \right)\>.
           \end{array}
       \end{equation}
       Here, the definitions $N\di{n}{0}=\sigma_n/\omega_n$ and
       $D\di{n}{0}=1/\omega_n$ may be used to obtain the model sequence
       $\Seq{\sigma_n}$
       for $\JD\di{}{k}$, that may be identified as kernel of that
       transformation, and also may be regarded as model sequence of
       the $k$-th level according to $\Seq{\sigma_n^{(k)}(\vec
       c^{(k)},\vec p\,{}^{(k)})}$ with $\vec
       c\di{}{k}=(a_0,b_0,\dots,a_{k-1},b_{k-1})$ and $\vec p\di{}{k}$
       corresponds to $\omega\di{n}{k}=1/D\di{n}{k}$ and the
       $\Set{\zeta\di{n}{\kappa}| 0\le\kappa\le k-2}$.

       We note this as a theorem:

       \begin{theorem}\label{thJDkernel}
         The kernel of $\JD\di{}{k}$ is given by the set of sequences
         $\Seq{\sigma_n}$ such that Eq.\
         (\ref{eqJDkern}) holds with $N\di{n}{0}=\sigma_n/\omega_n$ and
         $D\di{n}{0}=1/\omega_n$.
       \end{theorem}

    \subsubsection{Important Special Cases}\label{secEXiter}
      Here, we give the hierarchies of model sequences for sequence
      transformations derived via hierarchically consistent iteration.
      \paragraph{\J\ transformation}
      The most prominent example is the \J\ transformation
      (actually a large class of transformations). The corresponding
      hierarchy of model sequences provided by the kernels that are
      explicitly known according to the following theorem:

      \begin{theorem}\label{thJkernel} \cite{Homeier94ahc}
      {The kernel of the $\J\di{}{k}$ transformation is given by the
      sequences $\Seq{\sigma_n}$  with elements of the form
     \begin{equation} \label{eqgenkernel}
      \sigma_n = \sigma + \omega_n \sum_{j=0}^{k-1} c_j \psi_j(n)\>
     \end{equation}
     with
       \begin{equation} \label{eqJtranskernel}
          \begin{array}{>{\displaystyle}r@{}>{\displaystyle}l}
             \psi_0(n) &{}={} 1
          \\
             \psi_1(n) &{}={} \sum_{n_1=0}^{n-1} \delta\di{n_1}{0}
          \\
             \psi_2(n) &{}={} \sum_{n_1=0}^{n-1} \delta\di{n_1}{0}
                 \sum_{n_2=0}^{n_1-1} \delta\di{n_2}{1}
          \\
            {}           & {}\vdots  {}
          \\
             \psi_{k-1}(n) &{}={} \sum_{n>n_1>n_2>\cdots>n_{k-1}}
             \delta\di{n_1}{0}\delta\di{n_2}{1}\cdots
             \delta\di{n_{k-1}}{k-2}
           \end{array}
       \end{equation}
       with arbitrary constants $c_0,\cdots,c_{k-1}$. }
       \end{theorem}

       \paragraph{\I\ transformation} Since the \I\ transformation is a
       special case of the \J\ transformation (cp.\ Table
       \ref{tabJtransSpecial} and \cite{Homeier98aah}, its kernels
       (corresponding to the
      hierarchy of model sequences) are
      explicitly known according to the following theorem:

\begin{theorem}\label{thkernelI}\cite[Theorem 8]{Homeier98aah}
The kernel of the $\I\di{}{k}$ transformation
is given by the
sequences $\Seq{\sigma_n}$  with elements of the form
\begin{equation}\label{eqkernelII}
\begin{array}{r@{}l}
\displaystyle{}\sigma_n &\displaystyle{}{}= \sigma + \exp(-i\alpha n)\omega_n
           \Biggl[
             d_0 + d_1 \exp(2i\alpha n)
 \\
\displaystyle{}&\displaystyle{}{}                          +
                 \sum_{n_1=0}^{n-1} \sum_{n_2=0}^{n_1-1}
                 \exp(2i\alpha (n_1-n_2))
                 \Biggl(d_2 + d_3 \exp(2i\alpha n_2)
                 \Biggr)
                 \Delta_{n_2}^{(0)} + \cdots
 \\
\displaystyle{}&\displaystyle{}{}             +
\sum_{n>n_1>n_2>\cdots>n_{2k-2}}
\exp(2i\alpha [n_1-n_2+\dots +n_{2k-3}-n_{2k-2}])
 \\
\displaystyle{}&\displaystyle{}{} \times
        \Biggl(d_{2k-2} + d_{2k-1}
                 \exp(2i\alpha n_{2k-2})
                 \Biggr) \prod_{j=0}^{k-2}
\Delta_{n_{2j+2}}^{(j)}
           \Biggr]
\end{array}
\end{equation}
with constants $d_0,\cdots,d_{2k-1}$.
Thus,
we have
$s=
\I\di{n}{k'}(\alpha,\Seq{s_n},\Seq{\omega_n},\Set{\Delta\di{n}{k}})$
for $k'\ge k$ for sequences of this form.
\end{theorem}

\subsection{A Two-Step Approach}
In favorable cases, one may use a two-step approach for the construction of 
sequence transformations:
\begin{description}
\item[Step 1] Use asymptotic analysis of the remainder $R_n=s_n-s$ of the given problem to
              find the adequate model sequence (or hierarchy of model sequences) for large
              	$n$.
\item[Step 2] Use the methods described in Sections \ref{secmodann} or \ref{seciter} to
construct the sequence transformation adapted to the problem.
\end{description}
This is, of course, a mathematically promising approach.
A good example for the two-step approach is the derivation of the $d\di{}{m}$ transformations
by Levin and Sidi \cite{LevinSidi81} (compare also Section \ref{appdm}).

But there are two difficulties with this approach. 

The first difficulty is a practical one. In many cases, the problems to be treated in
applications are simply too complicated to allow to perform Step 1 of the two-step approach.

The second difficulty is a more mathematical one. The optimal system of functions $f_j(n)$ used in the asymptotic expansion, 
\begin{equation}
s_n -s \sim \sum_{j=0}^{\infty} c_j f_j(n)
\end{equation}
with $f_{j+1}(n)=o(f_j(n))$,
i.e.,
the optimal \emph{asymptotic scale} \cite[p. 2]{Wimp81},  is not clear \emph{a priori}. 
For instance, as the work of
Weniger has shown, sequence transformations like the Levin transformation that  are based 
on expansions in powers of $1/n$, i.e.,
the asymptotic scale $\phi_j(n)=1/(n+\beta)^j$, are not always superior to, and even often worse 
than those based
upon
factorial series, like Weniger's \S\
transformation that is based on the asymptotic scale $\psi_j(n)=1/(n+\beta)_j$. 
To find an optimal asymptotic scale in combination with nonlinear sequence transformations 
seems to be an open mathematical problem.

Certainly, the proper choice of remainder estimates \cite{HomeierWeniger95} is also crucial in the context of
Levin-type sequence transformations. See also Section \ref{secApplic}.

{\quad}\section{Properties of Levin-Type
Transformations}\label{secProper}

  \subsection{Basic Properties}
    Directly from the definition in Eqs.\ (\ref{eqLevinType}) and
    (\ref{eqsnd}), we obtain the following theorem. The proof is left
    to the interested reader.
    \begin{theorem}\label{thsimple}
       Any Levin-type sequence transformation \T\ is quasilinear, i.e.,
       we have
\begin{equation}
\T\di{n}{k}(\Seq{A s_n +B},\Seq{\omega_n}) =
A \T\di{n}{k}(\Seq{s_n},\Seq{\omega_n}) + B
\end{equation}
for arbitrary constants $A$ and $B$.
       It is multiplicatively invariant in $\omega_n$, i.e., we have
\begin{equation}
\T\di{n}{k}(\Seq{s_n},\Seq{C\omega_n}) =
\T\di{n}{k}(\Seq{s_n},\Seq{\omega_n})
\end{equation}
for arbitrary constants $C\ne 0$.
    \end{theorem}

      For a coefficient set $\Lambda$ define the sets
      $Y\di{n}{k}[\Lambda]$ by
       \begin{equation} \label{eqYLambdank}
         Y\di{n}{k}[\Lambda] = \Set{ (x_0,\dots,x_{k}) \in
         \Fset^{k+1}\>\vert\>
         \sum_{j=0}^{k} \lambda\di{n,j}{k} / x_j \ne 0}\>.
       \end{equation}
      Since $\T\di{n}{k}(\Seq{s_n},\Seq{\omega_n})$ for given
       coefficient set $\Lambda$ depends only on the
       $2k+2$ numbers
       $s_n,\dots,s_{n+k}$ and $\omega_n,\dots,\omega_{n+k}$,
       it  may be regarded as a mapping
       \begin{equation} \label{eqUnk}
         U\di{n}{k} \>:\> \Cset^{k+1} \times Y\di{n}{k}[\Lambda]
         \Longrightarrow \Cset \>,\quad
         (x,y) \mapsto U\di{n}{k}(x
                        \>\vert\> y)
       \end{equation}
       such that
       \begin{equation}
       \T\di{n}{k}=U\di{n}{k}(s_n,\dots,s_{n+k}
                     \>\vert\>\omega_n,\dots,\omega_{n+k})
       \end{equation}

The following theorem is a generalization of theorems for the \J\
transformation \cite[Theorem 5]{Homeier94ahc} and the \I\
transformation \cite[Theorem 5]{Homeier98aah}.

\begin{theorem}\label{thgenprop}
\begin{enumerate}
\item[(I-0)]
The $\T\di{}{k}$ transformation  can
be regarded as continous mapping $U\di{n}{k}$ on
{\rm $\Cset^{k+1}\times Y\di{n}{k}[\Lambda]$} where
{\rm $Y\di{n}{k}[\Lambda]$} is defined in Eq.\ (\ref{eqYLambdank}).
\item[(I-1)] According to Theorem \ref{thsimple},
$U\di{n}{k}$ is a homogeneous function of first
degree in the first
$(k+1)$ variables and a homogeneous function of degree zero in the
last
$(k+1)$ variables. Hence, for all vectors
{\rm $\vec x\in \Cset^{k+1}$} and {\rm $\vec y\in Y\di{n}{k}[\Lambda]$}
and for all complex constants
$s$ and $t\ne 0$ the equations
\begin{equation}
\begin{array}{r@{}l}
\displaystyle{}&\displaystyle{}U\di{n}{k}(s \vec x\,\vert\,\vec y)
= s
U\di{n}{k}( \vec x\,\vert\,\vec y)\>,\\
\displaystyle{}&\displaystyle{}U\di{n}{k}( \vec x\,\vert\,
t \vec y) = U\di{n}{k}( \vec x\,\vert\,\vec y)\>
\end{array}
\end{equation}
hold.
\item[(I-2)] $U\di{n}{k}$ is linear in the first  $(k+1)$
variables. Thus, for all vectors {\rm $\vec x\in \Cset^{k+1}$,
$\vec x'\in \Cset^{k+1}$,} und {\rm $\vec y\in Y\di{n}{k}[\Lambda]$}
\begin{equation}
U\di{n}{k}( \vec x+\vec x'\,\vert\,\vec y)=U\di{n}{k}( \vec x\,\vert\,\vec
y)+U\di{n}{k}( \vec x'\,\vert\,\vec y)\>
\end{equation}
holds.
\item[(I-3)] For all constant vectors {\rm $\vec
c=(c,c,\dots,c)\in\Cset^{k+1}$} and all vectors {\rm $\vec
y\in Y\di{n}{k}[\Lambda]$} we have
\begin{equation}\label{eqFouriercc}
U\di{n}{k} (\vec c \,\vert\, \vec y) = c \>.
\end{equation}
\end{enumerate}
\end{theorem}
\begin{proof} These are immediate consequences of the definitions.
\end{proof}

   \subsection{The Limiting Transformation}

We note that if a limiting transformation $\Ringel\T[\Ringel\Lambda]$
exists, it is also of Levin-type, and thus, the above theorems apply to
the limiting transformation as well.

Also, we have the following result for the kernel of the limiting
transformation:

    \begin{theorem}\label{thlimkernel}
      Suppose that for a Levin-type sequence transformation
      $\T\di{}{k}$ of order $k$ there exists a limiting transformation
      $\Ringel\T\di{}{k}$ with characteristic
      polynomial $\Ringel\Pi\in\Pset^{k}$ given by
      \begin{equation} \label{eqRingelPi}
         \Ringel\Pi\di{}{k}(z) = \sum_{j=0}^k \Ringel\lambda\di{j}{k}
         z^j = \prod_{\ell=1}^{M} (z-\zeta_{\ell})^{m_{\ell}}
      \end{equation}
      where the zeroes $\zeta_\ell\ne 0$ have multiplicities
      $m_{\ell}$.
      Then the kernel of the limiting transformation consists of all
      sequences $\Seq{s_n}$ with elements of the form
       \begin{equation} \label{eqlimkernel}
         \sigma_n=\sigma +\omega_n  \sum_{\ell=1}^{M} \zeta_{\ell}^{n}
         P_{\ell}(n)
       \end{equation}
       where $P_{\ell}\in \Pset^{m_{\ell}-1}$ are arbitrary
       polynomials and $\Seq{\omega_n}\in \Ringel \Yset\di{}{k}$.
    \end{theorem}
    \begin{proof}
      This follows directly from the observation that for such
      sequences $(\sigma_n-\sigma)/\omega_n$ is nothing but a finite linear
      combination of the solutions $\varphi\di{n,\ell,j_{\ell}}{k} =
      n^{j_{\ell}} \zeta_\ell^n$ with $\ell=1,\dots,M$ and
      $j_{\ell}=0,\dots,m_{\ell}-1$ of the recursion relation
       \begin{equation} \label{eqlimPrec}
         \sum_{j=0}^{k} \Ringel\lambda\di{j}{k} v_{n+j} = 0\>,
       \end{equation}
       and thus, it is annihilated by $\P[\Ringel\Pi\di{}{k}]$.
    \end{proof}

  \subsection{Application to Power Series}

   Here, we generalize some results of Weniger \cite{Weniger94} that
   regard the application of Levin-type sequence transformations to
   power series.

   We use the definitions in Eq.\ (\ref{eq33}).

   Like Pad\'e approximants, Levin-type sequence transformations yield
   rational approximants when applied to the partial sums $f_n(z)$ of a power
   series $f(z)$ with terms $a_j=c_j z^j$. These approximations offer a practical way for the
   analytical continuation of power series to regions outside of their
   circle of convergence. Furthermore, the poles of the rational
   approximations model the singularities of $f(z)$. They may also be
   used to approximate further terms beyond the last one used in
   constructing the rational approximant.

   When applying a Levin-type sequence transformation \T\ to a power
   series, remainder estimates $\omega_n=m_n z^{\gamma+n}$ will be
   used. We note that $t$ variants correspond to $m_n=c_n$, $\gamma=0$,
   $u$ variants correspond to $m_n=c_n(n+\beta)$, $\gamma=0$,
   $\widetilde t$ variants to $m_n=c_{n+1}$, $\gamma=1$. Thus, for
   these variants, $m_n$ is independent of $z$ (Case A). For $v$
   variants, we have $m_n=c_{n+1}c_n/(c_n-c_{n+1}z)$, and $\gamma=1$.
   In this case, $1/m_n\in \Pset\di{}{1}$ is a linear function of $z$
   (Case B).

   Application of \T\ yields after some simplification
   \begin{equation} \label{eqTpow1}
     \T\di{n}{k}(\Seq{f_n(z)},\Seq{m_n z^{\gamma+n}}) =
     \frac
          {\displaystyle \sum_{\ell=0}^{n+k} z^\ell
          \sum_{j=\max(0,k-\ell)}^{k}
          \frac{\lambda\di{n,j}{k}}{m_{n+j}} c_{\ell-(k-j)} }
          {\displaystyle \sum_{j=0}^{k}
          \frac{\lambda\di{n,j}{k}}{m_{n+j}} z^{k-j} }
=     \frac
          {P\di{n}{k}[T] (z)}
          {Q\di{n}{k}[T] (z)}\>,
   \end{equation}
   where in Case A, we have $P\di{n}{k}[T]\in \Pset^{n+k}$,
   $Q\di{n}{k}[T]\in\Pset^{k}$, and in Case B, we have
   $P\di{n}{k}[T]\in \Pset^{n+k+1}$,
   $Q\di{n}{k}[T]\in\Pset^{k+1}$. One needs the $k+1+\gamma$ partial sums
   $f_n(z),\dots,f_{n+k+\gamma}(z)$ to compute these rational approximants. This
   should be compared to the fact that for the computation of the
   Pad\'e approximant $[n+k+\gamma/k+\gamma]$ one needs the
   $2k+2\gamma+1$ partial
   sums $f_n(z),\dots,f_{n+2k+2\gamma}(z)$.

   We show that Taylor expansion of these rational approximants
   reproduces all terms of power series that have been used to calculate
   the rational approximation.

\begin{theorem}\label{thpow}
We have
  \begin{equation} \label{eqTpowth}
    \T\di{n}{k}(\Seq{f_n(z)},\Seq{m_n z^{\gamma+n}})-f(z) =
    O(z^{n+k+1+\tau})
  \end{equation}
  where $\tau=0$ for $t$ and $u$ variants corresponding to $m_n=c_n$,
  $\gamma=0$, or
   $m_n=c_n(n+\beta)$, $\gamma=0$, respectively, while $\tau=1$ holds
   for the $v$ variant corresponding to
   $m_n=c_{n+1}c_n/(c_n-c_{n+1}z)$, $\gamma=1$, and for the
   $\widetilde t$ variants corresponding to $m_n=c_{n+1}$,
   $\gamma=1$, one obtains $\tau=1$ if \T\ is convex.
\end{theorem}
\begin{proof}
   Using the identity
   \begin{equation} \label{eqTpow2}
     \T\di{n}{k}(\Seq{f_n(z)},\Seq{m_n z^{\gamma+n}}) =
     f(z)+\T\di{n}{k}(\Seq{f_n(z)-f(z)},\Seq{m_n z^{\gamma+n}}) =
   \end{equation}
   that follows from Theorem \ref{thsimple}, we obtain after some easy
   algebra
   \begin{equation} \label{eqTpow3}
     \T\di{n}{k}(\Seq{f_n(z)},\Seq{m_n z^{\gamma+n}})-f(z) =
     z^{n+k+1}
     \frac
          {\displaystyle \sum_{\ell=0}^{\infty} z^{\ell}
          \sum_{j=0}^{k}
          \frac{\lambda\di{n,j}{k}}{m_{n+j}} c_{\ell+n+j+1} }
          {\displaystyle \sum_{j=0}^{k}
          \frac{\lambda\di{n,j}{k}}{m_{n+j}} z^{k-j} }\>.
   \end{equation}
   This shows that the right hand side is at least $O(z^{n+k+1})$ since
   the denominator is $O(1)$ due to $\lambda\di{n,k}{k}\ne 0$. For
   the $\widetilde t$ variant, the term corresponding to $\ell=0$ in
   the numerator is $\sum_{j=0}^{k}\lambda\di{n,j}{k}=\Pi\di{n}{k}(1)$
   that vanishes for convex \T. For the $v$ variant, that term is
   $\sum_{j=0}^{k}\lambda\di{n,j}{k} (c_{n+j}-c_{n+j+1}z)/c_{n+j}$ that
   simplifies to $(-z) \sum_{j=0}^{k}\lambda\di{n,j}{k}
   c_{n+j+1}/c_{n+j}$ for convex \T. This finishes the proof.
\end{proof}

{\quad}\section{Convergence Acceleration Results \\ for
Levin-Type Transformations}\label{secConRes}

  \subsection{General Results}
    We note that Germain-Bonne \cite{GermainBonne73} developed  a
    theory of the regularity and convergence acceleration properties of
    sequence transformations that was later extended by Weniger
    \cite[Section 12]{Weniger89}, \cite[Section 6]{Weniger94} to
    sequence transformations that depend explicitly on $n$ and on an
    auxiliary sequence of remainder estimates. The essential results of
    this theory apply to convergence acceleration of linearly
    convergent sequences. Of course, this theory can be applied to
    Levin-type sequence transformations. However, for the latter
    transformations, many results can be obtained more easily and also,
    one may obtain results of a general nature that are also applicable
    to other convergence types like logarithmic convergence. Thus, we
    are not going to use the Germain-Bonne-Weniger theory in the
    present article.

    Here, we present some general convergence acceleration results for
    Levin-type sequence transformations that have a limiting
    transformation. The results, however, do not completely determine
    which transformation provides the best extrapolation results
    for a given problem sequence
    since the results are asymptotic in nature, but in practice, one is
    interested in obtaining good extrapolation results from as few
    members of the problem sequence as possible. Thus, it may well be
    that transformations with the same asymptotic behavior of the
    results perform rather differently in practice.

    Nevertheless, the results presented below provide a first
    indication which results one may expect for large classes of
    Levin-type sequence transformations.

    First, we present some results that show that the limiting
    transformation essentially determines for which sequences
    Levin-type sequence transformations are accelerative. The speed of
    convergence will be analyzed later.

    \begin{theorem}\label{thconver0}
      Assume that the following asymptotic relations hold for large $n$:
      \begin{equation}
       \lambda\di{n,j}{k} \sim \Ringel\lambda\di{j}{k}\>,
       \qquad
       \Ringel\lambda\di{k}{k}\ne 0\>,
      \end{equation}
      \begin{equation}
        \frac{s_n-s}{\omega_n} \sim \sum_{\nu=1}^{A} c_{\nu}
        \zeta_{\nu}^{n}\>, \quad c_{\nu}\zeta_{\nu} \ne 0\>,\quad
        \Ringel\Pi\di{}{k}(\zeta_{\nu})=0\>,\quad
      \end{equation}
      \begin{equation}
        \frac{\omega_{n+1}}{\omega_n} \sim \rho\ne 0\>, \quad
        \Ringel\Pi\di{}{k}(1/\rho)\ne 0\>.
      \end{equation}
      Then, $\Seq{\T\di{n}{k}}$ accelerates $\Seq{s_n}$ to $s$, i.e.,
      we have
      \begin{equation}
        \lim_{n\to\infty} \frac{\T\di{n}{k}-s}{s_n-s} = 0\>.
      \end{equation}
    \end{theorem}
    \begin{proof}
      Rewriting
      \begin{equation}
       \frac{\T\di{n}{k}-s}{s_n-s} =
       \frac{\omega_n}{s_n-s}
       \frac{\displaystyle
            \sum_{j=0}^{k} \lambda\di{n,j}{k}
            \frac{s_{n+j}-s}{\omega_{n+j}}
            }
            {\displaystyle
              \sum_{j=0}^{k} \lambda\di{n,j}{k}
              \frac{\omega_{n}}{\omega_{n+j}}
            }
      \end{equation}
      one may perform the limit for $n\to\infty$ upon using the
      assumptions according to
      \begin{equation}
       \frac{\T\di{n}{k}-s}{s_n-s} \to
       \frac{\displaystyle
             \sum_{j=0}^{k} \Ringel\lambda\di{j}{k}
             \sum_{\nu} c_{\nu} \zeta_{\nu}^{n+j}
             }
            {\displaystyle
             \sum_{\nu} c_{\nu} \zeta_{\nu}^{n}
             \sum_{j=0}^{k} \Ringel\lambda\di{j}{k} \rho^{-j}
             }
        = \frac{\displaystyle
               \sum_{\nu} c_{\nu}\zeta_{\nu}^{n}
               \Ringel\Pi\di{}{k}(\zeta_{\nu})
               }
               {\displaystyle
               \Ringel\Pi\di{}{k}(1/\rho)
               \sum_{\nu} c_{\nu}\zeta_{\nu}^{n}
               } =0
      \end{equation}
      since $\omega_n/\omega_{n+j}\to\rho^{-j}$.
    \end{proof}

    Thus, the zeroes $\zeta_{\nu}$ of the characteristic polynomial of
    the limiting transformation are of particular importance.

    It should be noted that the above assumptions correspond to a more
    complicated convergence type than linear or logarithmic convergence
    if $\abs{\zeta_1}=\abs{\zeta_2}\ge\abs{\zeta_3}\ge \dots$. This is
    the case, for
    instance, for the $\H\di{}{k}$ transformation where the limiting
    transformation has the characteristic polynomial
    $P\di{}{2k}(\alpha)$ with $k$-fold zeroes
    at $\exp(\mbox{i}\alpha)$ and $\exp(-\mbox{i}\alpha)$. Another example is the
    $\I\di{}{k}$ transformation where the limiting transformation has
    characteristic polynomials $Q\di{}{2k}(\alpha)$ with zeroes at
    $\exp(\pm \mbox{i} \alpha)/\Theta_j$, $j=0,\dots,k-1$.

    Specializing to $A=1$ in Theorem \ref{thconver0}, we obtain the
    following corollary:

    \begin{corollar}\label{thconver0a}
      Assume that the following asymptotic relations hold for large $n$:
      \begin{equation}
       \lambda\di{n,j}{k} \sim \Ringel\lambda\di{j}{k}\>,
       \qquad
       \Ringel\lambda\di{k}{k}\ne 0\>,
      \end{equation}
      \begin{equation}
        \frac{s_n-s}{\omega_n} \sim c
        q^{n}\>, \quad cq \ne 0\>,\quad
        \Ringel\Pi\di{}{k}(q)=0\>,\quad
      \end{equation}
      \begin{equation}
        \frac{\omega_{n+1}}{\omega_n} \sim \rho\ne 0\>, \quad
        \Ringel\Pi\di{}{k}(1/\rho)\ne 0\>.
      \end{equation}
      Then, $\Seq{\T\di{n}{k}}$ accelerates $\Seq{s_n}$ to $s$, i.e.,
      we have
      \begin{equation}
        \lim_{n\to\infty} \frac{\T\di{n}{k}-s}{s_n-s} = 0\>.
      \end{equation}
    \end{corollar}
    Note that the assumptions of Corollary \ref{thconver0a} imply
    \begin{equation}
      \frac{s_{n+1}-s}{s_n-s} =
      \frac{s_{n+1}-s}{\omega_{n+1}}\frac{\omega_n}{s_n-s}
      \frac{\omega_{n+1}}{\omega_n} \sim \rho\frac{cq^{n+1}}{cq^n} =
      \rho q
    \end{equation}
    and thus, Corollary \ref{thconver0a} corresponds to linear
    convergence for $0<\abs{\rho q}<1$ and to logarithmic convergence
    for $\rho q=1$.

    Many important sequence transformations have convex limiting
    transformations, i.e., the characteristic polynomials satisfy
    $\Ringel\Pi\di{}{k}(1)=0$. In this case, they accelerate linear
    convergence. More exactly, we have the following corollary:

    \begin{corollar}\label{thconver0b}
      Assume that the following asymptotic relations hold for large $n$:
      \begin{equation}
       \lambda\di{n,j}{k} \sim \Ringel\lambda\di{j}{k}\>,
       \qquad
       \Ringel\lambda\di{k}{k}\ne 0\>,
      \end{equation}
      \begin{equation}
        \frac{s_n-s}{\omega_n} \sim c
        \>, \quad c \ne 0\>,\quad
        \Ringel\Pi\di{}{k}(1)=0\>,\quad
      \end{equation}
      \begin{equation}
        \frac{\omega_{n+1}}{\omega_n} \sim \rho\ne 0\>, \quad
        \Ringel\Pi\di{}{k}(1/\rho)\ne 0\>.
      \end{equation}
      Then, $\Seq{\T\di{n}{k}}$ accelerates $\Seq{s_n}$ to $s$, i.e.,
      we have
      \begin{equation}
        \lim_{n\to\infty} \frac{\T\di{n}{k}-s}{s_n-s} = 0\>.
      \end{equation}
      Hence, any Levin-type sequence transformation with a convex
      limiting transformation accelerates linearly convergent sequences
      with
    \begin{equation}
      \lim_{n\to\infty}
      \frac{s_{n+1}-s}{s_n-s} =
      \rho \>,\qquad 0<\abs{\rho}<1
    \end{equation}
      such that $\Ringel\Pi\di{}{k}(1/\rho)\ne 0$ for suitably chosen
      remainder estimates $\omega_n$ satisfying $(s_n-s)/\omega_n\to
      c\ne 0$.
    \end{corollar}
    \begin{proof} Specializing Corollary \ref{thconver0a} to
    $q=1$, it suffices to
    prove the last assertion. Here, the proof follows from the
    observation that
    $
      {(s_{n+1}-s)}/{(s_n-s)} \sim
      \rho
    $
    and $(s_n-s)/\omega_n\sim
      c$ imply $\omega_{n+1}/\omega_n\sim\rho$ for large $n$ in view of
      the assumptions.
    \end{proof}

    Note that Corollary \ref{thconver0b} applies for instance to
    suitable variants of the
    Levin transformation, the \pJ\ transformation and, more generally, of
    the \J\ transformation. In particular, it applies to $t$, $\tilde
    t$, $u$ and $v$ variants, since in the case of linear convergence,
    one has $\forwarddiff s_{n}/\forwarddiff
    s_{n-1}\sim \rho$ which entails $(s_n-s)/\omega_{n}\sim c$ for
    all these variants by simple algebra.

    Now, some results for the speed of convergence are given.
    Matos \cite{Matos97} presented convergence theorems for
    sequence transformations based on annihilation difference operators
    with characteristic polynomials with constants coefficients that
    are close in spirit to the theorems given below. However, it should be
    noted that the theorems presented here apply to large classes of
    Levin-type transformations that have a limiting transformation (the
    latter, of course, has a characteristic polynomial with constants
    coefficients).

    \begin{theorem}\label{thconver1}{\quad}\\
    \begin{description}
    \item[(C-1)]
      Suppose that for a Levin-type sequence transformation
      $\T\di{}{k}$ of order $k$ there
      is a limiting transformation $\Ringel\T\di{}{k}$ with characteristic
      polynomial $\Ringel\Pi\in\Pset^{k}$ given by Eq.\ (\ref{eqRingelPi})
      where the multiplicities $m_{\ell}$ of the zeroes $\zeta_\ell\ne 0$
      satisfy $m_1\le m_2 \le \dots \le m_M$. Let
      \begin{equation}\label{eqc1}
        \lambda\di{n,j}{k} \frac{n^{m_1-1}}{(n+j)^{m_1-1}}\sim
        \Ringel\lambda\di{j}{k}
        \left(
          \sum_{t=0}^{\infty}  \frac{e\di{t}{k}}{(n+j)^t}
        \right)
        \>,\qquad e\di{0}{k} =1
      \end{equation}
      for $n\to\infty$.
    \item[(C-2)]
      Assume that  $\Seq{s_n}\in\Sset^{\Kset}$ and
      $\Seq{\omega_n}\in\Oset^{\Kset}$. Assume further that for
      $n\to\infty$ the asymptotic expansion
       \begin{equation} \label{eqasy1}
         \frac{s_n-s}{\omega_n} \sim \sum_{\ell=1}^{M} \zeta_{\ell}^{n}
         \sum_{r=0}^{\infty} c_{\ell,r} n^{-r}
       \end{equation}
      holds, and put
       \begin{equation} \label{eqrell}
         r_{\ell}= \min \Set{r\in \Nset_0\>\vert\>
         f_{\ell,r+m_1}\ne 0}
       \end{equation}
       where
       \begin{equation}\label{eqc2}
         f_{\ell,v} = \sum_{r=0}^{v} e\di{v-r}{k} \; c_{\ell,r}
       \end{equation}
       and
       \begin{equation} \label{eqPiell}
         B_{\ell}= (-1)^{m_\ell}
                   \frac{d^{m_\ell}\Ringel\Pi\di{}{k}}
                        {dx^{m_\ell}}({\zeta_\ell})
       \end{equation}
      for $\ell=1,\dots,M$.
    \item[(C-3)]
      Assume that the following limit exists and satisfies
       \begin{equation} \label{eqasy2}
         0\ne \lim_{n\to\infty}\frac{\omega_{n+1}}{\omega_n} = \rho
         \not\in\Set{\zeta_{\ell}^{-1}\>\vert\>\ell=1,\dots,M}\>.
       \end{equation}
    \end{description}
      Then we have
       \begin{equation} \label{eqasy3}
         \frac{\T\di{n}{k}(\Seq{s_n},\Seq{\omega_n})-s}{\omega_n} \sim
         \frac{\displaystyle
         \sum_{\ell=1}^M f_{\ell,r_\ell+m_1} \zeta_{\ell}^{n+m_\ell}
         \binom{r_\ell+m_\ell}{r_\ell}
         \frac{B_\ell}
              {n^{ r_{\ell} + m_{\ell} - m_1}
              }
         }{\Ringel\Pi\di{}{k}(1/\rho)}\> \frac{1}{n^{2m_1}}\>.
       \end{equation}
       Thus, $\Seq{\T\di{n}{k}(\Seq{s_n},\Seq{\omega_n})}$ accelerates
       $\Seq{s_n}$ to $s$ at least with order $2m_1$, i.e.,
       \begin{equation}\label{eqc3}
         \frac{\T\di{n}{k}-s}{s_n-s} = O(n^{-2m_1-\tau})\;,\qquad \tau\ge
         0\>,
       \end{equation}
       if
       $c_{\ell,0}\ne 0$ for all $\ell$.
    \end{theorem}
    \begin{proof}
      We rewrite $\T\di{n}{k}(\Seq{s_n},\Seq{\omega_n})=\T\di{n}{k}$ as
      defined in Eq.\ (\ref{eqT}) in the form
      \begin{equation}
                \T\di{n}{k} -s =
                \omega_n\frac{\displaystyle\sum_{j=0}^{k} \lambda\di{n,j}{k}
                       \frac{s_{n+j}-s}{\omega_{n+j}}}
                      {\displaystyle\sum_{j=0}^{k} \lambda\di{n,j}{k}
                       \frac{\omega_n}{\omega_{n+j}}}
                \sim \omega_n\frac{\displaystyle\sum_{j=0}^{k}
                   \Ringel\lambda\di{j}{k} \sum_{t=0}^{\infty}  \frac{e\di{t}{k}}{(n+j)^t}
                   \frac{(n+j)^{m_1-1}}{n^{m_1-1}}
                       \frac{s_{n+j}-s}{\omega_{n+j}}}
                      {\displaystyle\sum_{j=0}^{k} \Ringel\lambda\di{j}{k}
                      \frac{1}{\rho^{j}}}
      \end{equation}
      for
      large $n$ where we used Eq.\ (\ref{eqc1}) in the numerator, and in the
      denominator the relation $\omega_{n}/\omega_{n+j}\to \rho^{-j}$ that follows by repeated application of Eq.\
      (\ref{eqasy2}). Insertion of (\ref{eqasy1}) now yields
      \begin{equation}
                \T\di{n}{k} -s
                \sim
                \frac{\omega_n}
                     {n^{m_1-1}\Ringel\Pi\di{}{k}(1/\rho)}
                \sum_{\ell=1}^{M} \sum_{r=0}^{\infty} f_{\ell,r+m_1}
                \displaystyle\sum_{j=0}^{k}
                   \Ringel\lambda\di{j}{k}
                   \frac{\zeta_{\ell}^{n+j}}{(n+j)^{r+1}}
      \end{equation}
      where Eq.\ (\ref{eqc2}) was used. Also the fact was used
      that $\P[\Ringel\Pi\di{}{k}]$ annihilates any
      linear combination of the solutions $\varphi\di{n,\ell,j_{\ell}}{k} =
      n^{j_{\ell}} \zeta_\ell^n$ with $\ell=1,\dots,M$ and
      $j_{\ell}=0,\dots,m_1-1$ of the recursion relation
       (\ref{eqlimPrec}) since each $\zeta_\ell$ is a zero with
       multiplicity exceeding $m_1-1$.  Invoking Lemma
       \ref{thlemma} given in  Appendix \ref{applemma} one
       obtains
      \begin{equation}
                \T\di{n}{k} -s
                \sim
                \frac{\omega_n}
                     {n^{m_1-1}\Ringel\Pi\di{}{k}(1/\rho)}
                \sum_{\ell=1}^{M} \sum_{r=0}^{\infty} f_{\ell,r+m_1}
\zeta_\ell^{n+m_\ell} \binom{r+m_\ell}{r} \frac{(-1)^{m_\ell}}{n^{r+m_\ell+1}}
  \frac{d^{m_\ell}\Ringel\Pi\di{}{k}}{dx^{m_\ell}}(\zeta_\ell)
      \end{equation}
       The proof of Eq.\ (\ref{eqasy3}) is completed taking
       leading terms in the sums over $r$. Since $s_n-s \sim
       \omega_nZ^n  \sum_{\ell\in I} (\zeta_{\ell}/Z)^{n}
         c_{\ell,0} $ where $Z=\max\Set{\vert \zeta_\ell\vert
         \>\vert \ell=1,\dots,M}$, and $I=\Set{\ell=1,\dots,M
         \vert\>Z=\vert \zeta_\ell\vert}$, Eq.\ (\ref{eqc3}) is
         obtained where $\tau=\min\Set{r_\ell+m_\ell-m_1\>\vert\>\ell
         \in I}$.
    \end{proof}

If $\omega_{n+1}/\omega_n\sim\rho$, where
$\Ringel\Pi\di{}{k}(1/\rho)=0$, i.e., if (C-3) of Theorem \ref{thconver1}
does not hold, then the denominators vanish
asymptotically. In this case, one has to investigate whether the
numerators or the denominators vanish faster.

    \begin{theorem}\label{thconver2}
    Assume that (C-1) and (C-2) of Theorem \ref{thconver1} hold.
    \begin{description}
    \item[(C-3')]
      Assume that for $n\to\infty$ the asymptotic relation
       \begin{equation} \label{eqasy3a}
         \frac{\omega_{n+1}}{\omega_n} \sim \rho
         \exp(\epsilon_n)\>, \quad \rho \ne 0
       \end{equation}
       holds where
       \begin{equation}
         \frac{1}{\lambda!}\frac{d^{\lambda}\Ringel\Pi\di{}{k}}{d\,x^{\lambda}} (1/\rho)
         = \left\{
            \begin{array}{ll}
              0 & \mbox{for } \lambda=0,\dots,\mu-1 \\
              C\ne 0 & \mbox{for } \lambda=\mu
            \end{array}
           \right.
       \end{equation}
       and
       \begin{equation}
         \epsilon_n \to 0\>,\qquad \frac{\epsilon_{n+1}}{\epsilon_n}
         \to 1
       \end{equation}
       for large $n$. Define $\delta_n$ via
       $\exp({-\epsilon_n})=1+\delta_n\rho$.
    \end{description}
      Then we have for large $n$
       \begin{equation} \label{eqasy3b}
         \frac{\T\di{n}{k}(\Seq{s_n},\Seq{\omega_n})-s}{\omega_n} \sim
         \frac{\displaystyle
         \sum_{\ell=1}^M f_{\ell,r_\ell+m_1} \zeta_{\ell}^{n+m_\ell}
         \binom{r_\ell+m_\ell}{r_\ell}
         \frac{B_\ell}
              {n^{ r_{\ell} + m_{\ell} - m_1}
              }
         }{C (\delta_n)^{\mu}}\> \frac{1}{n^{2m_1}}\>.
       \end{equation}
    \end{theorem}
    \begin{proof}
       The proof proceeds as the proof of Theorem \ref{thconver1} but
       in the denominator we use
       \begin{equation}\label{eqdenlim}
         \sum_{j=0}^{k}\lambda\di{n,j}{k} \frac{\omega_n}{\omega_{n+j}}
         \sim
         C (\delta_n)^\mu
       \end{equation}
       that follows from Lemma \ref{thlemma1} given in  Appendix
       \ref{applemma}.
    \end{proof}

Thus, the effect of the sequence transformation in this case
essentially depends on the question whether $(\delta_n)^{-\mu} n^{-2m_1}$
goes to 0 for large $n$ or not. In many important
    cases like the Levin transformation and the \pJ\ transformations,
    we have $M=1$ and $m_1=k$. We note that Theorem \ref{thconver2} becomes
especially important in the case of logarithmic convergence since for
instance for $M=1$ one observes that
    $
      {(s_{n+1}-s)}/{(s_n-s)} \sim 1
    $
    and
    $(s_n-s)/\omega_n\sim
      \zeta_1^{n} c_{1,0}\ne 0$
    imply $\omega_{n+1}/\omega_n\sim 1/\zeta_1$ for large $n$ such that
    the denominators vanish asymptotically. In this
    case, we have $\mu=m_1$ whence $(\delta_n)^{-\mu}
    n^{-2m_1}=O(n^{-m_1})$ if $\delta_n=O(1/n)$. This reduction of the
    speed of convergence of the acceleration process from $O(n^{-2k})$
    to $O(n^{-k})$ in the case of logarithmic convergence is a generic
    behavior that is reflected in a number of theorems regarding
    convergence acceleration properties of Levin-type sequence
    transformations. Examples are Sidi's theorem for the Levin
    transformation given below (Theorem \ref{thSidi}), and for the \pJ\
    transformation the Corollaries \ref{thCo1} and \ref{thCo2} given
    below. Compare also \cite[Theorems 13.5, 13.9,
13.11, 13.12, 14.2]{Weniger89}.

    The following theorem was given by Matos \cite{Matos97} where the
    proof may be found. To
    formulate it, we define that a sequence $\Seq{u_n}$ has
    \emph{property M} if it satisfies
    \begin{equation}\label{eqpropM}
              \frac{u_{n+1}}{u_n} \sim
              1+\frac{\alpha}{n} + r_n\> \mbox{with } r_n=o(1/n)\>,
              \forwarddiff^{\ell} r_n = o(\forwarddiff^\ell (1/n))
              \>\mbox{for } n\to\infty\>.
    \end{equation}

    \begin{theorem}\label{thmatos} \cite[Theorem 13]{Matos97}
      Let $\Seq{s_n}$ be a sequence such that
       \begin{equation} \label{eqMatos6}
         s_n-s = \omega_n (a_1 g\di{1}{1}(n)+\cdots + a_k
         g\di{1}{k}(n)+\rho_n)
       \end{equation}
       with $g\di{1}{j+1}(n)=o(g\di{1}{j}(n))$,
       $\rho_n=o(g\di{1}{k}(n))$ for $n\to\infty$. Let us consider an
       operator $L$ of the form (\ref{eqMatos5}) for which we know a
       basis of solutions $\Seq{u\di{n}{j}}$, $j=1,\dots,k$, and each
       one can be written as
       \begin{equation} \label{eqMatos7}
         u\di{n}{j} \sim \sum_{m=1}^{\infty} \alpha\di{m}{j}
         g\di{m}{j}(n)\>, \qquad g\di{m+1}{j}(n) = o(g\di{m}{j}(n))
       \end{equation}
       as $n\to\infty$ for all $m\in\Nset$ and $j=1,\dots,k$. Suppose
       that
       \begin{equation} \label{eqMatos8}
          \begin{array}{>{\displaystyle}r>{\displaystyle}l}
           a) & g\di{2}{j+1}(n) = o(g\di{2}{j}(n))\mbox{\ for }
           n\to\infty\>,j=1,\dots,k-1\\
           b) & g\di{2}{1}(n)=o(g\di{1}{k}(n))\>, \mbox{\ and } \rho_n
           \sim K g\di{2}{1}(n)\>\quad \mbox{for } n\to\infty \\
           c) & \Seq{g\di{m}{j}(n)} \>\mbox{has property M for } \\
              & \>m\in\Nset\>, j=1,\dots,k.
           \end{array}
       \end{equation}
       Then
       \begin{enumerate}
        \item If $\Seq{\omega_n}$ satisfies $\lim_{n\to\infty}
         \omega_{n}/\omega_{n+1}=\lambda\ne 1$, the sequence
         transformation $\T\di{n}{k+1}$ corresponding to the operator $L$
         accelerates the convergence of $\Seq{s_n}$. Moreover, the
         acceleration can be measured by
         \begin{equation} \label{eqMatos9}
           \frac{\T\di{n}{k+1}-s}{s_n-s} \sim C n^{-k}
           \frac{g\di{2}{1}(n)}{g\di{1}{1}(n)}\>, \qquad n\to\infty\>.
         \end{equation}
        \item If  $\Seq{1/\omega_n}$ has property M, then the speed of
        convergence of $\T\di{n}{k+1}$ can be measured by
         \begin{equation} \label{eqMatos10}
           \frac{\T\di{n}{k+1}-s}{s_n-s} \sim C
           \frac{g\di{2}{1}(n)}{g\di{1}{1}(n)}\>, \qquad n\to\infty\>.
         \end{equation}
       \end{enumerate}
    \end{theorem}

  \subsection{Results for Special Cases}
  In the case that peculiar properties of a Levin-type sequence
  transformation are used, more stringent theorems can often be proved
  as regards convergence acceleration using this particular
  transformation.

   In the case of the Levin transformation, Sidi proved the following
   theorem:

    \begin{theorem}\label{thSidi}\cite{Sidi90}
    \cite[Theorem 2.32]{BrezinskiRedivoZaglia91}
      If $s_n=s+\omega_n f_n$ where $f_n \sim \sum_{j=0}^{\infty}
      \beta_j/n^j$ with $\beta_0\ne 0$ and $\omega_n \sim \sum_{j=0}^{\infty}
      \delta_j/n^{j+a}$ with $a>0$, $\delta_0\ne 0$ for $n\to\infty$
      then, if $\beta_k\ne 0$
      \begin{equation}
        \L\di{n}{k} -s \sim \frac{\delta_0 \beta_k}{\binom{-a}{k}}
        \cdot n^{-a-k} \qquad (n\to\infty).
      \end{equation}
    \end{theorem}

   For the $W$ algorithm and the $d\di{}{1}$ transformation that may be
   regarded as direct generalizations of the Levin transformation, Sidi has
   obtained a large number of results. The interested reader is
   referred to the literature (See \cite{Sidi95,Sidi95caf} and
       references therein).

   Convergence results for the Levin transformation, the Drummond
   transformation and the Weniger transformations may be found in
   Section 13 of Weniger's report \cite{Weniger89}.

   Results for the \J\ transformation and in particular, for the \pJ\
   transformation are given in \cite{Homeier96aan,Homeier96Hab}. Here,
   we recall the following theorems:

   \begin{theorem}\label{thaltern}
Assume that the following holds:
\begin{description}
\item[(A-0)] The sequence $\Seq{s_n}$ has the (anti)limit $s$.
\item[(A-1a)] For every $n$, the elements of the sequence $\Seq{\omega_n}$ are strictly
alternating in sign and do not vanish.
\item[(A-1b)] For all $n$ and $k$, the elements of the sequence $\Seq{\delta\di{n}{k}}=\Seq{\Delta
r\di{n}{k}}$ are of the same sign and do not vanish.
\item[(A-2)] For all $n\in\Nset_0$ the ratio $(s_n-s)/\omega_n$ can be
expressed as a series of the form
\begin{equation}
\frac{s_n-s}{\omega_n} = c_0 + \sum_{j=1}^{\infty} c_j
\sum_{n>n_1>n_2>\cdots>n_j}
\delta\di{n_1}{0}\delta\di{n_2}{1}\cdots\delta\di{n_j}{j-1}
\end{equation}
with $c_0\ne 0$.
\end{description}
\noindent Then the following holds for
$s\di{n}{k}=\J\di{n}{k}(\Seq{s_n},\Seq{\omega_n},\Seq{\delta\di{n}{k}})$:
\begin{description}
\item[a)] The error $s\di{n}{k}-s$ satisfies
\begin{equation}
s\di{n}{k} -s =
\frac{b\di{n}{k}}
     {\nabla\di{n}{k-1}\nabla\di{n}{k-2}\cdots\nabla\di{n}{0}[1/\omega_n]}
\end{equation}
with
\begin{equation}
b\di{n}{k}= c_k + \sum_{j=k+1}^{\infty} c_j
\sum_{n>n_{k+1}>n_{k+2}>\cdots>n_{j}}
\delta\di{n_{k+1}}{k}\delta\di{n_{k+2}}{k+1}\cdots\delta\di{n_j}{j-1}\>.
\end{equation}
\item[b)] The error $s\di{n}{k}-s$ is bounded in magnitude according to
\begin{equation}
\vert s\di{n}{k}-s\vert \le
\vert
      \omega_n b\di{n}{k}
      \delta\di{n}{0}\delta\di{n}{1}\cdots\delta\di{n}{k-1}
\vert\>.
\end{equation}
\item[c)] For large $n$ the estimate
\begin{equation}
\frac{s\di{n}{k}-s}{s_n-s} = O(\delta\di{n}{0}\delta\di{n}{1}\cdots\delta\di{n}{k-1})
\end{equation}
holds if $b\di{n}{k}=O(1)$ and $(s_n-s)/\omega_n=O(1)$ as $n\to\infty$.
\end{description}
\end{theorem}

\begin{theorem}\label{thziel}
Define $s\di{n}{k}=\J\di{n}{k}
(\Seq{s_n},\Seq{\omega_n},\Seq{\delta\di{n}{k}})$ and $\omega\di{n}{k}=1/D\di{n}{k}$ where the $D\di{n}{k}$ are
defined as in Eq.\ (\ref{eqJtransRec}). Put  $e\di{n}{k} = 1- {\omega\di{n+1}{k}}/{\omega\di{n}{k}}$
and $b\di{n}{k} = {(s\di{n}{k}-s)}/{\omega\di{n}{k}}$. Assume that (A-0) of Theorem \ref{thaltern} holds
and that the following conditions are satisfied:
\begin{description}
      \item[(B-1)] Assume that
      \begin{equation}
      \lim_{n\to\infty}
      \frac{b\di{n}{k}}{b\di{n}{0}} = B_k
      \end{equation}
      exists and is finite.
      \item[(B-2)] Assume that
      \begin{equation}
      \Omega_k=\lim_{n\to\infty} \frac{\omega\di{n+1}{k}}{\omega\di{n}{k}} \ne 0
      \end{equation}
      and
      \begin{equation}
      F_k=\lim_{n\to\infty} \frac{\delta\di{n+1}{k}}{\delta\di{n}{k}}\ne 0\>,
      \end{equation}
      exist for all $k\in\Nset_0$. Hence the limits
      $\Phi_k=\lim_{n\to\infty}\Phi\di{n}{k}$ (compare Eq.\
      (\ref{eqPhink}))
      exist for all $k\in\Nset_0$.
\end{description}
Then, the following holds:
  \begin{description}
    \item[a)] If $\Omega_0 \not\in \Set{\Phi_0=1, \Phi_1, \dots,
      \Phi_{k-1}}$, then
      \begin{equation}
      \lim_{n\to\infty} \frac{s\di{n}{k}-s}{s_n-s}
         \left\{
           \prod_{l=0}^{k-1} \delta\di{n}{l}
         \right\} ^{-1}
      =
         B_k \frac{[\Omega_0]^k}
                 {\displaystyle \prod_{l=0}^{k-1} (\Phi_l - \Omega_0)}
      \end{equation}
      and, hence,
      \begin{equation}
      \frac{s\di{n}{k}-s}{s_n-s} = O(\delta\di{n}{0}\delta\di{n}{1}\cdots\delta\di{n}{k-1})
      \end{equation}
      holds in the limit $n\to\infty$.
    \item[b)] If $\Omega_l=1$ for $l\in\Set{0,1,2,\dots,k}$ then
      \begin{equation}
      \lim_{n\to\infty} \frac{s\di{n}{k}-s}{s_n-s}
         \left\{
           \prod_{l=0}^{k-1} \frac{\delta\di{n}{l}}{e\di{n}{l}}
         \right\} ^{-1}
      =
         B_k
      \end{equation}
      and, hence,
      \begin{equation}
      \frac{s\di{n}{k}-s}{s_n-s} =
       O\left(
                  \prod_{l=0}^{k-1} \frac{\delta\di{n}{l}}{e\di{n}{l}}
        \right)
      \end{equation}
      holds in the limit $n\to\infty$.
      \end{description}
\end{theorem}

This theorem has the following two corollaries for the \pJ\
transformation \cite{Homeier96aan}:

\begin{corollar}\label{thCo1}
Assume that the following holds:
\begin{description}
\item{(C-1)} Let $\beta>0\>,p\ge 1$ and $\delta\di{n}{k}=\forwarddiff
[(n+\beta+(p-1)k)^{-1}]$. Thus, we deal with the \pJ\
transformation and, hence, the equations $F_k=\lim_{n\to\infty}
\delta\di{n+1}{k}/\delta\di{n}{k}=1$ and $\Phi_k=1$ hold for all $k$.
\item{(C-2)} Assumptions (A-2) of Theorem \ref{thaltern} and (B-1) of
Theorem \ref{thziel} are satisfied for the particular choice (C-1) for
$\delta\di{n}{k}$.
\item{(C-3)} The limit $\Omega_0=\lim_{n\to\infty}
\omega_{n+1}/ \omega_{n}$ exists, and it satisfies
$\Omega_0\not\in \Set{0,1}$. Hence, all the limits
$\Omega_k=\lim_{n\to\infty}\omega\di{n+1}{k}/\omega\di{n}{k}$ exist for $k\in\Nset$
exist and satisfy $\Omega_k=\Omega_0$.
\end{description}
Then
the transformation
$s\di{n}{k}=\pJ\di{n}{k}(\beta,\Seq{s_n},\Seq{\omega_n})$ satisfies
      \begin{equation}
      \lim_{n\to\infty} \frac{s\di{n}{k}-s}{s_n-s}
         \left\{
           \prod_{l=0}^{k-1} \delta\di{n}{l}
         \right\} ^{-1}
      =
         B_k \left\{\frac{\Omega_0}
                 { 1 - \Omega_0}
             \right\}^k
      \end{equation}
      and, hence,
      \begin{equation}\label{eqorderlin}
      \frac{s\di{n}{k}-s}{s_n-s} = O\left( (n+\beta)^{-2k}
                              \right)
      \end{equation}
      holds in the limit $n\to\infty$.
\end{corollar}

Note that Corollary \ref{thCo1} can be applied in the case of linear convergence
because then $0<\abs{\Omega_0}<1$ holds.

Corollary \ref{thCo1} allows to conclude that
in the case of linear convergence, the \pJ\ transformations
should be superior to Wynn's epsilon algorithm \cite{Wynn56a}.
Consider for instance the case that
\begin{equation}
s_n \sim s + \lambda^n n^{\theta} \sum_{n=0}^{\infty}
c_j/n^{j}\>,
\qquad c_0\ne 0\>,\qquad n\to\infty
\end{equation}
is an asymptotic expansion of the sequence elements $s_n$.
Assuming $\lambda\ne 1$ and $\theta\not\in\{0,1,\dots,k-1\}$ it
follows that (\cite[p.\ 127]{Wimp81}; \cite[ p.\ 333,
Eq.\ (13.4-7)]{Weniger89})
\begin{equation}
\frac
 {\epsilon\di{2k}{n}-s}
 {s_n -s}
 =
 O\left(n^{-2k}\right)\>, \qquad n\to \infty \>.
\end{equation}
This is the same order of convergence acceleration as in Eq.\
(\ref{eqorderlin}). But it should be noted  that for the computation of
$\epsilon\di{2k}{n}$ the $2k+1$ sequence elements
$\{s_{n}, \dots, s_{n+2k}\}$ are required. But for the computation of
\pJ$\di{n}{k}$ only  the $k+1$ sequence elements
$\{s_{n}, \dots, s_{n+k}\}$ are required in the case of the $t$ and $u$
variants, and additionally $s_{n+k+1}$ in the case of the
$\tilde t$ variant. Again, this is similar to Levin-type accelerators
(\cite[p.\ 333]{Weniger89}).

The following corollary
\ref{thCo2} applies to the case of logarithmic convergence:

\begin{corollar}\label{thCo2}
Assume that the following holds:
\begin{description}
\item{(D-1)} Let $\beta>0\>,p\ge 1$ and $\delta\di{n}{k}=\Delta
[(n+\beta+(p-1)k)^{-1}]$. Thus, we deal with the \pJ\
transformation and, hence, the equations $F_k=\lim_{n\to\infty}
\delta\di{n+1}{k}/\delta\di{n}{k}=1$ and $\Phi_k=1$ hold for all $k$
\item{(D-2)} Assumptions (A-2) of Theorem \ref{thaltern} and (B-1) of
Theorem \ref{thziel} are satisfied for the particular choice (C-1) for
$\delta\di{n}{k}$.
\item{(D-3)} Some constants $a\di{l}{j}$, $j=1,2$, exist such
that
\begin{equation}
e\di{n}{l} = 1- \omega\di{n+1}{l}/\omega\di{n}{l}=
 \frac{a\di{l}{1}}{n+\beta}
+\frac{a\di{l}{2}}{(n+\beta)^2}
+O((n+\beta)^{-3})
\end{equation}
holds for $l=0$. This implies
that this equation, and hence, $\Omega_l=1$  holds for
$l\in\Set{0,1,2,\dots,k}$. Assume further that $a\di{l}{1}\ne 0$ for
$l\in\Set{0,1,2,\dots,k-1}$.
\end{description}
Then
the transformation $s\di{n}{k}=\pJ\di{n}{k}(\beta,\Seq{s_n},\Seq{\omega_n})$ satisfies
      \begin{equation}
      \lim_{n\to\infty} \frac{s\di{n}{k}-s}{s_n-s}
         \left\{
           \prod_{l=0}^{k-1} \frac{\delta\di{n}{l}}{e\di{n}{l}}
         \right\} ^{-1}
      =
         B_k
      \end{equation}
      and, hence,
      \begin{equation}
      \frac{s\di{n}{k}-s}{s_n-s} =
       O\left(
                (n+\beta)^{-k}
        \right)
      \end{equation}
      holds in the limit $n\to\infty$.
\end{corollar}

For convergence acceleration results regarding the \H\ and \I\
transformations, see \cite{Homeier92} and \cite{Homeier98aah}.

{\quad}\section{Stability Results for  Levin-Type
Transformations}\label{secStaRes}

  \subsection{General Results}
We remind the reader of the definition of the \emph{stability indices}
${\boldsymbol \Gamma}\di{n}{k}(\T)\ge 1$ as given in Eq.\ (\ref{eqGammank}). We
consider the sequence $\Seq{\omega_n}\in\Oset^{\Kset}$ as given.
We call the transformation \T\  \emph{stable along the path}
$\mathcal{P}=\{(n_\ell,k_\ell)\vert [ n_\ell>n_{\ell-1} \mbox{ and }
k_\ell\ge k_{\ell-1} ]\mbox{ or }
[ n_\ell \ge n_{\ell-1} \mbox{ and }
k_\ell> k_{\ell-1} ]
\}$ in the \T\ table if the limit of its stability index along the path $\mathcal{P}$
exists and is bounded, i.e., if
\begin{equation}
\lim_{\ell\to\infty} {\boldsymbol \Gamma}\di{n_\ell}{k_{\ell}}(\T)
= \lim_{\ell\to\infty}
\sum_{j=0}^{k_\ell} \vert\gamma\di{n_\ell,j}{k_\ell}(\vec \omega_n) \vert < \infty\>
\end{equation}
where the $\gamma\di{n,j}{k}(\vec \omega_n)$ are defined in Eq.\
(\ref{eqgammanjk}).
The transformation \T\ is called \emph{S-stable}, if it is stable along  all paths
$\mathcal{P}\di{}{k}=\{(n,k)\vert n=0,1,\dots\}$ for fixed $k$, i.e.
along all columns in the \T\ table.

The case of stability along diagonal paths is much more
difficult to treat analytically unless Theorem \ref{thstabgen}
applies. Up to now it seems that such diagonal stability issues
have only been analysed by Sidi for the case of the $d^{(1)}$
transformation (see \cite{Sidi95caf} and references therein).
We will treat only S-stability in the sequel.

The higher the stability index ${\boldsymbol \Gamma}(\T)$ is,
the smaller is the numerical stability of the transformation $\T$: If
$\epsilon_j$ is the numerical error of $s_{j}$,
\begin{equation}
\epsilon_j = s_j - \float{s_j}
\end{equation}
then the difference between the true value $\T\di{n}{k}$ and the
numerically computed approximation $\float{\T\di{n}{k}}$ may
be bounded according to
\begin{equation}
\left\vert \T\di{n}{k}-\float{\T\di{n}{k}}\right\vert
\le {\boldsymbol \Gamma}\di{n}{k}(\T) \> \left(\max_{j\in\{0,1,\dots,k\}}
\vert \epsilon_{n+j}\vert \right)\>.
\end{equation}
Compare also \cite{Sidi95caf}.

  \begin{theorem}\label{thstab}
     If the Levin-type sequence transformation $\T\di{}{k}$ has a limiting
     transformation $\Ringel\T\di{}{k}$ with characteristic polynomial
     $\Ringel\Pi\di{}{k}\in \Pset\di{}{k}$ for all $k\in \Nset$, and if
     $\Seq{\omega_n}\in\Oset^{\Kset}$ satisfies $\omega_{n+1}/\omega_n \sim \rho\ne 0$ for large $n$ with
     $\Ringel\Pi\di{}{k}(1/\rho) \ne 0$ for all $k\in \Nset$ then the
     transformation \T\ is S-stable. If additionally, the coefficients
     $\Ringel\lambda\di{j}{k}$ of the characteristic polynomial
     alternate in sign, i.e., if
     $\Ringel\lambda\di{j}{k}=(-1)^{j}\vert
     \Ringel\lambda\di{j}{k}\vert/\tau_k$ with $\vert\tau_k\vert=1$, then the limits
     $\Ringel\Gamma\di{}{k}(\T)=\lim_{n\to\infty} \Gamma\di{n}{k}(\T)$
     obey
\begin{equation}\label{eqGammabeh}
\Ringel{\boldsymbol \Gamma}\di{}{k}(\T)
=\tau_k\frac{ \Ringel\Pi\di{}{k}(-1/\vert\rho\vert)
}
{\vert\Ringel\Pi\di{}{k}(1/\rho)\vert}
\end{equation}
\end{theorem}
  \begin{proof}
   We have for fixed $k$
   \begin{equation}
\gamma\di{n,j}{k}(\vec \omega_n)
=
\lambda\di{n,j}{k} \frac{\omega_n}{\omega_{n+j}}
\left[\sum_{j'=0}^{k}\lambda\di{n,j'}{k}\frac{\omega_n}{\omega_{n+j'}}\right]^{-1}
\sim
 \Ringel\lambda\di{j}{k} \rho^{-j}
\left[\sum_{j'=0}^{k}\Ringel\lambda\di{j'}{k}\rho^{-j'}\right]^{-1}
=
\frac{ \Ringel\lambda\di{j}{k} \rho^{-j}}
{\Ringel\Pi\di{}{k}(1/\rho)}
\end{equation}
whence
\begin{equation}
\lim_{n\to\infty} {\boldsymbol \Gamma}\di{n}{k}(\T)
= \lim_{n\to\infty}
\sum_{j=0}^{k} \vert\gamma\di{n,j}{k}(\vec \omega_n) \vert
=
\frac{ \sum_{j=0}^{k} \vert\Ringel\lambda\di{j}{k}\vert\> \vert\rho\vert^{-j}
}
{\vert\Ringel\Pi\di{}{k}(1/\rho)\vert}
< \infty\>.
\end{equation}
If the $\Ringel\lambda\di{j}{k}$ alternate in sign, we obtain for these
limits
\begin{equation}
\Ringel{\boldsymbol \Gamma}\di{}{k}(\T)
=
\tau_k\frac{ \sum_{j=0}^{k} \Ringel\lambda\di{j}{k} (-\vert\rho\vert)^{-j}
}
{\vert\Ringel\Pi\di{}{k}(1/\rho)\vert}\>.
\end{equation}
This implies Eq.\ (\ref{eqGammabeh}).
  \end{proof}

\begin{corollar}
     Assume that the Levin-type sequence transformation $\T\di{}{k}$ has a limiting
     transformation $\Ringel\T\di{}{k}$ with characteristic polynomial
     $\Ringel\Pi\di{}{k}\in \Pset\di{}{k}$ and the coefficients
     $\Ringel\lambda\di{j}{k}$ of the characteristic polynomial
     alternate in sign, i.e., if
     $\Ringel\lambda\di{j}{k}=(-1)^{j}\vert
     \Ringel\lambda\di{j}{k}\vert/\tau_k$ with $\vert\tau_k\vert=1$
     for all $k\in \Nset$. The sequence
     $\Seq{\omega_n}\in\Oset^{\Kset}$ is assumed to be alternating and
     to satisfy $\omega_{n+1}/\omega_n
     \sim \rho<0 $ for large $n$.  Then the
     transformation \T\ is S-stable. Additionally the limits are
     $\Ringel{\boldsymbol \Gamma}\di{}{k}(\T)=1$.
\end{corollar}
\begin{proof}
Since
\begin{equation}
\sum_{j'=0}^{k}\frac{\lambda\di{n,j'}{k}}{\tau_k^{-1}}\frac{\omega_n}{\omega_{n+j'}}
\sim
\frac{\Ringel\Pi\di{}{k}(1/\rho)}{\tau_k^{-1}}
=
\sum_{j'=0}^{k}
  \frac{\Ringel\lambda\di{j'}{k}}
       {\tau_k^{-1}}
  (-1)^{j'}
  \vert\rho\vert^{-j'}
=\sum_{j'=0}^{k}\vert\Ringel\lambda\di{j'}{k}\vert
\>\vert\rho\vert^{-j'}\ge \vert\Ringel\lambda\di{k}{k}\vert /\vert \rho\vert^k >0
\end{equation}
$1/\rho$ cannot be a zero of $\Ringel\Pi\di{}{k}$. Then, Theorem
\ref{thstab} entail that \T\ is S-stable. Furthermore, Eq.\
(\ref{eqGammabeh}) is applicable and yields
$\Ringel{\boldsymbol \Gamma}\di{}{k}(\T)=1$.
\end{proof}

This result can be improved if all the coefficients
$\lambda\di{n,j}{k}$ are alternating:
\begin{theorem}\label{thstabgen}
     Assume that the Levin-type sequence transformation $\T\di{}{k}$ has a
characteristic polynomials
     $\Pi\di{n}{k}\in \Pset\di{}{k}$ with alternating  coefficients
     $\lambda\di{n,j}{k}$
     i.e.,
     $\lambda\di{n,j}{k}=(-1)^{j}\vert
     \lambda\di{n,j}{k}\vert/\tau_k$ with $\vert\tau_k\vert=1$ for all $n\in\Nset_0$ and $k\in \Nset$. The sequence
     $\Seq{\omega_n}\in\Oset^{\Kset}$ is assumed to be alternating and
     to satisfy $\omega_{n+1}/\omega_n
     <0 $ for all $n\in\Nset_0$.  Then we have
     ${\boldsymbol\Gamma}\di{n}{k}(\T)=1$. Hence, the
     transformation \T\ is stable along all paths for such remainder
     estimates.
\end{theorem}
  \begin{proof}
   We have for fixed $n$ and $k$
   \begin{equation}
\gamma\di{n,j}{k}(\vec \omega_n)
=
\frac{\lambda\di{n,j}{k} \frac{\omega_n}{\omega_{n+j}}}
{\sum_{j'=0}^{k}\lambda\di{n,j'}{k}\frac{\omega_n}{\omega_{n+j'}}}
=
\frac{\lambda\di{n,j}{k} \tau_k(-1)^j\vert{\omega_n}/{\omega_{n+j}}\vert}
{\sum_{j'=0}^{k}\lambda\di{n,j'}{k}\tau_k(-1)^{j'}
\left\vert\frac{\omega_n}{\omega_{n+j'}}\right\vert}
=
\frac{\vert\lambda\di{n,j}{k}\vert
\>\vert{\omega_n}/{\omega_{n+j}}\vert}
{\sum_{j'=0}^{k}\vert \lambda\di{n,j'}{k}\vert\>
\left\vert\frac{\omega_n}{\omega_{n+j'}}\right\vert}
\ge 0\>.
\end{equation}
Note that the denominators cannot vanish and are bounded from below by
$\vert\lambda\di{n,k}{k}\omega_{n}/\omega_{n+k}\vert >0$.
Hence, we have $\gamma\di{n,j}{k}(\vec
\omega_n)=\vert\gamma\di{n,j}{k}(\vec \omega_n)\vert$ and consequently,
${\boldsymbol\Gamma}\di{n}{k}(\T)=1$ since
$\sum_{j=0}^{k} \gamma\di{n,j}{k}(\vec
\omega_n)=1$ according to Eq.\ (\ref{eqgammanjk}).
\end{proof}

  \subsection{Results for Special Cases}
     Here, we collect some special results on the stability of
     various Levin-type sequence transformations that have been
     reported in \cite{Homeier98ots} and generalize some results
     of Sidi on the S-stability of the $d\di{}{1}$
     transformation.

\begin{theorem}\label{theorem2}
If the sequence $\omega_{n+1}/\omega_n$ possesses a limit according to
\begin{equation}
\lim_{n\to\infty} \omega_{n+1}/\omega_n = \rho \ne 0
\end{equation}
and if $\rho\not\in\{1,\Phi_1,\dots,\Phi_{k},\dots\}$ such that the limiting
transformation exists, the ${\mathcal{J}}$ transformation
is S-stable with the same limiting stability indices as the transformation
$\Ringel{\mathcal{J}}$, i.e., we have
\begin{equation}\label{A}
\lim_{n\to\infty} {\boldsymbol{\Gamma}}\di{n}{k}=
\frac{\sum_{j=0}^k \vert \lambda\di{j}{k}\rho^{k-j}\vert}
{\prod_{j'=0}^{k-1}
\vert\Phi_{j'}-\rho\vert}<\infty\>.
\end{equation}
If all $\Phi_k$ are positive then
\begin{equation}\label{B}
\lim_{n\to\infty} {\boldsymbol{\Gamma}}\di{n}{k}=
\prod_{j=0}^{k-1}
\frac{\Phi_{j}+\vert\rho\vert}
{\vert\Phi_{j}-\rho\vert}<\infty\>
\end{equation}
holds.
\end{theorem}

As corollaries, we get the following results
\begin{corollar}\label{cor1}
If the sequence $\omega_{n+1}/\omega_n$ possesses a limit according to
\begin{equation}
\lim_{n\to\infty} \omega_{n+1}/\omega_n = \rho \not\in \{0,1\}
\end{equation}
the ${\PJ{p}}$ transformation for $p>1$ and $\beta>0$
is S-stable and we have
\begin{equation}
\lim_{n\to\infty} {\boldsymbol{\Gamma}}\di{n}{k}=
\frac{\sum_{j=0}^k  {k \choose j} \vert \rho^{k-j}\vert}
{\vert 1-\rho\vert^k}=\frac{(1+\vert \rho\vert)^k}
{\vert 1-\rho\vert^k}<\infty\>.
\end{equation}
\end{corollar}
\begin{corollar}\label{cor2}
If the sequence $\omega_{n+1}/\omega_n$ possesses a limit according to
\begin{equation}
\lim_{n\to\infty} \omega_{n+1}/\omega_n = \rho \not\in \{0,1\}
\end{equation}
the Weniger $\mathcal{S}$ transformation  \cite[Sec. 8]{Weniger89} for
$\beta>0$
is S-stable and we have
\begin{equation}
\lim_{n\to\infty} {\boldsymbol{\Gamma}}\di{n}{k}(\mathcal{S})=
\frac{\sum_{j=0}^k  {k \choose j} \vert \rho^{k-j}\vert}
{\vert 1-\rho\vert^k}=\frac{(1+\vert \rho\vert)^k}
{\vert 1-\rho\vert^k}<\infty\>.
\end{equation}
\end{corollar}

\begin{corollar}\label{cor3}
If the sequence $\omega_{n+1}/\omega_n$ possesses a limit according to
\begin{equation}
\lim_{n\to\infty} \omega_{n+1}/\omega_n = \rho \not\in \{0,1\}
\end{equation}
the Levin $\mathcal{L}$ transformation \cite{Levin73,Weniger89}
is S-stable and we have
\begin{equation}
\lim_{n\to\infty} {\boldsymbol{\Gamma}}\di{n}{k}(\mathcal{L})=
\frac{\sum_{j=0}^k  {k \choose j} \vert \rho^{k-j}\vert}
{\vert 1-\rho\vert^k}=\frac{(1+\vert \rho\vert)^k}
{\vert 1-\rho\vert^k}<\infty\>.
\end{equation}
\end{corollar}

\begin{corollar}\label{cor4}
Assume that the elements of the sequence $\{t_n\}_{n\in\Nset}$ satisfy
$t_n\ne 0$ for all $n$ and $t_n\ne t_{n'}$ for all $n\ne n'$.
If the sequence $t_{n+1}/t_n$ possesses a limit
\begin{equation}
\lim_{n\to\infty} t_{n+1}/t_n = \tau \mbox{ with } 0<
\tau < 1
\end{equation}
and if the sequence $\omega_{n+1}/\omega_n$ possesses a limit according to
\begin{equation}
\lim_{n\to\infty} \omega_{n+1}/\omega_n = \rho \not \in
\{0,1,\tau^{-1},\dots,\tau^{-k},\dots\}
\end{equation}
then the generalized Richardson extrapolation process $\mathcal{R}$
introduced by Sidi \cite{Sidi82} that is identical to the $\mathcal{J}$
transformation with $\delta\di{n}{k}=t_{n}-t_{n+k+1}$ as shown in
\cite{Homeier94ahc}, i.e., the $W$ algorithm is S-stable and we have
\begin{equation}
\lim_{n\to\infty} {\boldsymbol{\Gamma}}\di{n}{k}(\mathcal{R})=
\frac{\sum_{j=0}^k \vert \widetilde{\lambda}\di{j}{k}\rho^{k-j}\vert}
{\prod_{j'=0}^{k-1}
\vert \tau^{-j'}-\rho\vert}=
\prod_{j'=0}^{k-1}\frac{
1 +\tau^{j'}\vert\rho\vert}{
\vert 1- \tau^{j'}\rho\vert}<\infty\>.
\end{equation}
Here,
\begin{equation}
\widetilde{\lambda}\di{j}{k}=(-1)^{k-j} \sum_{j_0+j_1+\dots+j_{k-1}=j,
\atop j_0\in\{0,1\},\dots,j_{k-1}\in\{0,1\}} \prod_{m=0}^{k-1}
(\tau)^{-m\,j_m}\>
\end{equation}
such that
\begin{equation}
\sum_{j=0}^k \widetilde{\lambda}\di{j}{k}\rho^{k-j} = \prod_{j=0}^{k-1}
(\tau^{-j}-\rho) =
\tau^{-k(k-1)/2}
\prod_{j=0}^{k-1}(1-\tau^{j}\rho)\>.
\end{equation}
\end{corollar}

Note that the preceeding corollary is essentially the same as a result of Sidi
\cite[Theorem 2.2]{Sidi95caf}
that now appears as a special case of the more general
Theorem \ref{theorem2} that applies to a much wider class of sequence
transformations. As noted above, Sidi  has also derived conditions under
which the $d\di{}{1}$ transformation is stable along the paths
$\mathcal{P}_n=\{(n,k)\vert k=0,1,\dots\}$ for fixed $n$. For
details and more references see \cite{Sidi95caf}. Analogous work for
the $\mathcal{J}$ transformation is in progress.

An efficient algorithm for the computation of the stability index of the
$\mathcal{J}$ transformation can be given in the case
$\delta\di{n}{k}>0$. Since the $\mathcal{J}$ transformation is
invariant under $\delta\di{n}{k}\to \alpha\di{}{k}\delta\di{n}{k}$ for
any $\alpha\di{}{k}\ne 0$ according to \cite[Theorem 4]{Homeier94ahc},
$\delta\di{n}{k}>0$ can always be achieved if for  given
$k$, all $\delta\di{n}{k}$ have the same sign. This is the case, for
instance, for the $\PJ{p}$ transformation
\cite{Homeier94ahc,Homeier96aan}.

\begin{theorem}
Define
\begin{equation}
F\di{n}{0} = (-1)^{n} \vert D\di{n}{0} \vert\>,\qquad
F\di{n}{k+1}=(F\di{n+1}{k}-F\di{n}{k})/ \delta\di{n}{k}\>
\end{equation}
and $\widehat F\di{n}{0}=F\di{n}{0}$, $\widehat
F\di{n}{k}=(\delta\di{n}{0}\cdots \delta\di{n}{k-1})\,F\di{n}{k}$.
If all $\delta\di{n}{k}>0$ then
\begin{enumerate}
\item
  $F\di{n}{k}=(-1)^{n+k}\vert F\di{n}{k}\vert$,
\item
  \label{itemii} $\lambda\di{n,j}{k}=(-1)^{j+k}\vert
  \lambda\di{n,j}{k}\vert$, and
\item
  \begin{equation}\boldsymbol{\Gamma}\di{n}{k}=\frac{\vert \widehat
  F\di{n}{k}\vert}{\vert \widehat D\di{n}{k}\vert}=\frac{\vert
  F\di{n}{k}\vert}{\vert D\di{n}{k}\vert}\>.
  \end{equation}
\end{enumerate}
\end{theorem}
This generalizes Sidi's method for the computation of stability
indices \cite{Sidi95caf} to a larger class of sequence
transformations.

{\quad}\section{Application of Levin-type Sequence
Transformations}\label{secApplic}

\subsection{Practical Guidelines}

Here, we address shortly the following questions:

\begin{description}
\item[When should one try to use sequence transformations?]
One can only hope for good convergence acceleration,
extrapolation, or summation results if a) the $s_n$ have some
asymptotic structure for large $n$ and are not erratic or random, b)
a sufficiently large number of decimal digits is available. Many
problems can be successfully tackled if 13-15 digits are
available but some require a much larger number of digits in
order to overcome some inevitable rounding errors, especially
for the acceleration of logarithmically convergent sequences.
The asymptotic information that is required for a successful
extrapolation is often hidden in the last digits of the problem
data.
\item[How should the transformations be applied?] The
recommended mode of application is that one computes the highest
possible order $k$ of the transformation from the data. In the
case of triangular recursive schemes like that of the \J\
transformation and the Levin transformation, this means that one
computes as transformed sequence $\Seq{\T\di{0}{n}}$. For
L-shaped recursive schemes as in the case of the \H, \I, and \K\
transformations, one usually computes as transformed sequence
$\Seq{\T\di{n-2\Ent{n/2}}{\Ent{n/2}}}$. The error $\epsilon$ of the current
estimate can usually be approximated {\it a posteriori} using
sums of magnitudes of differences of a few entries of the \T\
table, e.g.,
\begin{equation}
\epsilon \approx \vert \T\di{1}{n} - \T\di{0}{n}\vert +
\vert \T\di{0}{n-1} - \T\di{0}{n}\vert\>
\end{equation}
for transformations with triangular recursive schemes.
Such a simple approach works surprisingly well in practice.
The loss of decimal digits can be estimated computing
stability indices. An example is given below.
\item[What happens if one of the denominator vanishes?] The occurrence of
zeroes in the $D$ table for specific combinations of $n$ and $k$
is usually no problem since the recurrences for numerators and
denominators still work in this case. Thus, no special devices
are required to jump over such singular points in the \T\ table.
\item[Which transformation and which variant should be chosen?]
This depends on the type of convergence of the problem sequence.
For linearly convergent sequences, $t$, $\widetilde t$, $u$ and
$v$ variants of the Levin transformation, the \pJ\
transformation, especially the \PJ{2}\ transformation are
usually a good choice \cite{Homeier96aan} as long as one is not
too close to a singularity or to a logarithmically convergent
problem. Especially well-behaved is usually the application to
alternating series since then, the stability is very good as
discussed above. For the summation of alternating divergent
sequences and series, usually the $t$ and the $\widetilde t$
variants of the Levin transformation, the \PJ{2}\ and the
Weniger \S\ and \M\ transformations provide often surprisingly
accurate results. In the case of logarithmic convergence, $t$
and $\widetilde t$ variants become useless, and the order of
acceleration is dropping from $2k$ to $k$ when the
transformation is used columnwise. If a Kummer-related series is
available (cp.\ Section \ref{secseqser}), then $K$ and $lu$
variants leading to linear sequence transformations can be
efficient \cite{HomeierWeniger95}.
Similarly, linear variants can be based on some good asymptotic estimates ${}^{asy}\omega_n$,
that have to be obtained via a separate
analysis \cite{HomeierWeniger95}.
In the case of logarithmcic convergence, it pays to consider
special devices like using subsequences $\Seq{s_{\xi_n}}$ where
the $\xi_n$ grow exponentially like $\xi_n=\Ent{\sigma
\xi_{n-1}}+1$ like in the $d$ transformations. This choice can
be also used in combination with the \F\ transformation.
Alternatively, one can use some other transformations like the
condensation transformation
\cite{PelzlKing98,JentschuraMohrSoffWeniger99cav} or interpolation to
generate a linearly convergent sequence \cite{Homeier99tlt},
before applying an usually nonlinear sequence transformation. A
somewhat different approach is possible if one can obtain a few
terms $a_n$ with large $n$ easily \cite{Homeier99cao}.
\item[What to do near a singularity?] When extrapolating power
series or, more generally, sequences depending on certain
parameters, quite often extrapolation becomes difficult near the
singularities of the limit function. In the case of linear
convergence, one can often transform to a problem with a larger
distance to the singularity: If Eq.\ (\ref{eqbekaverf_rho})
holds, then the subsequence $\Seq{s_{\tau
n}}$ satisfies
\begin{equation}
        \lim_{n\to\infty}
        ({s_{\tau (n+1)}-s})/({s_{\tau n}-s}) = \rho^{\tau}
\end{equation}
This is a method of Sidi that has can, however, be applied to
large classes of sequence transformations \cite{Homeier98ots}.
\item[What to do for more complicated convergence type?]
Here, one should try to rewrite the problem sequence as a
sum of sequences with more simple convergence behavior. Then,
nonlinear sequence transformations are used to extrapolate each
of these simpler series, and to sum the extrapolation results to
obtain an estimate for the original problem. This
is for instance often possible for (generalized) Fourier series  where it leads to
  complex series that may asymptotically be regarded as power
  series. For details, the reader is referred to the literature
  \cite{BrezinskiRedivoZaglia91,Homeier93,%
  Homeier96Hab,TC-NA-97-1,TC-NA-97-3,%
  TC-NA-97-4,Homeier98oca,Homeier98aah,Sidi95}. If this approach
  is not possible one is forced to use more complicated sequence
  transformations like the $d\di{}{m}$ transformations or the
  (generalized) \H\ transformation. These more complicated
  sequence transformations, however, do require more numerical
  effort to achieve a desired accuracy.
\end{description}

\subsection{Numerical Examples}

\begin{table}
\caption{Comparison of the \F\ transformation and the W algorithm for
the series (\ref{eqexserFW})}
\label{tabFW}
\begin{tabular*}{\linewidth}{r@{\extracolsep{0.01 cm plus 1 fill}}rrrr}
\hline\hline
$n$ & $A_n$ & $B_n$ & $C_n$ & $D_n$ \\
\hline
 14 & 13.16 & 13.65 & 7.65 & 11.13\\
 16 & 15.46 & 15.51 & 9.43 & 12.77\\
 18 & 18.01 & 17.84 & 11.25 & 14.43\\
 20 & 21.18 & 20.39 & 13.10 & 16.12\\
 22 & 23.06 & 23.19 & 14.98 & 17.81\\
 24 & 25.31 & 26.35 & 16.89 & 19.53\\
 26 & 27.87 & 28.17 & 18.83 & 21.26\\
 28 & 30.83 & 30.59 & 20.78 & 23.00\\
 30 & 33.31 & 33.19 & 22.76 & 24.76\\
\hline
  $n$  & $E_n$ & $F_n$ & $G_n$ & $H_n$ \\
\hline
      14 & 14.07 & 13.18 & 9.75 & 10.47\\
      16 & 15.67 & 15.49 & 11.59 & 12.05\\
      18 & 17.94 & 18.02 & 13.46 & 13.66\\
      20 & 20.48 & 20.85 & 15.37 & 15.29\\
      22 & 23.51 & 23.61 & 17.30 & 16.95\\
      24 & 25.66 & 25.63 & 19.25 & 18.62\\
      26 & 27.89 & 28.06 & 21.23 & 20.31\\
      28 & 30.46 & 30.67 & 23.22 & 22.02\\
      29 & 31.82 & 32.20 & 24.23 & 22.89\\
      30 & 33.43 & 33.45 & 25.24 & 23.75\\
\hline\hline
\end{tabular*}\\
Plotted is the negative decadic logarithm of the relative error.\\
$A_n$: $\F\di{0}{n}(\Seq{S_n(z,a)},\Seq{(2+\ln(n+a))\forwarddiff
S_n(z,a)},\Seq{1+\ln(n+a)})$ \\
$B_n$: $W\di{0}{n}(\Seq{S_n(z,a)},\Seq{(2+\ln(n+a))\forwarddiff
S_n(z,a)},\Seq{1/(1+\ln(n+a))})$\\
$C_n$: $\F\di{0}{n}(\Seq{S_n(z,a)},\Seq{(n+1)\forwarddiff
S_n(z,a)},\Seq{1+n+a})$\\
$D_n$: $W\di{0}{n}(\Seq{S_n(z,a)},\Seq{(n+1)\forwarddiff
S_n(z,a)},\Seq{1/(1+n+a)})$\\
$E_n$: $\F\di{0}{n}(\Seq{S_n(z,a)},\Seq{\forwarddiff
S_n(z,a)},\Seq{1+\ln(n+a)})$ \\
$F_n$: $W\di{0}{n}(\Seq{S_n(z,a)},\Seq{\forwarddiff
S_n(z,a)},\Seq{1/(1+\ln(n+a))})$\\
$G_n$: $\F\di{0}{n}(\Seq{S_n(z,a)},\Seq{\forwarddiff
S_n(z,a)},\Seq{1+n+a})$\\
$H_n$: $W\di{0}{n}(\Seq{S_n(z,a)},\Seq{\forwarddiff
S_n(z,a)},\Seq{1/(1+n+a)})$\\

\end{table}

In Table \ref{tabFW}, we present results of the application of
certain variants of the \F\
transformation and the $W$ algorithm to the series
  \begin{equation} \label{eqexserFW}
    S(z,a) = 1+ \sum_{j=1}^{\infty} z^{j}\>\prod_{\ell=0}^{j-1}
    \frac{1}{\ln(a+\ell)}
  \end{equation}
  with partial sums
  \begin{equation}
    S_n(z,a) = 1+ \sum_{j=1}^{n} z^{j}\>\prod_{\ell=0}^{j-1}
    \frac{1}{\ln(a+\ell)}
  \end{equation}
  for $z=1.2$ and $a=1.01$. Since the terms $a_j$ satisfy
  $a_{j+1}/a_j=z/\ln(a+j)$, the ratio test reveals that $S(z,a)$
  converges for all $z$ and, hence, represents an analytic
  function. Nevertheless, only for $j\ge -a + \exp(\vert
  z\vert)$, the ratio of the terms becomes less than unity in
  absolute value. Hence, for larger $z$ the series converges
  rather slowly.

  It should be noted that for cases $C_n$ and $G_n$, the \F\
  transformation is identical to the Weniger transformation \S,
  i.e., to the \PJ{3}\ transformation,
  and for cases $C_n$ and $H_n$ the $W$ algorithm is identical
  to the Levin transformation. In the upper part of the table,
  we use u-type remainder estimates while in the lower part, we
  use $\widetilde t$ variants. It is seen that the choices
  $x_n=1+\ln(a+n)$ for the \F\ transformation and
  $t_n=1/(1+\ln(a+n))$ for the $W$ algorithm perform for both
  variants nearly identical (columns $A_n$, $B_n$, $E_n$ and
  $F_n$) and are superior to the
  choices $x_n=1+n+a$ and $t_n=1/(1+n+a)$, respectively, that
  correspond to the Weniger and the Levin transformation as
  noted above. For the latter two transformations, the Weniger
  ${}^{\widetilde t}\S$ transformation is slightly superior the
  ${}^{\widetilde t}\L$ transformation for this particular
  example (columns $G_n$ vs. $H_n$) while the situation is
  reversed for the u-type variants displayed in colums $C_n$ and
  $D_n$.


The next example is taken from \cite{Homeier98ots}, namely
the ``inflated Riemann $\zeta$ function'',
i.e., the series
\begin{equation}
\zeta(\epsilon,1,q) = \sum_{j=0}^{\infty} \frac{q^{j}}{(j+1)^{\epsilon}}
\end{equation}
that is a special case of the Lerch zeta function $\zeta(s,b,z)$ (cf.
\cite[p.\ 142, Eq.\ (6.9.7)]{Hansen75} and
\cite[Sec.\ 1.11]{ErdelyiMagnusOberhettingerTricomi53I}).
The partial sums are defined as
\begin{equation}
s_n = \sum_{j=0}^{n} \frac{q^{j}}{(j+1)^{\epsilon}}\>.
\end{equation}
The series converges linearly for $0<\vert q\vert < 1$ for any complex
$\epsilon$. In fact, we have in this case
$\rho=\lim_{n\to\infty}(s_{n+1}-s)/(s_n-s)=q$. We choose $q=0.95$ and
$\epsilon=-0.1+10i$. Note that for this value of $\epsilon$, there is a
singularity of $\zeta(\epsilon,1,q)$ at $q=1$ where the defining series
diverges since $\Re(\epsilon)<1$.

\begin{table}
\caption{Acceleration of $\zeta(-1/10+10i,1,95/100)$ with the $\J$
transformation}\label{tab1t1}
\tabcolsep=0.2cm
\begin{tabular}{@{}rrrrrrrrr@{}}
\hline\hline
n & \multicolumn{1}{c}{$A_n$} &\multicolumn{1}{c}{$B_n$} &
    \multicolumn{1}{c}{$C_n$} &\multicolumn{1}{c}{$D_n$} &
    \multicolumn{1}{c}{$E_n$} &\multicolumn{1}{c}{$F_n$} &
    \multicolumn{1}{c}{$G_n$} &\multicolumn{1}{c}{$H_n$} \\
\hline
  10 & 2.59e-05 & 3.46e+01 & 2.11e-05 & 4.67e+01 & 1.84e-05 & 4.14e+01 & 2.63e-05 & 3.90e+01 \\
  20 & 1.72e-05 & 6.45e+05 & 2.53e-05 & 5.53e+07 & 1.38e-04 & 3.40e+09 & 1.94e-05 & 2.47e+06 \\
  30 & 2.88e-05 & 3.52e+10 & 8.70e-06 & 2.31e+14 & 8.85e-05 & 1.22e+17 & 2.02e-05 & 6.03e+11 \\
  40 & 4.68e-06 & 1.85e+15 & 8.43e-08 & 1.27e+20 & 4.06e-06 & 2.78e+23 & 1.50e-06 & 9.27e+16 \\
  42 & 2.59e-06 & 1.46e+16 & 2.61e-08 & 1.51e+21 & 2.01e-06 & 4.70e+24 & 6.64e-07 & 8.37e+17 \\
  44 & 1.33e-06 & 1.10e+17 & 7.62e-09 & 1.76e+22 & 1.73e-06 & 7.85e+25 & 2.76e-07 & 7.24e+18 \\
  46 & 6.46e-07 & 8.00e+17 & 1.80e-09 & 2.02e+23 & 1.31e-05 & 1.30e+27 & 1.09e-07 & 6.08e+19 \\
  48 & 2.97e-07 & 5.62e+18 & 1.07e-08 & 2.29e+24 & 1.52e-04 & 2.12e+28 & 4.16e-08 & 5.00e+20 \\
  50 & 1.31e-07 & 3.86e+19 & 1.51e-07 & 2.56e+25 & 1.66e-03 & 3.43e+29 & 1.54e-08 & 4.05e+21 \\
\hline\hline
\end{tabular}\\
\begin{tabular*}{\linewidth}{@{}r@{\extracolsep\fill}rrrrrrrr@{}}
\multicolumn{9}{@{}l@{\extracolsep\fill}}{$A_n: \mbox{ Relative
error of }\PJ{{1}}\di{0}{n}(1,\Seq{s_{n}},\Seq{(n+1)(s_{n}-s_{n-1})})$}\\
\multicolumn{9}{@{}l@{\extracolsep\fill}}{$B_n: \mbox{ Stability
index of }\PJ{{1}}\di{0}{n}(1,\Seq{s_{n}},\Seq{(n+1)(s_{n}-s_{n-1})})$}\\
\multicolumn{9}{@{}l@{\extracolsep\fill}}{$C_n: \mbox{ Relative
error of }\PJ{{2}}\di{0}{n}(1,\Seq{s_{n}},\Seq{(n+1)(s_{n}-s_{n-1})})$}\\
\multicolumn{9}{@{}l@{\extracolsep\fill}}{$D_n: \mbox{ Stability
index of }\PJ{{2}}\di{0}{n}(1,\Seq{s_{n}},\Seq{(n+1)(s_{n}-s_{n-1})})$}\\
\multicolumn{9}{@{}l@{\extracolsep\fill}}{$E_n: \mbox{ Relative
error of }\PJ{{3}}\di{0}{n}(1,\Seq{s_{n}},\Seq{(n+1)(s_{n}-s_{n-1})})$}\\
\multicolumn{9}{@{}l@{\extracolsep\fill}}{$F_n: \mbox{ Stability
index of }\PJ{{3}}\di{0}{n}(1,\Seq{s_{n}},\Seq{(n+1)(s_{n}-s_{n-1})})$}\\
\multicolumn{9}{@{}l@{\extracolsep\fill}}{$G_n: \mbox{ Relative error of }
\J\di{0}{n}(\Seq{s_{n}},\Seq{(n+1)(s_{n}-s_{n-1})},
\Set{1/(n+1)-1/(n+k+2)})$}\\
\multicolumn{9}{@{}l@{\extracolsep\fill}}{$H_n: \mbox{ Stability index of }
\J\di{0}{n}(\Seq{s_{n}},\Seq{(n+1)(s_{n}-s_{n-1})},
\Set{1/(n+1)-1/(n+k+2)})$}\\
 & & & & & & & &
\end{tabular*}
\end{table}

\begin{table}
\caption{Acceleration of $\zeta(-1/10+10i,1,95/100)$ with the $\J$
transformation ($\tau=10$)}\label{tab1}
\tabcolsep=0.2cm
\begin{tabular}{@{}rrrrrrrrr@{}}
\hline\hline
n & \multicolumn{1}{c}{$A_n$} &\multicolumn{1}{c}{$B_n$} &
    \multicolumn{1}{c}{$C_n$} &\multicolumn{1}{c}{$D_n$} &
    \multicolumn{1}{c}{$E_n$} &\multicolumn{1}{c}{$F_n$} &
    \multicolumn{1}{c}{$G_n$} &\multicolumn{1}{c}{$H_n$} \\
\hline
  10 & 2.10e{-}05 & 2.08e{+}01 & 8.17e{-}06 & 3.89e{+}01 & 1.85e{-}05 & 5.10e{+}01 & 1.39e{-}05 & 2.52e{+}01 \\
  12 & 2.49e{-}06 & 8.69e{+}01 & 1.43e{-}07 & 3.03e{+}02 & 9.47e{-}06 & 8.98e{+}02 & 1.29e{-}06 & 1.26e{+}02 \\
  14 & 1.93e{-}07 & 3.11e{+}02 & 5.98e{-}09 & 1.46e{+}03 & 8.24e{-}07 & 4.24e{+}03 & 6.86e{-}08 & 5.08e{+}02 \\
  16 & 1.11e{-}08 & 9.82e{+}02 & 2.02e{-}11 & 6.02e{+}03 & 6.34e{-}08 & 2.09e{+}04 & 2.57e{-}09 & 1.77e{+}03 \\
  18 & 5.33e{-}10 & 2.87e{+}03 & 1.57e{-}12 & 2.29e{+}04 & 4.08e{-}09 & 9.52e{+}04 & 7.81e{-}11 & 5.66e{+}03 \\
  20 & 2.24e{-}11 & 7.96e{+}03 & 4.15e{-}14 & 8.26e{+}04 & 2.31e{-}10 & 4.12e{+}05 & 2.07e{-}12 & 1.73e{+}04 \\
  22 & 8.60e{-}13 & 2.14e{+}04 & 8.13e{-}16 & 2.89e{+}05 & 1.16e{-}11 & 1.73e{+}06 & 4.95e{-}14 & 5.08e{+}04 \\
  24 & 3.07e{-}14 & 5.61e{+}04 & 1.67e{-}17 & 9.87e{+}05 & 5.17e{-}13 & 7.07e{+}06 & 1.10e{-}15 & 1.46e{+}05 \\
  26 & 1.04e{-}15 & 1.45e{+}05 & 3.38e{-}19 & 3.31e{+}06 & 1.87e{-}14 & 2.84e{+}07 & 2.33e{-}17 & 4.14e{+}05 \\
  28 & 3.36e{-}17 & 3.69e{+}05 & 6.40e{-}21 & 1.10e{+}07 & 3.81e{-}16 & 1.13e{+}08 & 4.71e{-}19 & 1.16e{+}06 \\
  30 & 1.05e{-}18 & 9.30e{+}05 & 1.15e{-}22 & 3.59e{+}07 & 1.91e{-}17 & 4.43e{+}08 & 9.19e{-}21 & 3.19e{+}06 \\
\hline\hline
\end{tabular}\\
\begin{tabular*}{\linewidth}{@{}r@{\extracolsep\fill}rrrrrrrr@{}}
\multicolumn{9}{@{}l@{\extracolsep\fill}}{$A_n: \mbox{ Relative
error of }\PJ{{1}}\di{0}{n}(1,\Seq{s_{10\,n}},\Seq{(10\,n+1)(s_{10\,n}-s_{10\,n-1})})$}\\
\multicolumn{9}{@{}l@{\extracolsep\fill}}{$B_n: \mbox{ Stability
index of }\PJ{{1}}\di{0}{n}(1,\Seq{s_{10\,n}},\Seq{(10\,n+1)(s_{10\,n}-s_{10\,n-1})})$}\\
\multicolumn{9}{@{}l@{\extracolsep\fill}}{$C_n: \mbox{ Relative
error of }\PJ{{2}}\di{0}{n}(1,\Seq{s_{10\,n}},\Seq{(10\,n+1)(s_{10\,n}-s_{10\,n-1})})$}\\
\multicolumn{9}{@{}l@{\extracolsep\fill}}{$D_n: \mbox{ Stability
index of }\PJ{{2}}\di{0}{n}(1,\Seq{s_{10\,n}},\Seq{(10\,n+1)(s_{10\,n}-s_{10\,n-1})})$}\\
\multicolumn{9}{@{}l@{\extracolsep\fill}}{$E_n: \mbox{ Relative
error of }\PJ{{3}}\di{0}{n}(1,\Seq{s_{10\,n}},\Seq{(10\,n+1)(s_{10\,n}-s_{10\,n-1})})$}\\
\multicolumn{9}{@{}l@{\extracolsep\fill}}{$F_n: \mbox{ Stability
index of }\PJ{{3}}\di{0}{n}(1,\Seq{s_{10\,n}},\Seq{(10\,n+1)(s_{10\,n}-s_{10\,n-1})})$}\\
\multicolumn{9}{@{}l@{\extracolsep\fill}}{$G_n: \mbox{ Relative
error of }
$}\\\multicolumn{9}{@{}l@{\extracolsep\fill}}{$\phantom{G_n:\mbox{\ }}
\J\di{0}{n}(\Seq{s_{10\,n}},\Seq{(10\,n+1)(s_{10\,n}-s_{10\,n-1})},
\Set{1/(10\,n+10)-1/(10\,n+10\,k+10)})$}\\
\multicolumn{9}{@{}l@{\extracolsep\fill}}{$H_n: \mbox{ Stability
index of }
$}\\\multicolumn{9}{@{}l@{\extracolsep\fill}}{$\phantom{H_n:\mbox{\ }}
\J\di{0}{n}(\Seq{s_{10\,n}},\Seq{(10\,n+1)(s_{10\,n}-s_{10\,n-1})},
\Set{1/(n+1)-1/(n+k+2)})$}\\
 & & & & & & & &
\end{tabular*}
\end{table}

The results of applying $u$ variants of the $\PJ{p}$
transformation with $p=1,2,3$ and of the Levin transformation to the
sequence of partial sums is displayed in Table \ref{tab1t1}.
For each
of these four variants of the $\J$ transformation, we give the
relative error and the stability index.
The true value of the series
(that is used to compute the errors) was computed using a more accurate
method described below. It is seen that the $\PJ{2}$
transformation achieves the best results. The attainable accuracy for
this transformation is limited to about 9 decimal digits by the fact
that the stability index displayed in the column $D_n$ of Table
\ref{tab1t1}\ grows relatively fast. Note that for $n=46$,
the number of digits (as given by the negative decadic logarithm
of the relative error) and the decadic logarithm of
the stability index sum up to approximately 32 which corresponds
to the maximal number of decimal digits that could be achieved
in the run.
Since the stability index increases with $n$, indicating
decreasing stability, it is clear that for higher values of $n$
the accuracy will be lower.

The magnitude of the stability index is largely controlled by the value
of $\rho$, compare Corollary \ref{cor1}. If one can treat a related
sequence with a smaller value of $\rho$, the stability index will be
smaller and thus, the stability of the extrapolation will be greater.

Such a related sequence is given by putting $\check{s}_\ell=s_{\xi_\ell}$
for $\ell\in\Nset_0$, where the sequence $\xi_{\ell}$ is a monotonously
increasing sequence of nonnegative integers. In the case of linear
convergent sequences, the choice $\xi_{\ell}=\tau \ell$ with $\tau\in\Nset$
can be used as in the case of the $d\di{}{1}$ transformation. It is
easily seen that the new sequence also converges linearly with
$\rho=\lim_{n\to\infty}(\check{s}_{n+1}-s)/(\check{s}_n-s)=q^\tau$. For
$\tau>1$, both the effectiveness and the stability  of the various
transformations are increased as shown in Table \ref{tab1} for the case
$\tau=10$. Note that this value was chosen to display basic features
relevant to the stability analysis, and is not necessarily the optimal
value. As in Table \ref{tab1t1}, the
relative errors and the stability indices of some variants of the
$\J$ transformation are displayed. These are nothing but the
$\PJ{p}$ transformation for $p=1,2,3$ and the Levin
transformation as applied to the sequence $\Seq{\check{s}_n}$ with
remainder estimates $\omega_n=(n\tau+\beta) (s_{n\tau}-s_{n\tau -1})$
for $\beta=1$.
Since constant factors in the remainder estimates are irrelevant since
the $\J$ transformation is invariant under any scaling
$\omega_n\to\alpha\omega_n$ for $\alpha\ne 0$, the same results would
have been obtained for $\omega_n=(n+\beta/\tau) (s_{n\tau}-s_{n\tau
-1})$.

If the Levin transformation is applied to
the series with partial sums $\check{s}_n=s_{\tau\,n}$,
and if the remainder estimates $\omega_n =
(n+\beta/\tau)(s_{\tau\,n}-s_{(\tau\,n)-1}) $ are used, then one
obtains nothing but the $d^{(1)}$ transformation with
$\xi_\ell=\tau \ell$ for $\tau\in\Nset$.
\cite{Sidi95,Homeier98ots}

It is seen from Table \ref{tab1} that again the best accuracy is
obtained for the $\PJ{2}$ transformation. The $d\di{}{1}$
transformation is worse, but better than the $\PJ{p}$
transformations for $p=1$ and $p=3$. Note that the
stability indices are now much smaller and do not limit the achievable
accuracy for any of the transformations up to $n=30$. The true value of
the series was computed numerically by applying the
$\PJ{2}$
transformation to the further sequence $\Seq{s_{40n}}$ and using
64 decimal digits in the calculation. In this way, a
sufficiently accurate approximation was
obtained that was used to compute the relative errors in Tables
\ref{tab1t1} and \ref{tab1}. A comparison value was
computed using the
representation \cite[p.\ 29, Eq.\
(8)]{ErdelyiMagnusOberhettingerTricomi53I}
\begin{equation}
\zeta(s,1,q)= \frac{\Gamma(1-s)}{z} (\log 1/q)^{s-1} + z^{-1}
\sum_{j=0}^{\infty} \zeta(s-j) \frac{(\log q)^{j}}{j!}\>
\end{equation}
that holds for $\vert \log q\vert < 2\pi$ and $s\not\in\Nset$.
Here, $\zeta(z)$ denotes the Riemann zeta function. Both values
agreed to all relevant decimal digits.

\begin{table}
\caption{Stability indices for the $\PJ{2}$
transformation ($\tau=10$)}\label{tab2}
\tabcolsep=0.2cm
\begin{tabular}{@{}llllllll@{}}
\hline\hline
$n$
 & \multicolumn{1}{c}{$\boldsymbol{\Gamma}\di{n}{1}$}
 & \multicolumn{1}{c}{$\boldsymbol{\Gamma}\di{n}{2}$}
 & \multicolumn{1}{c}{$\boldsymbol{\Gamma}\di{n}{3}$}
 & \multicolumn{1}{c}{$\boldsymbol{\Gamma}\di{n}{4}$}
 & \multicolumn{1}{c}{$\boldsymbol{\Gamma}\di{n}{5}$}
 & \multicolumn{1}{c}{$\boldsymbol{\Gamma}\di{n}{6}$}
 & \multicolumn{1}{c}{$\boldsymbol{\Gamma}\di{n}{7}$}
    \\
\hline
  20 & 3.07 & 9.26 & 2.70$\,10^1$ & 7.55$\,10^1$ & 2.02$\,10^2$   & 5.20$\,10^2$ & 1.29$\,10^3$  \\
  30 & 3.54 & 1.19$\,10^1$ & 3.81$\,10^1$ & 1.16$\,10^2$ &   3.36$\,10^2$ & 9.36$\,10^2$ & 2.51$\,10^3$  \\
  40 & 3.75 & 1.33$\,10^1$ & 4.49$\,10^1$ & 1.44$\,10^2$ &   4.42$\,10^2$ & 1.30$\,10^3$ & 3.71$\,10^3$  \\
  41 & 3.77 & 1.34$\,10^1$ & 4.54$\,10^1$ & 1.46$\,10^2$ &   4.51$\,10^2$ & 1.34$\,10^3$ & 3.82$\,10^3$  \\
  42 & 3.78 & 1.35$\,10^1$ & 4.59$\,10^1$ & 1.49$\,10^2$ &   4.60$\,10^2$ & 1.37$\,10^3$ & 3.93$\,10^3$  \\
  43 & 3.79 & 1.36$\,10^1$ & 4.64$\,10^1$ & 1.51$\,10^2$ &   4.69$\,10^2$ & 1.40$\,10^3$ & 4.05$\,10^3$  \\
  44 & 3.80 & 1.37$\,10^1$ & 4.68$\,10^1$ & 1.53$\,10^2$ &   4.77$\,10^2$ & 1.43$\,10^3$ &   \\
  45 & 3.81 & 1.38$\,10^1$ & 4.73$\,10^1$ & 1.55$\,10^2$ &   4.85$\,10^2$ &  &   \\
  46 & 3.82 & 1.39$\,10^1$ & 4.77$\,10^1$ & 1.57$\,10^2$ &  &  &   \\
  47 & 3.83 & 1.39$\,10^1$ & 4.81$\,10^1$ &  &  &  &   \\
  48 & 3.84 & 1.40$\,10^1$ &  &  &  &  &  \\
  49 & 3.85 &  &  &  &  &  &  \\
\hline
 Extr.     & 4.01 & 1.59$\,10^1$ & 6.32$\,10^1$ & 2.52$\,10^2$ &
 1.00$\,10^3$ & 4.00$\,10^3$ & 1.59$\,10^4$ \\
\hline
 Cor.~\ref{cor1} & 3.98 & 1.59$\,10^1$ & 6.32$\,10^1$ & 2.52$\,10^2$ &
 1.00$\,10^3$ & 4.00$\,10^3$ & 1.59$\,10^4$ \\
\hline\hline
\end{tabular}
\end{table}

In Table \ref{tab2}, we display stability indices corresponding to
the acceleration of $\check{s}_n$ with the $\PJ{2}$
transformation columnwise, as obtainable by using the sequence elements
up to $\check{s}_{50}=s_{500}$. In the row labelled
\emph{Cor.\ \ref{cor1}}, we display
the limits of the $\boldsymbol{\Gamma}\di{n}{k}$ for large $n$, i.e.,
the quantities
\begin{equation}
\lim_{n\to\infty}\boldsymbol{\Gamma}\di{n}{k} =
\left(\frac{1+q^{\tau}}{1-q^{\tau}}\right)^k
\end{equation}
that are the limits according to Corollary \ref{cor1}. It is seen that
the values for finite $n$ are still relatively far off the limits. In
order to check numerically the validity of the corollary, we
extrapolated the values of all $\boldsymbol{\Gamma}\di{n}{k}$
for fixed $k$ with $n$ up to the
maximal $n$ for which there is an entry in the corresponding
column of Table \ref{tab2} using the $u$ variant of the
$\PJ{1}$ transformation. The results of the
extrapolation are displayed in the row labelled \emph{Extr} in Table
\ref{tab2} and coincide nearly perfectly with the values
expected according to Corollary \ref{cor1}.

As a final example, we consider the evaluation of the $F_m(z)$
functions that are used in quantum chemistry calculations via
the series representation
\begin{equation}\label{eqFm}
F_m(z) = \sum_{j=0}^{\infty} (-z)^j/j!(2m+2j+1)\>,
\end{equation}
with partial sums
\begin{equation}
s_n = \sum_{j=0}^{n} (-z)^j/j!(2m+2j+1)\>.
\end{equation}
In this case, for larger $z$, the convergence is rather slow
although the convergence finally is hyperlinear.
As a $K$ variant, one may use
\begin{equation}
{}^{K}\omega_n=\left(\sum_{j=0}^{n} (-z)^j/(j+1)!-(1-e^{-z})/z\right)\>.
\end{equation}
since $(1-e^{-z})/z$ is a Kummer related series.
The results for several variants in Table \ref{tabFm} show that
the $K$ variant is superior to $u$ and $t$ variants in this case.

\begin{table}
\caption{Extrapolation of series representation (\ref{eqFm}) of
the $F_m(z)$ function using the \PJ{2} transformation ($z=8\>,m=0$)
}\label{tabFm}
{
{\begin{tabular*}{\textwidth}{r@{\extracolsep{\fill}}rlll}
\hline\hline
$n$ & \multicolumn{1}{c}{$s_n$} &
\multicolumn{1}{c}{${}^u\omega_n$} &
\multicolumn{1}{c}{${}^{t}\omega_n$} &
\multicolumn{1}{c}{${}^{K}\omega_n$}\\
\hline
  5&     --13.3&        0.3120747&   0.3143352&   0.3132981\\
  6&       14.7&        0.3132882&   0.3131147&   0.3133070\\
  7&     --13.1&        0.3132779&   0.3133356&   0.3133087\\
  8&       11.4&        0.3133089&   0.3133054&   0.3133087\\
  9&      --8.0&        0.3133083&   0.3133090&   0.3133087\\
\hline\hline
\end{tabular*}
}}
\end{table}

Many further numerical examples are given in the literature
\cite{Weniger89,HomeierWeniger95,Homeier96aan,TC-NA-97-1,TC-NA-97-3,TC-NA-97-4,Homeier98aah}.

{\footnotesize

\begin{thebibliography}{100}

\bibitem{AbramowitzStegun70}
M.~Abramowitz and I.~Stegun.
\newblock {\em Handbook of Mathematical Functions}.
\newblock Dover Publications, Inc., New York, 1970.

\bibitem{Aitken26}
A.~C. {A}itken.
\newblock On {Bernoulli}'s numerical solution of algebraic equations.
\newblock {\em Proc. Roy. Soc. Edinburgh}, 46:289--305, 1926.

\bibitem{Baker75}
G.~A. {Baker, Jr.}
\newblock {\em Essentials of {P}ad{\'e} approximants}.
\newblock Academic {P}ress, New York, 1975.

\bibitem{BakerGravesMorris81a}
G.~A. {Baker, Jr.} and P.~Graves-Morris.
\newblock {\em {P}ad{\'e} approximants. {P}art {I}: {B}asic theory}.
\newblock Addison-Wesley, Reading, Mass., 1981.

\bibitem{BakerGravesMorris96}
G.~A. {Baker, Jr.} and P.~Graves-Morris.
\newblock {\em {P}ad{\'e} approximants}.
\newblock Cambridge U.P., Cambridge (GB), second edition, 1996.

\bibitem{BenderOrszag87}
C.~M. {B}ender and S.~A. Orszag.
\newblock {\em Advanced mathematical methods for scientists and engineers}.
\newblock McGraw-Hill, Singapore, 1987.

\bibitem{Brezinski77}
C.~Brezinski.
\newblock {\em Acc{\'e}l{\'e}ration de la convergence en analyse
  num{\'e}rique}.
\newblock Springer-Verlag, Berlin, 1977.

\bibitem{Brezinski78}
C.~Brezinski.
\newblock {\em Algorithmes d'acc{\'e}l{\'e}ration de la convergence --
  {\'E}tude {nu\-m{\'e}\-ri\-que}}.
\newblock {\'E}ditions Technip, {P}aris, 1978.

\bibitem{Brezinski80b}
C.~Brezinski.
\newblock A general extrapolation algorithm.
\newblock {\em Numer. Math.}, 35:175--180, 1980.

\bibitem{Brezinski80a}
C.~Brezinski.
\newblock {\em {P}ad{\'e}-type approximation and general orthogonal
  polynomials}.
\newblock Birk\-h{\"a}u\-ser, Basel, 1980.

\bibitem{Brezinski91c}
C.~Brezinski.
\newblock {\em A bibliography on continued fractions, {P}ad{\'e} approximation,
  extrapolation and related subjects}.
\newblock Prensas Universitarias de Zaragoza, Zaragoza, 1991.

\bibitem{Brezinski91b}
C.~Brezinski, editor.
\newblock {\em Continued fractions and {P}ad{\'e} approximants}.
\newblock North-Holland, Amsterdam, 1991.

\bibitem{BrezinskiMatos96}
C.~Brezinski and A.~C. Matos.
\newblock A derivation of extrapolation algorithms based on error estimates.
\newblock {\em J. Comput. Appl. Math.}, 66(1-2):5--26, 1996.

\bibitem{BrezinskiRedivoZaglia91}
C.~Brezinski and M.~{Redivo Zaglia}.
\newblock {\em Extrapolation methods. {T}heory and practice}.
\newblock North-Holland, Amsterdam, 1991.

\bibitem{BrezinskiRedivoZaglia93a}
C.~Brezinski and M.~{Redivo Zaglia}.
\newblock A general extrapolation algorithm revisited.
\newblock {\em Adv. Comput. Math.}, 2:461--477, 1994.

\bibitem{CioslowskiWeniger93}
J.~Cioslowski and E.~J. Weniger.
\newblock Bulk properties from finite cluster calculations. {VIII} {B}enchmark
  calculations on the efficiency of extrapolation methods for the {HF} and
  {MP}2 energies of polyacenes.
\newblock {\em J. Comput. Chem.}, 14:1468--1481, 1993.

\bibitem{CizekVinetteWeniger91}
J.~{\v C}{\'\i}{\v z}ek, F.~Vinette, and E.~J. Weniger.
\newblock Examples on the use of symbolic computation in physics and chemistry:
  Applications of the inner projection technique and of a new summation method
  for divergent series.
\newblock {\em Int. J. Quantum Chem. Symp.}, 25:209--223, 1991.

\bibitem{CizekVinetteWeniger93}
J.~{\v C}{\'\i}{\v z}ek, F.~Vinette, and E.~J. Weniger.
\newblock On the use of the symbolic language {M}aple in physics and chemistry:
  {S}everal examples.
\newblock In R.~A. de~Groot and J.~Nadrchal, editors, {\em Pro\-ceed\-ings of
  the Fourth International Conference on Computational Physics {P}HYSICS
  COMPUTING '92}, pages 31--44, Singapore, 1993. World Scientific.

\bibitem{Drummond72}
J.~E. Drummond.
\newblock A formula for accelerating the convergence of a general series.
\newblock {\em Bull. Austral. Math. Soc.}, 6:69--74, 1972.

\bibitem{ErdelyiMagnusOberhettingerTricomi53I}
A.~Erd{\'e}lyi, W.~Magnus, F.~Oberhettinger, and F.~G. Tricomi.
\newblock {\em Higher transcendental functions}, volume~I.
\newblock McGraw-Hill, New York, 1953.

\bibitem{FesslerFordSmith83a}
T.~Fessler, W.~F. Ford, and D.~A. Smith.
\newblock {HURRY:} an acceleration algorithm for scalar sequences and series.
\newblock {\em ACM Trans. Math. Software}, 9:346--354, 1983.

\bibitem{FordSidi87}
W.~F. Ford and A.~Sidi.
\newblock An algorithm for a generalization of the {R}ichardson extrapolation
  process.
\newblock {\em SIAM J. Numer. Anal.}, 24(5):1212--1232, 1987.

\bibitem{GermainBonne73}
B.~Germain-Bonne.
\newblock Transformations de suites.
\newblock {\em Rev. Fran\c{c}aise Automat. Informat. Rech. Operat.}, 7
  (R-1):84--90, 1973.

\bibitem{GravesMorris72}
P.~R. Graves-Morris, editor.
\newblock {\em {P}ad{\'e} approximants}.
\newblock {T}he Institute of {P}hysics, London, 1972.

\bibitem{GravesMorris73}
P.~R. Graves-Morris, editor.
\newblock {\em {P}ad{\'e} approximants and their applications}.
\newblock Academic {P}ress, London, 1973.

\bibitem{GravesMorrisSaffVarga84}
P.~R. Graves-Morris, E.~B. Saff, and R.~S. Varga, editors.
\newblock {\em Rational approximation and interpolation}.
\newblock Springer-Verlag, Berlin, 1984.

\bibitem{Grotendorst91}
J.~Grotendorst.
\newblock A {M}aple package for transforming series, sequences and functions.
\newblock {\em Comput. {P}hys. Commun.}, 67:325--342, 1991.

\bibitem{GrotendorstSteinborn86}
J.~Grotendorst and E.~O. Steinborn.
\newblock Use of nonlinear convergence accelerators for the efficient
  evaluation of {{GTO}} molecular integrals.
\newblock {\em J. Chem. {P}hys.}, 84:5617--5623, 1986.

\bibitem{GrotendorstWenigerSteinborn86}
J.~Grotendorst, E.~J. Weniger, and E.~O. Steinborn.
\newblock Efficient evaluation of infinite-series representations for overlap,
  two-center nuclear attraction and {C}oulomb integrals using nonlinear
  convergence accelerators.
\newblock {\em {P}hys. Rev. A}, 33:3706 -- 3726, 1986.

\bibitem{Hansen75}
E.~R. Hansen.
\newblock {\em A table of series and products}.
\newblock Prentice-Hall, Englewood-Cliffs, 1975.

\bibitem{HasegawaSidi96}
T.~Hasegawa and A.~Sidi.
\newblock An automatic integration procedure for infinite range integrals
  involving oscillatory kernels.
\newblock {\em Numer. Algor.}, 13:1--19, 1996.

\bibitem{Haavie79}
T.~H{\aa}vie.
\newblock Generalized {N}eville type extrapolation schemes.
\newblock {\em BIT}, 19:204--213, 1979.

\bibitem{Homeier90}
H.~H.~H. Homeier.
\newblock {\em {Integraltransformationsmethoden} und
  {Qua\-dra\-tur\-ver\-fah\-ren} f{\"u}r {Molek{\"u}lintegrale} mit
  {B}-{Funktionen}}, volume 121 of {\em Theorie und Forschung}.
\newblock S. Roderer Verlag, Regensburg, 1990.
\newblock {\AlsoDoctoraldissertation, Universit\"at Regensburg}.

\bibitem{Homeier92}
H.~H.~H. Homeier.
\newblock A {Levin}--type algorithm for accelerating the convergence of
  {Fourier} series.
\newblock {\em Numer. Algo.}, 3:245--254, 1992.

\bibitem{Homeier93}
H.~H.~H. Homeier.
\newblock Some applications of nonlinear convergence accelerators.
\newblock {\em Int. J. Quantum Chem.}, 45:545--562, 1993.

\bibitem{Homeier94ahc}
H.~H.~H. Homeier.
\newblock A hierarchically consistent, iterative sequence transformation.
\newblock {\em Numer. Algo.}, 8:47--81, 1994.

\bibitem{Homeier94nca}
H.~H.~H. Homeier.
\newblock Nonlinear convergence acceleration for orthogonal series.
\newblock In R.~Gruber and M.~Tomassini, editors, {\em Proceedings of the 6th
  {J}oint {EPS--APS} {I}nternational {C}onference on {P}hysics {C}omputing,
  {P}hysics {C}omputing '94}, pages 47--50. European Physical Society, Boite
  Postale 69, CH-1213 Petit-Lancy, Genf, Schweiz, 1994.

\bibitem{Homeier95}
H.~H.~H. Homeier.
\newblock Determinantal representations for the {$\mathcal J$} transformation.
\newblock {\em Numer. Math.}, 71(3):275--288, 1995.

\bibitem{Homeier96aan}
H.~H.~H. Homeier.
\newblock Analytical and numerical studies of the convergence behavior of the
  {$\mathcal J$} transformation.
\newblock {\em J. Comput. Appl. Math.}, 69:81--112, 1996.

\bibitem{Homeier96Hab}
H.~H.~H. Homeier.
\newblock {\em {Extrapolationsverfahren f\"ur Zahlen-, Vektor- und
  Matrizenfolgen und ihre Anwendung in der Theoretischen und Physikalischen
  Chemie}}.
\newblock Habilitation thesis, Universit\"at Regensburg, 1996.

\bibitem{TC-NA-97-3}
H.~H.~H. Homeier.
\newblock Extended complex series methods for the convergence acceleration of
  {Fourier} series.
\newblock Technical Report TC-NA-97-3, Institut f{\"u}r {P}hysikalische und
  {T}heoretische {C}hemie, {U}niversit{\"a}t {R}egensburg, {D}-93040
  {R}egensburg, 1997.

\bibitem{TC-NA-97-4}
H.~H.~H. Homeier.
\newblock On an extension of the complex series method for the convergence
  acceleration of orthogonal expansions.
\newblock Technical Report TC-NA-97-4, Institut f{\"u}r {P}hysikalische und
  {T}heoretische {C}hemie, {U}niversit{\"a}t {R}egensburg, {D}-93040
  {R}egensburg, 1997.

\bibitem{TC-NA-97-1}
H.~H.~H. Homeier.
\newblock On properties and the application of {Levin}-type sequence
  transformations for the convergence acceleration of {Fourier} series.
\newblock Technical Report TC-NA-97-1, Institut f{\"u}r {P}hysikalische und
  {T}heoretische {C}hemie, {U}niversit{\"a}t {R}egensburg, {D}-93040
  {R}egensburg, 1997.

\bibitem{Homeier98aah}
H.~H.~H. Homeier.
\newblock An asymptotically hierarchy-consistent iterative sequence
  transformation for convergence acceleration of {Fourier} series.
\newblock {\em Numer. Algo.}, 18:1--30, 1998.

\bibitem{Homeier98oca}
H.~H.~H. Homeier.
\newblock On convergence acceleration of multipolar and orthogonal expansions.
\newblock {\em Internet J. Chem.}, 1(Article 28):no pp. given, 1998.
\newblock Online computer file: URL: http://www.ijc.com/articles/1998v1/28/.
  Proceedings of the {4$^{th}$ Electronic Computational Chemistry Conference}.

\bibitem{Homeier98ots}
H.~H.~H. Homeier.
\newblock On the stability of the {$\mathcal{J}$} transformation.
\newblock {\em Numer. Algo.}, 17:223--239, 1998.

\bibitem{Homeier99cao}
H.~H.~H. Homeier.
\newblock Convergence acceleration of logarithmically convergent series
  avoiding summation.
\newblock {\em Appl. Math. Lett.}, 12:29--32, 1999.

\bibitem{Homeier99tlt}
H.~H.~H. Homeier.
\newblock Transforming logarithmic to linear convergence by interpolation.
\newblock {\em Appl. Math. Lett.}, 12:13--17, 1999.

\bibitem{TC-PC-95-1}
H.~H.~H. Homeier and B.~Dick.
\newblock Zur {B}erechnung der {L}inienform spektraler {L}{\"o}cher, ({E}ngl.:
  {O}n the computation of the line shape of spectral holes).
\newblock Technical Report TC-PC-95-1, Institut f{\"u}r {P}hysikalische und
  {T}heoretische {C}hemie, {U}niversit{\"a}t {R}egensburg, {D}-93040
  {R}egensburg, 1995.
\newblock Poster CP 6.15, 59. Physikertagung Berlin 1995, Abstract:
  Verhandlungen der Deutschen Physikalischen Gesellschaft, Reihe VI, Band 30,
  1815, (Physik-Verlag GmbH, D-69469 Weinheim, 1995).

\bibitem{HomeierWeniger95}
H.~H.~H. Homeier and E.~J. Weniger.
\newblock On remainder estimates for {Levin}-type sequence transformations.
\newblock {\em Comput. Phys. Commun.}, 92:1--10, 1995.

\bibitem{JentschuraMohrSoffWeniger99cav}
U.~Jentschura, P.~J. Mohr, G.~Soff, and E.~J. Weniger.
\newblock Convergence acceleration via combined nonlinear-condensation
  transformations.
\newblock {\em Comput. Phys. Commun.}, 116:28 -- 54, 1999.

\bibitem{Khovanskii63}
A.~N. Khovanskii.
\newblock {\em {T}he application of continued fractions and their
  generalizations to problems in approximation theory}.
\newblock Noordhoff, Groningen, 1963.

\bibitem{Levin73}
D.~{L}evin.
\newblock Development of non-linear transformations for improving convergence
  of sequences.
\newblock {\em Int. J. Comput. Math. B}, 3:371--388, 1973.

\bibitem{LevinSidi81}
D.~Levin and A.~Sidi.
\newblock Two new classes of nonlinear transformations for accelerating the
  convergence of infinite integrals and series.
\newblock {\em Appl. Math. Comput.}, 9:175--215, 1981.

\bibitem{Longman81}
I.~M. Longman.
\newblock Difficulties of convergence acceleration.
\newblock In M.~G. de~Bruin and H.~van Rossum, editors, {\em Pad{\'e}
  approximation and its applications Amsterdam 1980}, pages 273--289.
  Springer-Verlag, Berlin, 1981.

\bibitem{LorentzenWaadeland92}
L.~Lorentzen and H.~Waadeland.
\newblock {\em Continued fractions with applications}.
\newblock North-Holland, Amsterdam, 1992.

\bibitem{LucasStone95}
S.~K. Lucas and H.~A. Stone.
\newblock Evaluating infinite integrals involving {Bessel} functions of
  arbitrary order.
\newblock {\em J. Comput. Appl. Math.}, 64:217--231, 1995.

\bibitem{MagnusOberhettingerSoni66}
W.~Magnus, F.~Oberhettinger, and R.~P. Soni.
\newblock {\em Formulas and theorems for the special functions of mathematical
  physics}.
\newblock Springer-Verlag, New York, 1966.

\bibitem{Matos97}
A.~C. Matos.
\newblock Linear difference operators and acceleration methods.
\newblock Publication ANO-370, Laboratoire d'Analyse {Num\'erique} et
  d'Optimisation, {Universit\'e} des Sciences et Technologies de Lille, France,
  submitted, 1997.

\bibitem{Michalski98}
K.~A. Michalski.
\newblock Extrapolation methods for {Sommerfeld} integral tails.
\newblock {\em {IEEE} Transactions on Antennas and Propagation},
  46(10):1405--1418, 1998.

\bibitem{Mosig89}
J.~R. Mosig.
\newblock Integral equation technique.
\newblock In T.~Itoh, editor, {\em Numerical Techniques for Microwave and
  Millimeter--Wave Passive Structures}, pages 133--213. Wiley, New York, 1989.

\bibitem{NikishinSorokin91}
E.~M. Nikishin and V.~N. Sorokin.
\newblock {\em Rational approximations and orthogonality}.
\newblock American Mathematical Society, Providence, Rhode Island, 1991.

\bibitem{Oleksy96}
C.~Oleksy.
\newblock A convergence acceleration method of {Fourier} series.
\newblock {\em Comput. Phys. Commun.}, 96:17--26, 1996.

\bibitem{Overholt65}
K.~J. Overholt.
\newblock Extended {A}itken acceleration.
\newblock {\em BIT}, 5:122--132, 1965.

\bibitem{PelzlKing98}
P.~J. Pelzl and F.~W. King.
\newblock Convergence accelerator approach for the high-precision evaluation of
  three-electron correlated integrals.
\newblock {\em Phys. Rev. E}, 57(6):7268--7273, 1998.

\bibitem{PetrushevPopov87}
P.~P. Petrushev and V.~A. Popov.
\newblock {\em Rational approximation of real functions}.
\newblock Cambridge U. P., Cambridge, 1987.

\bibitem{Ross87}
B.~Ross.
\newblock {\em Methods of summation}.
\newblock Descartes Press, Koriyama, 1987.

\bibitem{SaffVarga77}
E.~B. Saff and R.~S. Varga, editors.
\newblock {\em {Pad{\'e}} and rational approximation}.
\newblock Academic Press, New York, 1977.

\bibitem{Schumaker81}
L.~Schumaker.
\newblock {\em Spline functions: {Basic} theory}.
\newblock Wiley, New York, 1981.

\bibitem{Shanks55}
D.~Shanks.
\newblock Non-linear transformations of divergent and slowly convergent
  sequences.
\newblock {\em J. Math. and Phys. (Cambridge, Mass.)}, 34:1--42, 1955.

\bibitem{Sidi79}
A.~Sidi.
\newblock Convergence properties of some nonlinear sequence transformations.
\newblock {\em Math. Comput.}, 33:315--326, 1979.

\bibitem{Sidi79b}
A.~Sidi.
\newblock Some properties of a generalization of the {R}ichardson extrapolation
  process.
\newblock {\em J. Inst. Math. Appl.}, 24:327--346, 1979.

\bibitem{Sidi82}
A.~Sidi.
\newblock An algorithm for a special case of a generalization of the
  {R}ichardson extrapolation process.
\newblock {\em Numer. Math.}, 38:299--307, 1982.

\bibitem{Sidi88a}
A.~Sidi.
\newblock Generalization of {R}ichardson extrapolation with application to
  numerical integration.
\newblock In H.~Bra{\ss} and G.~H{\"a}mmerlin, editors, {\em Numerical
  integration}, volume III, pages 237--250. Birkh{\"a}user, Basel, 1988.

\bibitem{Sidi88c}
A.~Sidi.
\newblock A user-friendly extrapolation method for oscillatory infinite
  integrals.
\newblock {\em Math. Comp.}, 51:249--266, 1988.

\bibitem{Sidi90}
A.~Sidi.
\newblock On a generalization of the {R}ichardson extrapolation process.
\newblock {\em Numer. Math.}, 47:365--377, 1990.

\bibitem{Sidi95}
A.~Sidi.
\newblock Acceleration of convergence of (generalized) {F}ourier series by the
  {$d$}-transformation.
\newblock {\em Ann. Numer. Math.}, 2:381--406, 1995.

\bibitem{Sidi95caf}
A.~Sidi.
\newblock Convergence analysis for a generalized {Richardson} extrapolation
  process with an application to the $d^{(1)}$ transformation on convergent and
  divergent logarithmic sequences.
\newblock {\em Math. Comput.}, 64(212):1627--1657, 1995.

\bibitem{SmithFord79}
D.~A. Smith and W.~F. Ford.
\newblock Acceleration of linear and logarithmic convergence.
\newblock {\em SIAM J. Numer. Anal.}, 16:223--240, 1979.

\bibitem{SmithFord82}
D.~A. Smith and W.~F. Ford.
\newblock Numerical comparisons of nonlinear convergence accelerators.
\newblock {\em Math. Comput.}, 38(158):481--499, 1982.

\bibitem{SteinbornHomeierFernandezRicoEmaLopezRamirez99}
E.~O. Steinborn, H.~H.~H. Homeier, J.~{Fern\'andez Rico}, I.~Ema, R.~L{\'o}pez,
  and G.~Ram{\'\i}rez.
\newblock An improved program for molecular calculations with {{\em B}}
  functions.
\newblock {\em J. Mol. Struct. (Theochem)}, 490:201--217, 1999.

\bibitem{SteinbornWeniger90}
E.~O. Steinborn and E.~J. Weniger.
\newblock Sequence transformations for the efficient evaluation of infinite
  series representations of some molecular integrals with exponentially
  decaying basis functions.
\newblock {\em J. Mol. Struct. ({T}heochem)}, 210:71--78, 1990.

\bibitem{Wall73}
H.~S. Wall.
\newblock {\em Analytic theory of continued fractions}.
\newblock Chelsea, New York, 1973.

\bibitem{Weniger89}
E.~J. Weniger.
\newblock Nonlinear sequence transformations for the acceleration of
  convergence and the summation of divergent series.
\newblock {\em Comput. Phys. Rep.}, 10:189--371, 1989.

\bibitem{Weniger90}
E.~J. Weniger.
\newblock On the summation of some divergent hypergeometric series and related
  perturbation expansions.
\newblock {\em J. Comput. Appl. Math.}, 32:291--300, 1990.

\bibitem{Weniger91}
E.~J. Weniger.
\newblock On the derivation of iterated sequence transformations for the
  acceleration of convergence and the summation of divergent series.
\newblock {\em Comput. Phys. Commun.}, 64:19--45, 1991.

\bibitem{Weniger92}
E.~J. Weniger.
\newblock Interpolation between sequence transformations.
\newblock {\em Numer. Algor.}, 3:477--486, 1992.

\bibitem{Weniger94}
E.~J. Weniger.
\newblock {\em {Verallgemeinerte Summationsprozesse als numerische Hilfsmittel
  f\"ur quantenmechanische und quantenchemische Rechnungen}}.
\newblock Habilitationsschrift, Universit{\"a}t Regensburg, 1994.

\bibitem{Weniger96cotw}
E.~J. Weniger.
\newblock Computation of the {Whittaker} function of the second kind by summing
  its divergent asymptotic series with the help of nonlinear sequence
  transformations.
\newblock {\em Computers in Physics}, 10(5):496--503, 1996.

\bibitem{Weniger96cots}
E.~J. Weniger.
\newblock Construction of the strong coupling expansion for the ground state
  energy of the quartic, sextic, and octic anharmonic oscillator via a
  renormalized strong coupling expansion.
\newblock {\em Phys. Rev. Lett.}, 77(14):2859--2862, 1996.

\bibitem{Weniger96acr}
E.~J. Weniger.
\newblock A convergent renormalized strong coupling perturbation expansion for
  the ground state energy of the quartic, sextic, and octic anharmonic
  oscillator.
\newblock {\em Ann. Phys.}, 246(1):133--165, 1996.

\bibitem{Weniger96nste}
E.~J. Weniger.
\newblock Erratum: Nonlinear sequence transformations: {A} computational tool
  for quantum mechanical and quantum chemical calculations.
\newblock {\em Int. J. Quantum Chem.}, 58:319--321, 1996.

\bibitem{Weniger96nst}
E.~J. Weniger.
\newblock Nonlinear sequence transformations: {A} computational tool for
  quantum mechanical and quantum chemical calculations.
\newblock {\em Int. J. Quantum Chem.}, 57:265--280, 1996.

\bibitem{WenigerCizek90}
E.~J. Weniger and J.~{\v C}{\'\i}{\v z}ek.
\newblock Rational approximations for the modified {B}essel function of the
  second kind.
\newblock {\em Comput. Phys. Commun.}, 59:471--493, 1990.

\bibitem{WenigerCizekVinette91}
E.~J. Weniger, J.~{\v C}{\'\i}{\v z}ek, and F.~Vinette.
\newblock Very accurate summation for the infinite coupling limit of the
  perturbation series expansions of anharmonic oscillators.
\newblock {\em Phys. Lett. A}, 156:169--174, 1991.

\bibitem{WenigerCizekVinette93}
E.~J. Weniger, J.~{\v C}{\'\i}{\v z}ek, and F.~Vinette.
\newblock {T}he summation of the ordinary and renormalized perturbation series
  for the ground state energy of the quartic, sextic and octic anharmonic
  oscillators using nonlinear sequence transformations.
\newblock {\em J. Math. Phys.}, 34:571--609, 1993.

\bibitem{WenigerLiegener90}
E.~J. Weniger and C.-M. Liegener.
\newblock Extrapolation of finite cluster and crystal-orbital calculations on
  trans-polyacetylene.
\newblock {\em Int. J. Quantum Chem.}, 38:55--74, 1990.

\bibitem{WenigerSteinborn87}
E.~J. Weniger and E.~O. Steinborn.
\newblock Comment on ``molecular overlap integrals with exponential-type
  integrals''.
\newblock {\em J. Chem. Phys.}, 87:3709--3711, 1987.

\bibitem{WenigerSteinborn88}
E.~J. Weniger and E.~O. Steinborn.
\newblock Overlap integrals of {B} functions. {A} numerical study of infinite
  series representations and integral representations.
\newblock {\em {T}heor. Chim. Acta}, 73:323--336, 1988.

\bibitem{WenigerSteinborn89a}
E.~J. Weniger and E.~O. Steinborn.
\newblock Nonlinear sequence transformations for the efficient evaluation of
  auxiliary functions for {GTO} molecular integrals.
\newblock In M.~Defranceschi and J.~Delhalle, editors, {\em Numerical
  determination of the electronic structure of atoms, diatomic and polyatomic
  molecules}, pages 341--346, Dordrecht, 1989. Kluwer.

\bibitem{WernerBuenger84}
H.~Werner and H.~J. B{\"u}nger, editors.
\newblock {\em {Pad{\'e}} approximations and its applications, Bad Honnef
  1983}.
\newblock Springer-Verlag, Berlin, 1984.

\bibitem{Wimp81}
J.~Wimp.
\newblock {\em Sequence transformations and their applications}.
\newblock Academic Press, New York, 1981.

\bibitem{Wuytack79a}
L.~{W}uytack, editor.
\newblock {\em {Pad{\'e}} approximations and its applications}.
\newblock Springer-Verlag, Berlin, 1979.

\bibitem{Wynn56a}
P.~Wynn.
\newblock On a device for computing the {$e_m (S_n)$} transformation.
\newblock {\em Math. Tables Aids Comput.}, 10:91--96, 1956.

\end{thebibliography}
 \newcommand{\homeier}[1]{}

}

\begin{appendix}

{\quad}\section{Stieltjes Series and Functions}\label{appstieltjes}
A Stieltjes series is a formal expansion
\begin{equation}\label{eqbekaverf_stielt}
f(z) = \sum_{j=0}^{\infty} (-1)^j \, \mu_j \, z^j
\>
\end{equation}
with partial sums
\begin{equation}\label{eq39}
f_n(z) = \sum_{j=0}^{n} (-1)^j \, \mu_j \, z^j
\>.
\end{equation}
The coefficients $\mu_n$ are the moments of an uniquely given
positive measure $\psi(t)$  that has infinitely
many different values on $ 0 \le t < \infty$ \cite[p.\ 159]{BakerGravesMorris81a}:
\begin{equation}\label{eq40}
\mu_n = \int_{0}^{\infty} t^n \d \psi(t), \qquad n\in \Nset_0\>.
\end{equation}
Formally, the Stieltjes series can be identified with a Stieltjes
integral
\begin{equation}\label{eq41}
f(z) = \int_{0}^{\infty} \frac{\d \psi(t)}{1+zt}, \qquad
\abs{\mathrm{arg}(z)}<\pi\>.
\end{equation}
If such an integral exists for a
function $f$ then the function is called a
Stieltjes function. For every
Stieltjes function there exist a unique asymptotical Stieltjes series
(\ref{eqbekaverf_stielt}),
uniformly in every sector $\abs{\mathrm{arg}(z)}<\theta$ for all
$\theta<\pi$. For any Stieltjes series, however,  several
different corresponding  Stieltjes functions may exist. To ensure
uniqueness, additional criteria are necessary \cite[Sec. 4.3]{Weniger94}.

In the context of convergence acceleration and summation of divergent
series, it is important that for given $z$ the tails $f(z)-f_n(z)$ of a
Stieltjes series are bounded in absolute value by the
by the next term of the series,
\begin{equation}\label{eq42}
\vert f(z) -f_n(z) \vert \le  \mu_{n+1} \, z^{n+1}
\>{}\qquad z\ge 0\>.
\end{equation}
Hence, for Stieltjes series the remainder estimates may be chosen as
\begin{equation}\label{eq43}
\omega_n = (-1)^{n+1} \, \mu_{n+1} \, z^{n+1}
\>{}
\end{equation}
This corresponds to $\omega_n=\Delta f_n(z)$, i.e., to a $\widetilde t$
variant.

\section{Derivation of the Recursive Scheme (\ref{eqFtransRec})}\label{appB}
  We show that for the divided difference operator
  $\square\di{n}{k}=\square\di{n}{k}[\Seq{x_n}] $ the identity
  \begin{equation} \label{eqB1}
    \square\di{n}{k+1} ((x)_{\ell+1}g(x)) =
    \frac{(x_{n+k+1}+\ell)\square\di{n+1}{k} ((x)_{\ell}g(x))
         -(x_{n}+\ell)\square\di{n}{k} ((x)_{\ell}g(x))}
         {x_{n+k+1}-x_{n}}
  \end{equation}
  holds. The proof is based on the Leibniz formula for divided
  differences (see, e.g., \cite[p. 50]{Schumaker81}) that yields upon
  use of $(x)_{\ell+1}=(x+\ell)(x)_{\ell}$ and $\square\di{n}{k}
  (x)=x_n\delta_{k,0}+\delta_{k,1}$
  \begin{equation} \label{eqB2}
     \begin{array}{>{\displaystyle}r@{}>{\displaystyle}l}
        \square\di{n}{k+1} ((x)_{\ell+1}g(x))     & {}=      \ell
        \square\di{n}{k+1} ((x)_{\ell}g(x)) +
        \sum_{j=0}^{k+1}
           \square\di{n}{j}(x) \square\di{n+j}{k+1-j} ((x)_{\ell}g(x))\\
            & {}= (x_n+\ell) \square\di{n}{k+1} ((x)_{\ell}g(x)) +
                    \square\di{n+1}{k} ((x)_{\ell}g(x))
     \end{array}
  \end{equation}
  Using the recursion relation of the divided differences, one obtains
  \begin{equation} \label{eqB3}
    \square\di{n}{k+1} ((x)_{\ell+1}g(x)) =
    (x_n+\ell) \frac{\square\di{n+1}{k} ((x)_{\ell}g(x))
                      -\square\di{n}{k} ((x)_{\ell}g(x))}
                      {x_{n+k+1}-x_{n}}
      + \square\di{n+1}{k} ((x)_{\ell}g(x))\>.
  \end{equation}
  Simple algebra then yields Eq.\ (\ref{eqB1}).

  Comparison with Eq.\ (\ref{eqFtrans}) shows that using the
  interpolation conditions $g_n=g(x_n)=s_n/\omega_n$ and $\ell=k-1$ in
  Eq.\ (\ref{eqB1}) yields the recursion for the numerators in Eq.\
  (\ref{eqFtransRec}), while the recursion for the denominators in Eq.\
  (\ref{eqFtransRec}) follows for $\ell=k-1$ and using the interpolation
  conditions $g_n=g(x_n)=1/\omega_n$. In each case, the initial
  conditions follow directly from Eq.\ (\ref{eqFtrans}) in combination
  with the definition of the divided
  difference operator: For $k=0$, we use $(a)_{-1}=1/(a-1)$ and obtain
  $\square\di{n}{k} (x_n)_{k-1} g_n= (x_n)_{-1} g_n=g_n/(x_n-1)$.

\section{Two Lemmata}\label{applemma}

\begin{lemma}\label{thlemma}
  Define
  \begin{equation}
  A= \sum_{j=0}^{k} \Ringel\lambda\di{j}{k}
  \frac{\zeta^{n+j}}{(n+j)^{r+1}}\>
  \end{equation}
  where $\zeta$ is a zero of of multiplicity $m$ of
  $\Ringel\Pi\di{}{k}(z) = \sum_{j=0}^{k} \Ringel\lambda\di{j}{k} z^j$.
  Then
  \begin{equation}
  A \sim \zeta^{n+m} \binom{r+m}{r} \frac{(-1)^{m}}{n^{r+m+1}}
  \frac{d^m\Ringel\Pi\di{}{k}}{dx^m}(\zeta)
\>,\qquad (n\to\infty)\>.
  \end{equation}
\end{lemma}
\begin{proof}
Use
\begin{equation}
\frac{1}{a^{r+1}} = \frac{1}{r!} \int_{0}^{\infty} \exp(-a t) t^r
\,dt\>,\qquad a>0\>
\end{equation}
to obtain
\begin{equation}
A=
\frac{1}{r!}
  \int_{0}^{\infty} \sum_{j=0}^{k}
  \Ringel\lambda\di{j}{k}{\zeta^{n+j}}\exp(-(n+j) t) t^r
\,dt=
\frac{\zeta^n}{r!}
  \int_{0}^{\infty} \exp(-nt) \Ringel\Pi\di{}{k}(\zeta\exp(-t)) \, t^r
\,dt\>.
\end{equation}
Taylor expansion of the polynomial yields due to the zero at $\zeta$
\begin{equation}
\Ringel\Pi\di{}{k} (\zeta\exp(-t)) = \frac{(-\zeta)^m}{m!}
\left.\frac{d^m\Ringel\Pi\di{}{k}(x)}{dx^m}\right\vert_{x=\zeta}
t^m (1+O(t)) \>.
\end{equation}
Invoking Watson's lemma \cite[p.\ 263ff]{BenderOrszag87} completes the proof.
\end{proof}

\begin{lemma}\label{thlemma1}
Assume that assumption (C-3') of Theorem \ref{thconver2} holds. Further
assume $\lambda\di{n,j}{k}\to\Ringel\lambda\di{j}{k}$ for $n\to\infty$.
Then, Eq.\ (\ref{eqdenlim}) holds.
\end{lemma}
\begin{proof}
       We have
       \begin{equation}
         \frac{\omega_{n+j}}{\omega_n} \sim \rho^{j} \exp\left(\epsilon_n
         \sum_{t=0}^{j-1} \frac{\epsilon_{n+t}}{\epsilon_n}\right) \sim
         \rho^{j} \exp(j\epsilon_n)
       \end{equation}
       for large $n$. Hence,
       \begin{equation}
         \sum_{j=0}^{k}\lambda\di{n,j}{k} \frac{\omega_n}{\omega_{n+j}}
         \sim
         \sum_{j=0}^{k}\Ringel\lambda\di{j}{k}
         (\rho\exp(\epsilon_n))^{-j} =
         \Ringel\Pi\di{}{k}(1/\rho+\delta_n)
       \end{equation}
       Since the characteristic polynomial $\Ringel\Pi\di{}{k}(z)$ has
       a zero of order $\mu$ at $z=1/\rho$ according to the
       assumptions, Eq.\ (\ref{eqdenlim}) follows using Taylor
       expansion.
\end{proof}

\end{appendix}

\end{document}